\documentclass[%
 aip,
 amsmath,amssymb,
 reprint,%
]{revtex4-1}

\usepackage{graphicx}
\usepackage{dcolumn}
\usepackage{bm}

\usepackage{comment}
\usepackage{tikz}
\usepackage[export]{adjustbox}

\definecolor{rr1}{HTML}{880E4F}

\usepackage[utf8]{inputenc}
\usepackage[T1]{fontenc}
\usepackage{mathptmx}
\usepackage{etoolbox}

\usepackage[normalem]{ulem}

\makeatletter
\def\@email#1#2{%
 \endgroup
 \patchcmd{\titleblock@produce}
  {\frontmatter@RRAPformat}
  {\frontmatter@RRAPformat{\produce@RRAP{*#1\href{mailto:#2}{#2}}}\frontmatter@RRAPformat}
  {}{}
}%
\makeatother
\begin{document}

\preprint{AIP/123-QED}

\title[A dynamical study of Hilda asteroids in the Circular and Elliptic RTBP]{A dynamical study of Hilda asteroids in the Circular and Elliptic RTBP}

\author{Àngel Jorba}
\affiliation{ 
Departament de Matemàtiques i Informàtica, Universitat
  de Barcelona, Gran Via de les Corts Catalanes 585, 08007 Barcelona,
  Spain
}
\affiliation{Centre de Recerca Matem\`atica (CRM), Edifici C, Campus UAB, 08193 Bellaterra, Barcelona, Spain\hfill
}
\author{Begoña Nicolás}%

\affiliation{ 
Departamento de Matemática Aplicada, Universidade de Santiago de Compostela, Rua Lope Gómez De Marzoa, s/n, 15705, Santiago de Compostela (A Coru\~na), Spain
}%

\author{Óscar Rodríguez}
  \email{oscar.rodriguez@upc.edu.}
\affiliation{%
Departament de Matem\`atiques, Universitat Polit\`ecnica de Catalunya, Av. Diagonal 647,
08028 Barcelona, Spain}%

\date{\today}

\begin{abstract}
\vspace{4mm}

The Hilda group is a set of asteroids whose mean motion is in a 3:2 orbital resonance with Jupiter. In this paper we use the planar Circular Restricted Three-Body Problem (CRTBP) as a dynamical model and we show that there exists a family of stable periodic orbits that are surrounded by islands of quasi-periodic motions.
We have computed the frequencies of these quasi-periodic motions and we have shown how
the Hilda family fits inside these islands.
We have compared these results with the ones obtained using the Elliptic Restricted Three-Body Problem and they are similar, showing the suitability of the CRTBP model. It turns out that, to decide if a given asteroid belongs to the Hilda class, it is much better to look at its frequencies in the planar CRTBP rather than to use two-body orbital elements as it is commonly done today.
\end{abstract}

\maketitle

\begin{quotation}



This work is devoted to the analysis of the motion of a particular group of asteroids, from a dynamical systems point of view.
Typically asteroids are classified by means of a simple model that only accounts for the solar gravity. 
By including the effect of Jupiter we provide a more accurate way of classifying this group, based on the
geometrical invariant structures to which these asteroids belong.
\end{quotation}

\section{Introduction}
The Hilda asteroids constitute a group of more than 5000 asteroids located beyond the main asteroid belt of our Solar System, but still within Jupiter's orbit. This group is composed by two collisional families; the Hilda family and the Schubart family. Hilda family is named after (153) Hilda asteroid, discovered on 2 November 1875 by Johann Palisa. And Schubart family owes its name to (1911) Schubart asteroid, discovered on 25 October 1973, by Paul Wild. Meanwhile the diameter of (153) Hilda asteroid is of around 179 km, the one of (1911) Schubart asteroid of around 70 km. Both of them are within the largest ones in the group.

Hilda-type asteroids are mainly known by two characteristics. The first one is that they have mean motion in a 3:2 orbital resonance with Jupiter. And the second one is due the shape of their orbits, since they seem to successively approach three Lagrangian Points, $L_3$, $L_4$ and $L_5$, of the Sun-Jupiter system. In the neighbourhood of $L_4$ and $L_5$ of the Sun-Jupiter system another well-known group of asteroids is found, the Trojans. Unlike Trojans, that owe small orbital eccentricities and some have large inclinations, the eccentricities of Hilda asteroids are typically larger while their inclinations are not.

In spite of the considerable bibliography on this topic from astronomical approach, highlighting the work of J. Schubart \cite{Schubart1968,Schubart82}, little is known about the dynamical behaviour of the Hilda group of asteroids \cite{Morbidelli93,bookFM}. Therefore our aim in this work is to analyse Hilda's behaviour from a dynamical systems approach, by studying their orbits within Sun-Jupiter Circular Restricted Three-Body Problem (CRTBP) and Elliptical Restricted Three-Body Problem (ERTBP), both in the planar case, looking for the invariant objects that are considered to be responsible for their motion. One of the reasons for 
studying these two models is to analyse the level of importance of Jupiter eccentricity in this particular application.

Asteroids are usually classified in groups based on their (osculating) orbital elements
at a given time. This works quite well when the motion of the asteroid is close to Keplerian,
in the sense that their orbital elements vary very slowly and that they can be used to describe
its motion for a long time. However, the motion of some asteroids in the Solar system
is not so well described by a two-body problem in the sense that their orbital elements
have significant variations in time. Therefore, they can be classified into different
categories depending on the epoch. An example of this situation is given by the
Hilda group of asteroids, whose motion is strongly influenced by Jupiter.
The usual definition of this group is that they are the set of asteroids with
semi-major axis between $3.7$ and $4.2$ AU, eccentricity below $0.3$ and inclination lower
than $20$ degrees. It is clear that asteroids with orbital elements close to the boundary
of the group can have orbital elements crossing back and forth this boundary.

An improvement to classify asteroids is to use a better model than the two-body problem
to define the different families. For instance, the well-known
Circular Restricted Three-Body Problem (Circular RTBP or CRTBP, for short)
includes the main effect of Jupiter but, on the other hand, the model is not integrable so we
cannot easily identify the motion of the asteroid by a set of numbers,
like the orbital elements do in the two-body model.
In this paper we will discuss how to use the planar Circular RTBP to characterise
the members of the Hilda group of asteroids. As we will see in Section~\ref{sec:periodic}, there exists a
family of linearly stable periodic orbits that can be parameterised by the Jacobi constant
and acts as the ``skeleton'' of the Hilda group: these periodic orbits are surrounded
by a Cantor family of invariant tori~\cite{JorbaV97b} that, selecting a Jacobi constant
and using a suitable (two-dimensional) Poincar\'e section, are seen as an ``island''
surrounded by a chaotic sea. So, a more natural classification is to identify
the island that corresponds to the Hilda motion, so that the Hilda asteroids are
those inside that island. This is what is done in Section~\ref{sec:CRTBP}.

It is clear that a classification based on the dynamics of the circular RTBP
should be closer to the real motion rather than a classification based on a two-body model.
However, as the CRTBP is not the real system it is natural to ask how this
classification is affected by perturbations. For the concrete
case of the Hilda group the first relevant effect is due to the eccentricity
of Jupiter, which is a periodic perturbation with a frequency that is nearly resonant
with one of the frequencies (the mean motion) of these asteroids.
Therefore, a natural improvement to the planar CRTBP is to use the Elliptic RTBP (ERTBP).
For this reason, Section~\ref{sec:ERTBP} includes a similar study but now using
the planar ERTBP. The eccentricity of Jupiter appears in the model as a time-periodic
perturbation, with frequency 1, of the circular RTBP.

The paper is organised as follows. Sections~\ref{sec:circular} and~\ref{sec:elliptic}
summarises some basic information on the Circular and Elliptic RTBP, and
Section~\ref{sec:selec} contains more information on the Hilda group.
The CRTBP and ERTBP are introduced in three dimensions because we have to
transform the coordinates of the asteroids from the three-dimensional JPL reference frame (equatorial, J2000)
to our models, what is explained in Section~\ref{sec:CC}.
Section~\ref{sec:CRTBP} contains the dynamical study for the Circular RTBP and
Section~\ref{sec:ERTBP} the study for the Elliptic RTBP.
Conclusions and further work are summarised in Section~\ref{sec:Conclu}.
Finally, to facilitate the reading, an appendix on frequency analysis has
also been included.

\subsection{The Circular Restricted Three-Body Problem}\label{sec:circular}

The Circular Restricted Three-Body Problem (CRTBP) describes the motion of a particle of negligible mass under the gravitational field of two punctual massive bodies $P_1$ and $P_2$, called primaries, with masses $m_1$ and $m_2$, being $m_1\geq m_2$, that revolve in circular motion around their centre of mass.

It is usual to describe this problem in a rotating reference frame where the two primaries $P_1$ and $P_2$ are fixed and their motion is contained in the $xy$-plane. In this coordinate system, the period of revolution of the primaries is normalised to $2\pi$, the unit of mass is taken as the sum of the masses of the primaries, $m_1 + m_2$, and the length unit is adimensionalised using the distance between the primaries, that in this model is a constant quantity. With this, $P_1$ is placed at $(\mu,0,0)$ and $P_2$ at $(\mu-1,0,0)$, being $\mu=m_2/(m_1+m_2)$ the mass parameter.

In this configuration of the problem, five equilibrium points are found, called Lagrangian points; the collinear points $L_1$, $L_2$ and $L_3$ that are unstable, and the triangular points, $L_4$ and $L_5$ that are linearly stable for values of the mass parameter below the critical Routh mass
$\mu_R=(1-\sqrt{69}/9)/2 \approx0.03852$. Also, the following symmetries are found:
\begin{equation}
\begin{split}
(t,x,y,z,\dot{x},\dot{y},\dot{z}) &\rightarrow (-t,x,-y,-z,-\dot{x},\dot{y},\dot{z}),\\[1.ex]
(t,x,y,z,,\dot{x},\dot{y},\dot{z}) &\rightarrow (-t,x,-y,z,-\dot{x},\dot{y},-\dot{z}),
\end{split}
\label{eq:sym}
\end{equation}
where the dot ($~\dot{ }~$) denotes the derivative with respect to the adimensional time. Symmetries above mean that if the set of coordinates on the left of these expressions is a particular solution of the system, also is the set of coordinates on the right.

The CRTBP is an autonomous model with three degrees of freedom that presents a first integral, defined by
\begin{equation}\label{eq:CJ}
    C = 2\widetilde{\Omega}(x,y,z) - \left(\dot{x}^2 + \dot{y}^2 + \dot{z}^2\right),
\end{equation}
and known as Jacobi constant. Here, $\widetilde{\Omega}=\widetilde{\Omega}(x,y,z)$ is given by
\begin{equation}
\label{eq:Omega}
    \widetilde{\Omega} = \frac{1}{2}(x^2+y^2) + \frac{1-\mu}{r_1} + \frac{\mu}{r_2} + \frac{1}{2}\mu(1-\mu),
\end{equation}
where $r_1$ and $r_2$ are the distances to the primaries, that is
$r_1~=~\sqrt{(x-\mu)^2+y^2+z^2}$ and $r_2~=~\sqrt{(x-\mu+1)^2+y^2+z^2}$.
Note that the term $\frac{1}{2}\mu(1-\mu)$ included in Equation~\eqref{eq:Omega} ensures a value equal to $3$, exactly, for the Jacobi constant of the triangular points $L_4$ and $L_5$.
For a fixed value of the Jacobi constant $C$, the admissible regions of motion, known as Hill regions, are defined by
\[
    \mathcal{R}(C) = \left\{(x,y,z)\in\mathbb{R}^3\,|\,2\widetilde{\Omega}(x,y,z)\geq C\right\}.
\]

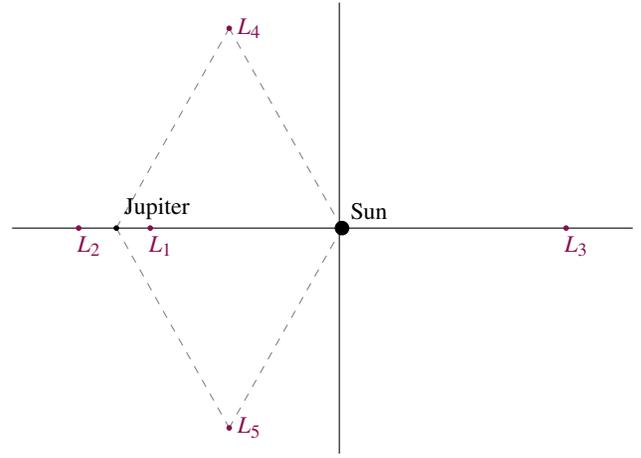
\begin{figure}
  \begin{center}
\begin{tikzpicture}[domain=0:2, scale = 3.0]
    \draw[black] (-1.45,0) -- (1.3,0); 
    \draw[black] (0,-1) -- (0,1); 
    \draw[gray, dashed] (-0.988,0)--(-0.488,0.886);
    \draw[gray, dashed] (0.012,0)--(-0.488,0.886);
    \draw[gray, dashed] (-0.988,0)--(-0.488,-0.886);
    \draw[gray, dashed] (0.012,0)--(-0.488,-0.886);
    \filldraw[rr1] (-0.488,0.886) circle (0.01) node[right]{$L_4$};
    \filldraw[rr1] (-0.488,-0.886) circle (0.01) node[right]{$L_5$};
    \filldraw[rr1] (-0.8377,0) circle (0.01) node[below]{\hspace{0.2 cm} $L_1$};
    \filldraw[rr1] (-1.1551,0) circle (0.01) node[below]{\hspace{0.2 cm} $L_2$};
    \filldraw[rr1] (1.005,0) circle (0.01) node[below]{ \hspace{0.2 cm} $L_3$};
    \filldraw[black] (-0.988,0) circle (0.01) node[above right]{Jupiter};
    \filldraw[black] (0.012,0) circle (0.03)  node[above right]{Sun};
\end{tikzpicture}
  \end{center}
    \caption{Sun-Jupiter system and its five equilibrium points.}
    \label{fig:ERTBP}
\end{figure}

In our particular case, we are interested in studying the Sun-Jupiter system. Therefore, $P_1$ corresponds to the Sun and $P_2$ to Jupiter, with mass parameter of value $\mu\approx9.53881\times 10^{-4}$. The primaries in the described configuration would be seen as in Figure~\ref{fig:ERTBP}, where the five Lagrangian points are also included. Notice that, triangular points are linearly stable in the Sun-Jupiter system since its mass parameter is below $\mu_R$.

\subsection{The Elliptic Restricted Three-Body Problem}\label{sec:elliptic}

The Elliptic Restricted Three-Body Problem (ERTBP) is a classical modification of the Circular RTBP, that accounts for the eccentricity of the elliptic orbits in which the primaries move.
Following the procedure detailed, for example in \cite{Szebehely1967}, we take adimensional units as in the CRTBP, but now the distance between the primaries, $r$, is not constant. It is given by
\begin{equation}\label{eq:r_ellip}
    r = \frac{a(1-e^2)}{1+e\cos f},
\end{equation}
where $a$ and $e$ are the semi-major axis and the eccentricity of the motion of either of the primaries around the other, respectively, and $f$ is the true anomaly. Therefore, the distance between the primaries in the ERTBP varies periodically with $f\in[0,2\pi)$.

Besides the normalisation of units, the motion of the infinitesimal particle in the ERTBP is typically described in a non-uniformly rotating and pulsating coordinate system, in which the primaries $P_1$ and $P_2$ are again fixed at positions $(\mu,0,0)$ and $(\mu-1,0,0)$, respectively, and their motion is contained in the $xy$-plane. In this way, the Elliptic RTBP can be seen as a $2\pi$ time-periodic perturbation of the Circular RTBP. However, due to the pulsating coordinates, the Lagrangian points continue being fixed points in the Elliptic RTBP configuration, as represented in  Figure~\ref{fig:ERTBP} for our case of interest, the Sun-Jupiter system, for which $e\approx0.04869$.

The rotating reference frame is defined using the angular velocity $\dot{f}$,
\begin{equation}\label{eq:df}
    \frac{df}{dt} = \frac{k\sqrt{m_1+m_2}}{a^{3/2}(1-e^2)^{3/2}}\left(1+e\cos f\right)^2,
\end{equation}
where $t$ is the dimensional time and $k$ is the Gauss' constant. Note that this expression comes from the angular momentum conservation.
Through some analysis, the equations of motion of the infinitesimal particle in the rotating and pulsating coordinate system can be written as
\begin{equation}\label{eq:ERTBP}
    \left\{
    \begin{aligned}
        x''  -2y' & = \frac{1}{1+e\cos f}\tilde{\Omega}_x,\\
        y''  +2x' & = \frac{1}{1+e\cos f}\tilde{\Omega}_y,\\
        z'' & = \frac{1}{1+e\cos f}\tilde{\Omega}_z - \frac{e\cos f}{1 + e\cos f}z,
    \end{aligned}
    \right.
\end{equation}
where $'$ denotes the derivative with respect to the true anomaly $f$ and
$\widetilde{\Omega}$ is given by Equation~\eqref{eq:Omega}, see \cite{Szebehely1967} for the derivation.
There are various formulations of the ERTBP but we have chosen this one because forcing a null eccentricity, $e=0$, in the expressions above, the resulting system corresponds to a Circular Restricted Three-Body Problem, described in previous section.

Similarly to the CRTBP, there is an invariant relation that we can use for the ERTBP (see \cite{Szebehely1967}),
    \begin{equation}
    \label{eq:C_ERTBP}
        \mathcal{C} = \frac{2\tilde{\Omega}}{1+e\cos f} - 2 e \int_{0}^{f} \frac{\tilde{\Omega} \sin f}{(1+e\cos f)^2} df - \left({x'}^2 + {y'}^2 + {z'}^2\right),
    \end{equation}
where  $\tilde{\Omega}$ is given by Equation~\eqref{eq:Omega}. 
Again, forcing null eccentricity, $e=0$, this equation would be equal to the Jacobi constant, introduced in Equation~\eqref{eq:CJ}.

The Elliptic Restricted Three-Body Problem, preserves some symmetries inherited from the CRTBP, like the symmetries when inverting the sense of the time (in this case, the sense of the true anomaly) and the symmetry with respect to the vertical coordinate,
\begin{equation}
\begin{split}
(f,x,y,z,x',y',z') &\rightarrow (-f,x,-y,-z,-x',y',z'),\\[1.ex]
(f,x,y,z,x',y',z') &\rightarrow (-f,x,-y,z,-x',y',-z').
\end{split}
\label{eq:symE}
\end{equation}

\subsection{The Hilda group}\label{sec:selec}

Hilda-type asteroids were classified according to their heliocentric orbital elements by \cite{Zellner1985} as those that possess semi-major axis between $3.7$ and $4.2$ AU, eccentricity lower than $0.3$ and inclination below $20$ degrees. Then, for our analysis we have used the Small Body Database from JPL (\url{https://ssd.jpl.nasa.gov/tools/sbdb_query.html}) to select those asteroids with orbital elements in that classification.
Among these constrains we find more than 5000 objects. However, we must be careful since some of the asteroids that satisfy the described conditions for some time spans do not share the dynamical behaviour of the Hilda group. Notice also that, since orbital elements vary in time, asteroids with orbital elements close to the boundary of these ranges, can have orbital elements crossing back and forth this boundary, and then, being classified (or not) as members of the Hilda group depending on the specific epoch we search for them in the databases.

For example, asteroids (164903) and (210340) would be classified as Hilda asteroids for some epochs but for others their eccentricities overcome the prescribed maximum value of $0.3$. Similarly, semi-major axis of asteroid 2015 PN318 exceeds the lower boundary defined for the Hilda group at some given times. Besides these asteroids that for some periods of time enter and leave smoothly the constrains of the Hilda classification, there are other examples of asteroids whose orbital elements suffer high variations, as the inclination of asteroid 2014 OM449, that reaches values larger than $20$ degrees in an abrupt way. 
This is one of the reasons why we are encouraged to introduce a new characterisation of the Hilda asteroids, according to their dynamics, instead of to their orbital elements.

To start the study,
we first need to start by selecting the asteroids from the database with suitable orbital elements. Then, since our change of coordinates is non-autonomous, we must decide an appropriate date for performing the change. It has been observed that the comparisons between simplified and realistic models through a change of coordinates similar to the one described, is sensible to the date at which the change of coordinates is performed. In particular, the comparisons seem to be better as the scale factor (instantaneous distance between the primaries) is closer to the mean value of their mutual distances \cite{JorbaNicolas2020}. Therefore, we look for a date at which the distance between Sun and Jupiter is closer to  $7.78479\times 10^8$~km, value taken from \url{https://nssdc.gsfc.nasa.gov/planetary/factsheet/jupiterfact.html}. The selected Modified Julian Date --the number of days since midnight on November 17, 1858-- corresponds to $58914$ (6$^{th}$ of March, 2020). 

Here we must clarify that the used database provides ephemeris at rounded dates. In our case, the coordinates of the selected asteroids were downloaded at  Modified Julian Date $59800$, and then propagated using a restricted N-body problem up to the date of interest. From now on, when we talk about restricted N-body problem we refer to the system conformed by the Sun, the eight planets, where Earth and Moon are taken as a single body with mass equal to the sum of both, and the asteroid that is gravitationally affected by the previous ones, but does not affect either of them.

\section{Change of coordinates}\label{sec:CC}

Our goal is to study the motion of the Hilda asteroids in both simplified systems, CRTBP and ERTBP. Since the former is a particular case of the latter, we describe a non-autonomous change of coordinates needed to translate the ephemeris of the asteroids, taken from the JPL database, to the ERTBP. Then, when the interest is focused on the CRTBP, the previous change of coordinates will be suitable just forcing null eccentricity for Jupiter.

First, we describe the cartesian change of coordinates for positions and velocities, similar to other changes of coordinates in the literature \cite{GomezLMS85rep,GomezJMS91rep,JorbaNicolas2021}, involving a rotation and a translation. Second, we include a Keplerian interpretation of that change, where the rotation is expressed in terms of the orbital elements of Jupiter with respect to the Sun. It is noteworthy that both changes provide exactly the same coordinates in the ERTBP.

\subsection{Cartesian change of coordinates}

Let the subscripts $S$ and $J$ denote the Sun and Jupiter, respectively, such that $m_S$ and $m_J$ represent their masses, $\bm{r}_J$ and $\bm{r}_S$ their positions and $\bm{v}_S$ and $\bm{v}_J$ their velocities in the real ecliptic heliocentrical system. The objective is to define a change of coordinates that allows us to translate the dynamics of the Sun-Jupiter system to the mathematical model described in (\ref{eq:ERTBP}), making suitable simplifications on the dynamics like assuming Sun and Jupiter to orbit in Keplerian orbits around their common barycentre.

For simplicity, the change of coordinates is described in the inverse sense, i.e. translating the coordinates, $\bm{x} = (x,y,z)^t$ and $\bm{x}' = (x',y',z')^t$, of the adimensional system with origin in the Sun-Jupiter centre of mass to the real ecliptic ones with origin in the Sun, $\bm{X}$ and $\dot{\bm{X}}$, measured in astronomical units ($AU$) and astronomical units per day ($AU/d$), respectively. Note that, in previous section, devoted to adimensional system, the dot denoted the derivative with respect to the adimensional time, and in this section, it denotes the derivative with respect to the dimensional time of the ecliptic system.
Consequently, the change for positions $\bm{x}$ to $\bm{X}$, is composed by a rotation plus a translation,
\begin{equation}\label{eq:change}
    \bm{X} = r \tilde{R}\bm{x} + \bm{b},   
\end{equation}
where $r\in\mathbb{R}$ is a scale factor, $\tilde{R}$ is an orthonormal matrix and $\bm{b}$ is a translation of the origin of coordinates from the Sun to the Sun-Jupiter barycenter,
\[
    \bm{b} = (1-\mu)\bm{r}_S + \mu \bm{r}_J.
\]
Before continuing, let us define 
\begin{equation*}
\bm{r} = \bm{r}_S-\bm{r}_J, \quad \bm{v} = \bm{v}_S-\bm{v}_J, \quad 
 \bm{h} = \bm{r}\times\bm{v},
\end{equation*}
as the relative position ($\bm{r}$), the relative velocity ($\bm{v}$) of the Sun with respect to Jupiter, and $\bm{h}$ as their angular momentum. 
Moreover, we denote with $~\hat{}~$ the unitary vector, i.e. $\hat{\bm{b}} = \bm{b}/||\bm{b}||$
and $b=||\bm{b}||$.

Since we are considering an instantaneous change of coordinates between the two systems, the scale factor is given by the instantaneous distance between the Sun and Jupiter, $r = ||\bm{r}||$. It is easy to see that the $x$ axis of the simplified system corresponds to the direction Sun-Jupiter, $z$ axis corresponds to the vector normal to the plane of motion and $y$ axis corresponds to a vector perpendicular to $x$ and contained in the plane of motion, such that the system is positively oriented. With this, the rotation matrix is defined as
\begin{equation}\label{eq:RotMat}
\tilde{R} = \begin{pmatrix}
        | & | &  |\\ 
        \bm{R}_1 & \bm{R}_2  & \bm{R}_3\\
        | & |  & |\\
    \end{pmatrix},
\end{equation}
where $\bm{R}_1=\hat{\bm{r}}$, $\bm{R}_2 = \hat{\bm{h}}\times\hat{\bm{r}}$ and $\bm{R}_3 = \hat{\bm{h}}$.

Once described the change of positions, we need to describe the change of velocities. For this, it is usual to derive the change of positions (\ref{eq:change}) with respect to the dimensional time of the ecliptic system, 
\begin{equation*}
    \dot{\bm{X}} 
    = \frac{d}{dt}\left(r \tilde{R}\bm{x} + \bm{b}\right)
    = \dot{r}\tilde{R}\bm{x} + r \dot{\tilde{R}}\bm{x} +r \tilde{R}\dot{\bm{x}} +\dot{\bm{b}},\\
\end{equation*}
where the dot denotes the derivative with respect to dimensional time.

However, in the elliptic problem, the velocities are considered as the derivatives of the positions with respect to the true anomaly, denoted by $\bm{x}'$. Then, the relation between $\dot{\bm{x}}$ and $\bm{x}'$ is given by 
\begin{equation}\label{eq:dfdt}
\dot{\bm{x}} = \frac{df}{dt} \bm{x}',    
\end{equation}
being the variation of the true anomaly with respect to time $\frac{df}{dt}$, defined in (\ref{eq:df}), of the form
$$\frac{df}{dt} = \frac{k\sqrt{m_S+m_J}}{a^{3/2}(1-e^2)^{3/2}}\left(1+e\cos f\right)^2,$$
where $a$, $e$ and $f$ correspond to the semi-major axis, eccentricity the true anomaly of the motion of either of the primaries, Sun or Jupiter, with respect to the other. 

Introducing (\ref{eq:dfdt}) in the change of velocities from $\dot{x}$ to $\dot{\bm{X}}$, the change from $\bm{x}'$ to $\dot{\bm{X}}$ can be written as
\begin{equation}\label{eq:changev}
    \dot{\bm{X}} = \dot{r}\tilde{R}\bm{x} + r \dot{\tilde{R}}\bm{x} +r \tilde{R}\frac{df}{dt}\bm{x}' +\dot{\bm{b}}.\\
\end{equation}
This expression involves the velocity of the barycentre $\dot{\bm{b}} = (1-\mu)\bm{v}_S + \mu \bm{v}_J$, the temporal derivative of the scale factor, 
$\dot{r} = \hat{\bm{r}}\cdot\bm{v}$, 
and the temporal derivative of the rotational matrix,
\begin{equation*}
    \dot{\bm{R}}_1 = \frac{r \bm{v}-\dot{r}\bm{r}}{r^2}, \ \ \ \
    \dot{\bm{R}}_2 = \bm{R}_3\times \dot{\bm{R}}_1, \ \ \ \
    \dot{\bm{R}}_3 = \bm{0}.  
\end{equation*}
We note that going from the complete Solar system to the Sun-Jupiter-asteroid ERTBP at a given time is like removing from the complete Solar system all the remaining planets (and asteroids) because, from that time on, Sun and Jupiter will move on ellipses.
Then, the motion is determined by a two-body system and the angular momentum is conserved, being its derivative zero, $\dot{\bm{h}}=\bm{0}$. As a consequence, the third column of matrix $\dot{\tilde{R}}$ is null. Also, this property must hold for equation (\ref{eq:dfdt}) to be true.

In order to perform the described change of coordinates, some values are needed, as the value of the true anomaly, $f$, the term $(1+\cos f)$ and the semi-major axis, $a$. For this, it is necessary to compute the eccentricity vector $\bm{e}$,
\[
    \bm{e} = \frac{\bm{v}\times\bm{h}}{\mu^*}- \hat{\bm{r}},
\]
being $\mu^*=k^2(m_S+m_J)$, since the true anomaly is obtained from
$\cos f = \hat{\bm{e}}\cdot \hat{\bm{r}}$. 
In the case in which $\bm{r}\cdot\bm{v}<0$, the value of the true anomaly must be $f=2\pi-f$. 
With some manipulations, the value of the term $(1+e\cos f)$ 
and the semi-major axis can be computed as
\[
1+e\cos f = \frac{\bm{v}\times\bm{h}}{\mu^*} \cdot \hat{\bm{r}},
\qquad
a 
= \frac{1}{\frac{2}{||\bm{r}||}-\frac{||\bm{v}||^2}{\mu^*}} 
\]

With all these considerations, the inverse change of coordinates of \eqref{eq:change} and \eqref{eq:changev}, from the ecliptical coordinates with the origin in the Sun ($\bm{X}$ and $\dot{\bm{X}}$) to the elliptic adimensional ones with the origin at the Sun-Jupiter barycentre ($\bm{x}$ and $\bm{x}'$), can be written as:
\begin{equation}\label{eq:ChangeCord}
    \left\{
    \begin{aligned}
    \bm{x} &= \frac{1}{r}\tilde{R}^{-1}(\bm{X}-\bm{b})\\[1.5ex]
    \bm{x}' & = \frac{1}{r}\frac{dt}{df}\tilde{R}^{-1}\left[
    \left(\dot{\bm{X}}-\dot{\bm{b}}\right) - \frac{\dot{r}}{r}(\bm{X}-\bm{b}) -
    \dot{\tilde{R}}\tilde{R}^{-1}(\bm{X}-\bm{b})
    \right],\\
    & =\frac{1}{r}\tilde{R}^{-1}\left[
    \frac{dt}{df}\left(\dot{\bm{X}}-\dot{\bm{b}}\right) - \frac{r'}{r}(\bm{X}-\bm{b}) -
    \tilde{R}'\tilde{R}^{-1}(\bm{X}-\bm{b})
    \right].
    \end{aligned}
    \right.
\end{equation}
Last expression is included for convenience when we use the Keplerian interpretation of the change of coordinates, explained below.

\subsection{Keplerian interpretation}

The change of coordinates previously detailed can be easily expressed in terms of the orbital elements of Jupiter with respect to the Sun. Let us consider that Jupiter orbit is given by the elements $(a,e,i,\Omega,\omega,f)$, that correspond to its semi-major axis, eccentricity, inclination, longitude of the ascendent node, argument of the periapsis and the true anomaly.
Then, the rotation matrix \eqref{eq:RotMat} used in the definition of the change of coordinates \eqref{eq:ChangeCord} can be written in terms of the orbital elements as,
\[
\tilde{R} = \tilde{R}_3(\Omega)\,\tilde{R}_1(i)\,\tilde{R}_3(\omega+f-\pi),
\]
where
\[
\begin{aligned}
    \tilde{R}_3(\alpha) &= 
    \begin{pmatrix}
    \cos\alpha & -\sin\alpha & 0\\
    \sin\alpha & \cos\alpha & 0\\
    0 & 0 & 1
    \end{pmatrix},\\
    \tilde{R}_1(\alpha) &= 
    \begin{pmatrix}
    1 & 0 & 0\\
    0 & \cos\alpha & -\sin\alpha\\
    0 & \sin\alpha & \cos\alpha
    \end{pmatrix},
\end{aligned}
\]
that corresponds to the product of three rotations. The first one accounts for the plane of motion of Jupiter's orbit, the second one for its inclination and the third one places Jupiter in its orbit. It is worth to mention that the term $\omega+f-\pi$ has this form due to the consideration that Jupiter occupies the position $(\mu-1,0,0)$ and the Sun the position $(\mu,0,0)$ both, in the ERTBP or in the CRTBP. In the case that one wants to consider the positions $(1-\mu,0,0)$ and $(-\mu,0,0)$ for the primaries, the angle of the third matrix must be $\omega+f$.

Since the scale factor corresponds to the instantaneous distance between the primaries, it also admits to be written in terms of the Keplerian coordinates as in equation (\ref{eq:r_ellip}), that again, is only true for $\dot{\bm{h}}=\bm{0}$.

Finally, note that writing the rotation matrix and the scale factor in terms of the orbital elements, allows to easily obtain the computation of $\tilde{R}^{-1}$ and the derivatives with respect to the true anomaly of $\tilde{R}$ and $r$.

\section{Analysis in the planar CRTBP}\label{sec:CRTBP}
\begin{figure}[!htbp]
    \centering
    \includegraphics[width=1.\columnwidth]{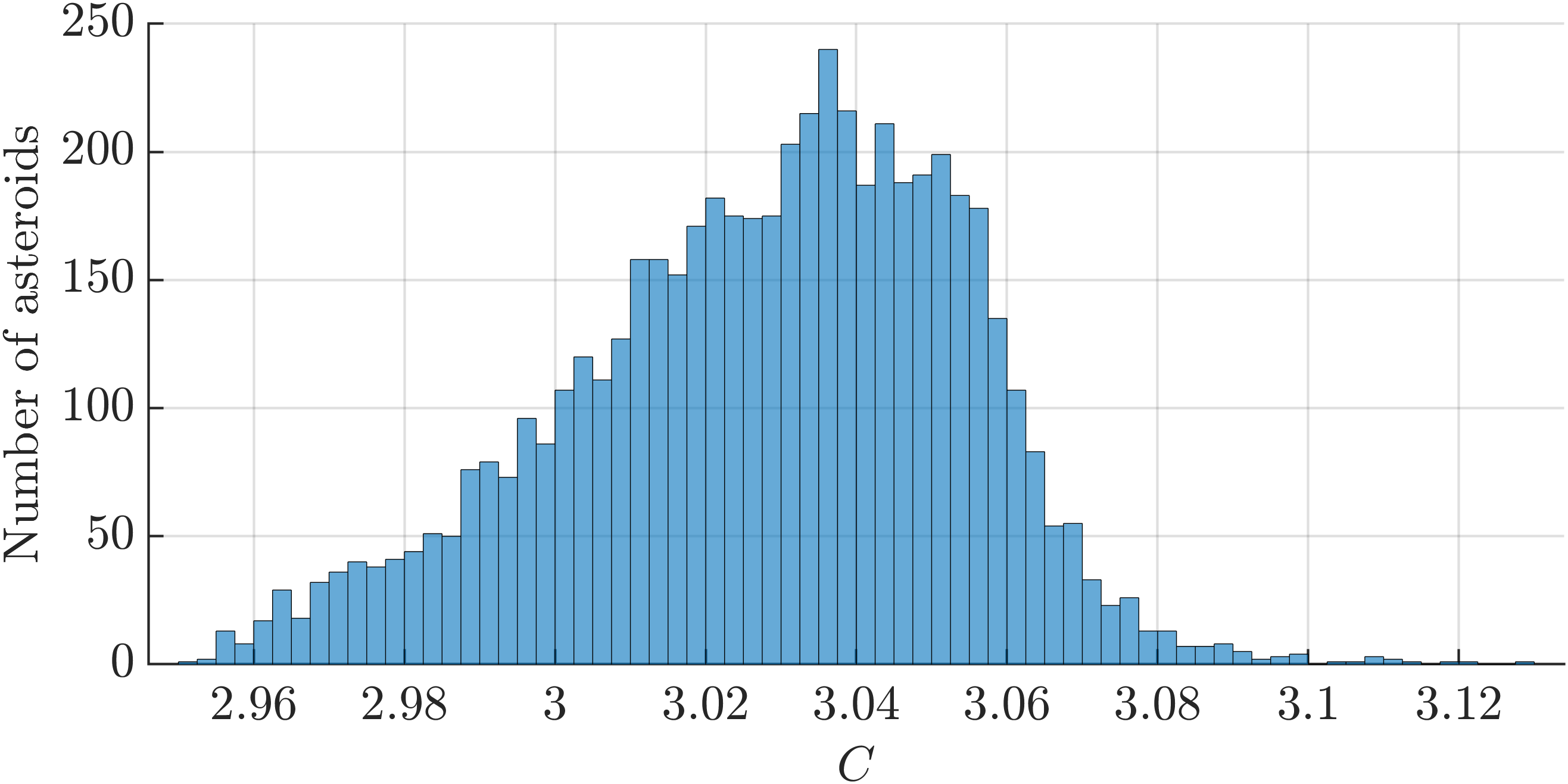}
    \caption{Histogram of the Jacobi constant, $C$, for the Hilda asteroids.}
    \label{fig:Cs}
\end{figure}

First we proceed with the analysis of Hilda asteroids in the planar Circular Restricted-Three Body Problem. For this, we apply the change of coordinates described in Section~\ref{sec:CC} forcing null eccentricity for the primaries and, then, we project the resulting Cartesian coordinates to the $xy$-plane.
Then, we compute their Jacobi constant~\eqref{eq:CJ} and we plot their distribution in Figure~\ref{fig:Cs}.
We observe that most of them have $C$ values between $2.98$ and $3.06$. In order to give some references, it is worth to mention that the Jacobi constant for the two most representative asteroids of this group has value of $C\approx3.05021$ for (153) Hilda and $C\approx3.04311$ for (1911) Schubart,
that the Jacobi constant of the triangular points $L_4$ and $L_5$ is exactly $3$ and the value for the collinear point $L_3$, for the mass parameter of Sun-Jupiter system, is approximately $C_{L_3}\approx3.00272$.

Here it is worth to mention that in the planar Sun-Jupiter CRTBP, $L_3$ is a centre $\times$ saddle point, which implies that there exists a continuous family of Lyapunov periodic orbits growing in the centre direction of the point. This family of periodic orbits constitutes the centre manifold of $L_3$. Moreover, each periodic solution in that family has stable and unstable invariant manifolds, that could explain the approaches of Hilda asteroids to the neighbourhood of $L_3$. For this reason, in some occasions the centre manifold of $L_3$ (and its invariant manifolds) were thought to be responsible for the dynamics of the Hilda group of asteroids. However, the Jacobi constant of the periodic orbits in this Lyapunov family decreases when we move away from $L_3$, which means that these periodic orbits (and their invariant manifolds) have a too low value of $C$ to be related to the Hilda group of asteroids, at least in the frame of the planar Sun-Jupiter CRTBP.

\begin{figure}[!htbp]
    \centering
    \includegraphics[width=1\columnwidth]{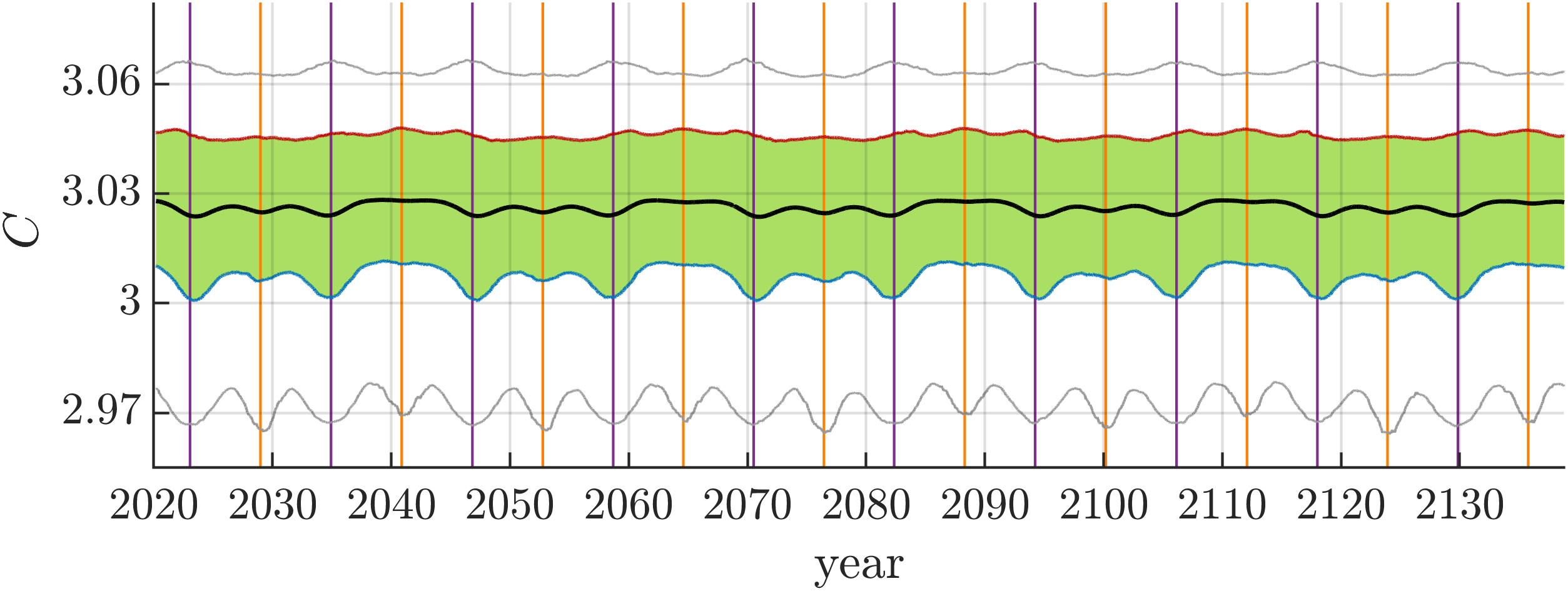}
    \caption{Evolution in time of the Jacobi constant for Hilda asteroids in a restricted $N$-body problem. Black line marks the medium value ($Q_2$), blue one the first quartile ($Q_1$) and red one the third ($Q_3$). Green area corresponds to $IQR=Q_3-Q_1$. Gray lines mark percentiles $P_5$ and $P_{95}$.
    }
    \label{fig:CsEvol}
\end{figure}

In order to strengthen the relation between Hilda asteroids and the values of Jacobi constant in the Sun-Jupiter CRTBP we plot, in Figure~\ref{fig:CsEvol}, the evolution in time of the Jacobi constant in a restricted planetary $N$-body problem. All the selected Hilda asteroids, whose coordinates are taken in the year 2020, are propagated forward in time under the gravitational effect of all the planets and Sun until year 2140. Vertical lines in purple and in orange correspond to times at which the true anomaly of Jupiter takes values $f=0$ and $f=\pi$, respectively. Then, the value of $C$ has been computed for all of the asteroids at different times, marking with a black line the medium value and in green area the interquartile range ($IQR$). As we can observe, the Jacobi constant in this realistic model, and at different times, varies in the same ranges found for the Sun-Jupiter CRTBP. Then, it is justified to relate the Hilda asteroids with the range of $C$ shown in Figure~\ref{fig:Cs}.

\subsection{Periodic solutions}\label{sec:periodic}

\begin{figure}
    \centering
    \includegraphics[width=1\columnwidth]{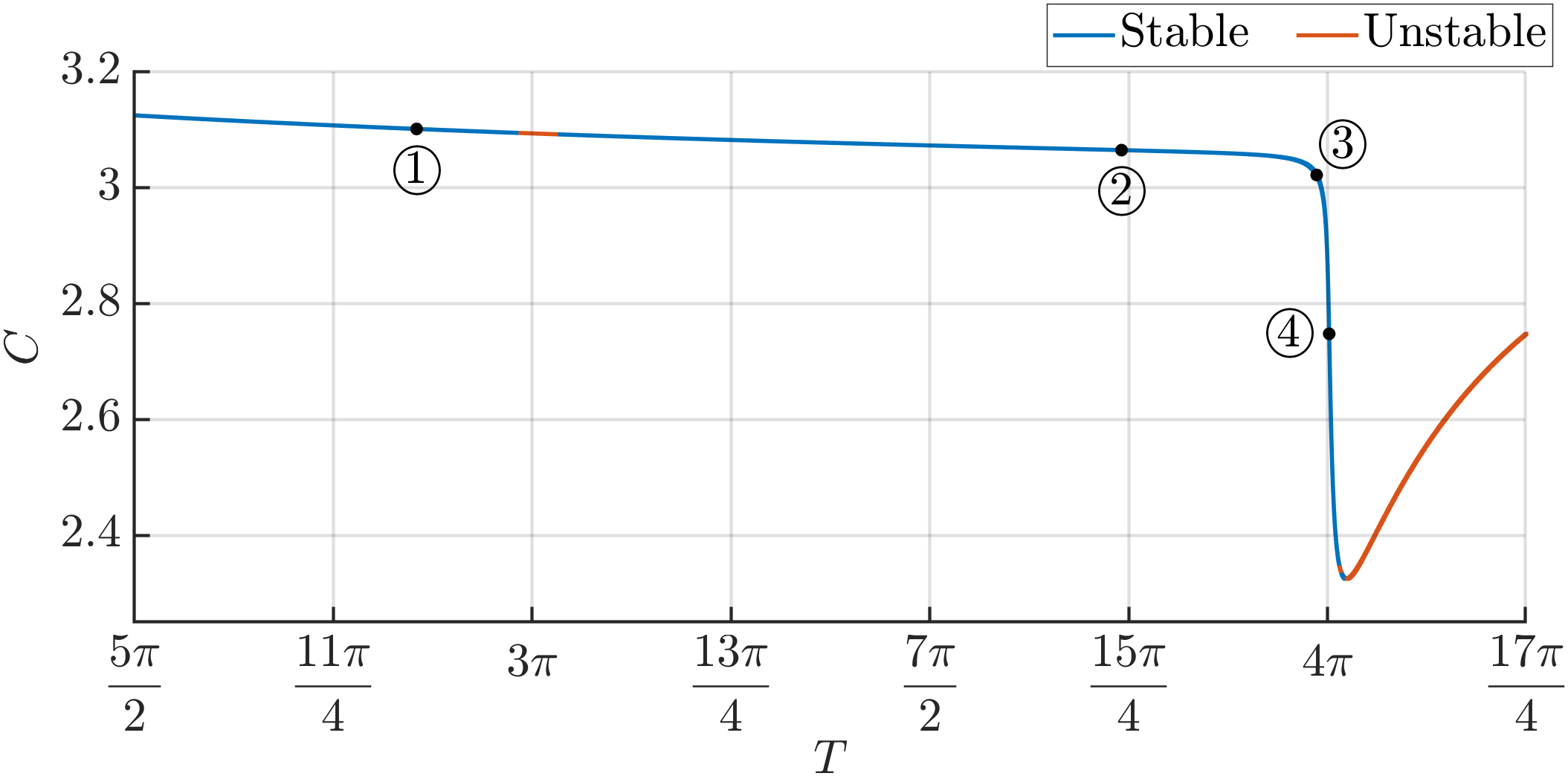}
    \includegraphics[width=0.48\columnwidth]{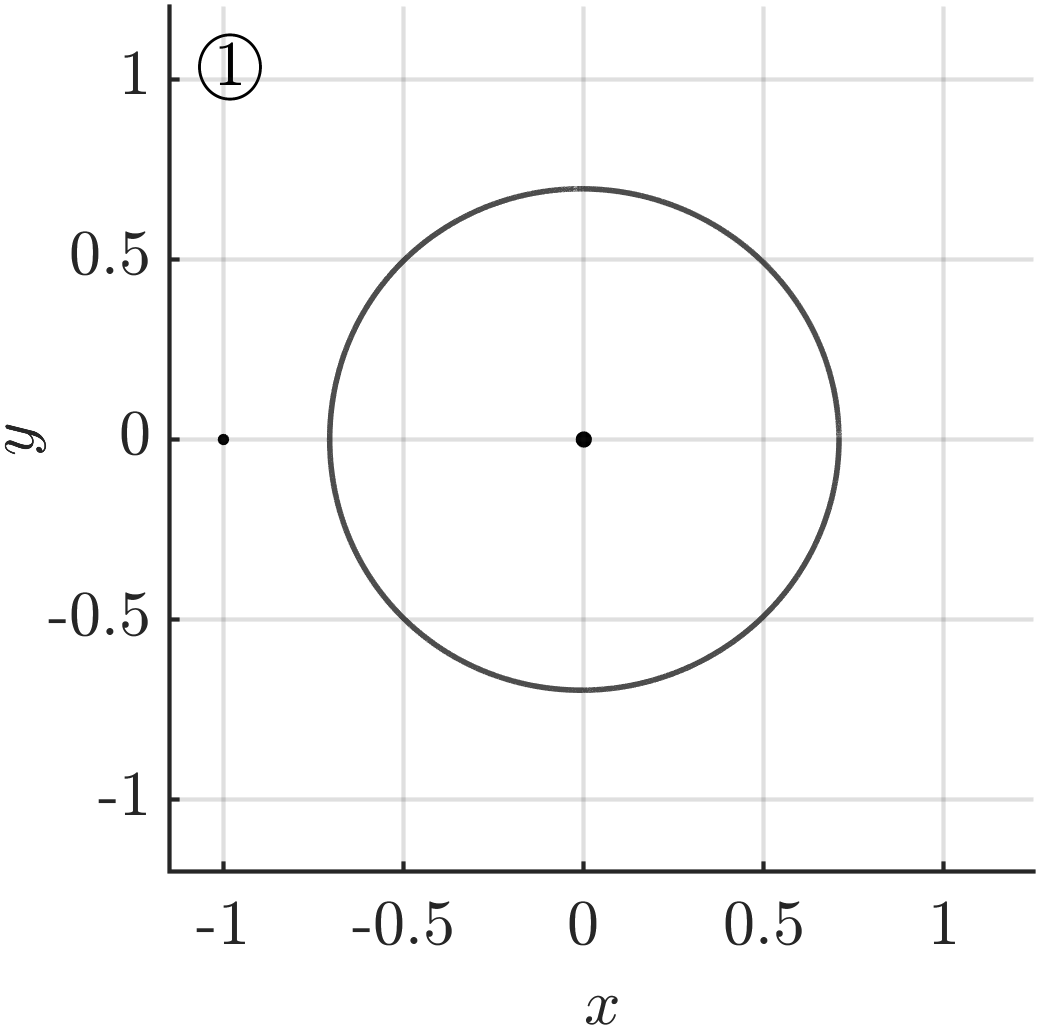}
    \includegraphics[width=0.48\columnwidth]{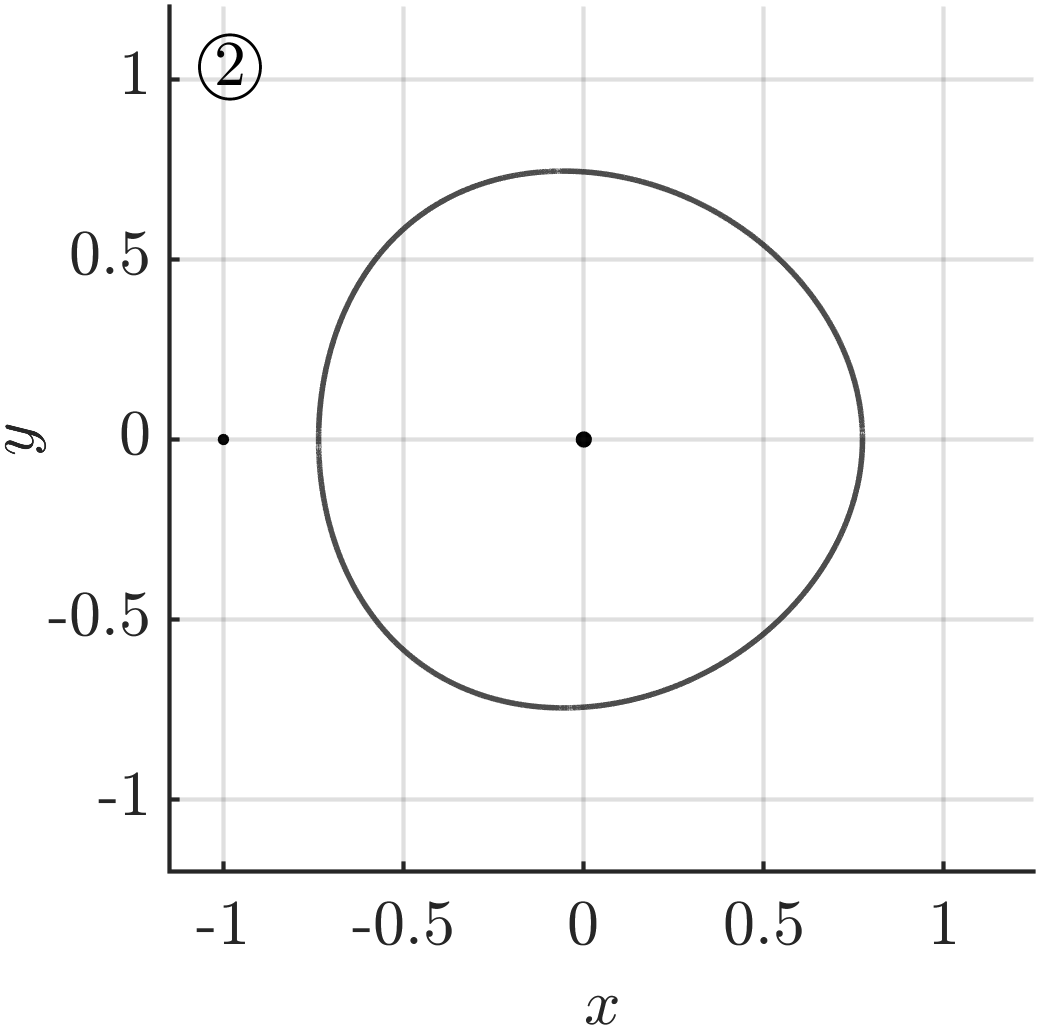}
    \includegraphics[width=0.48\columnwidth]{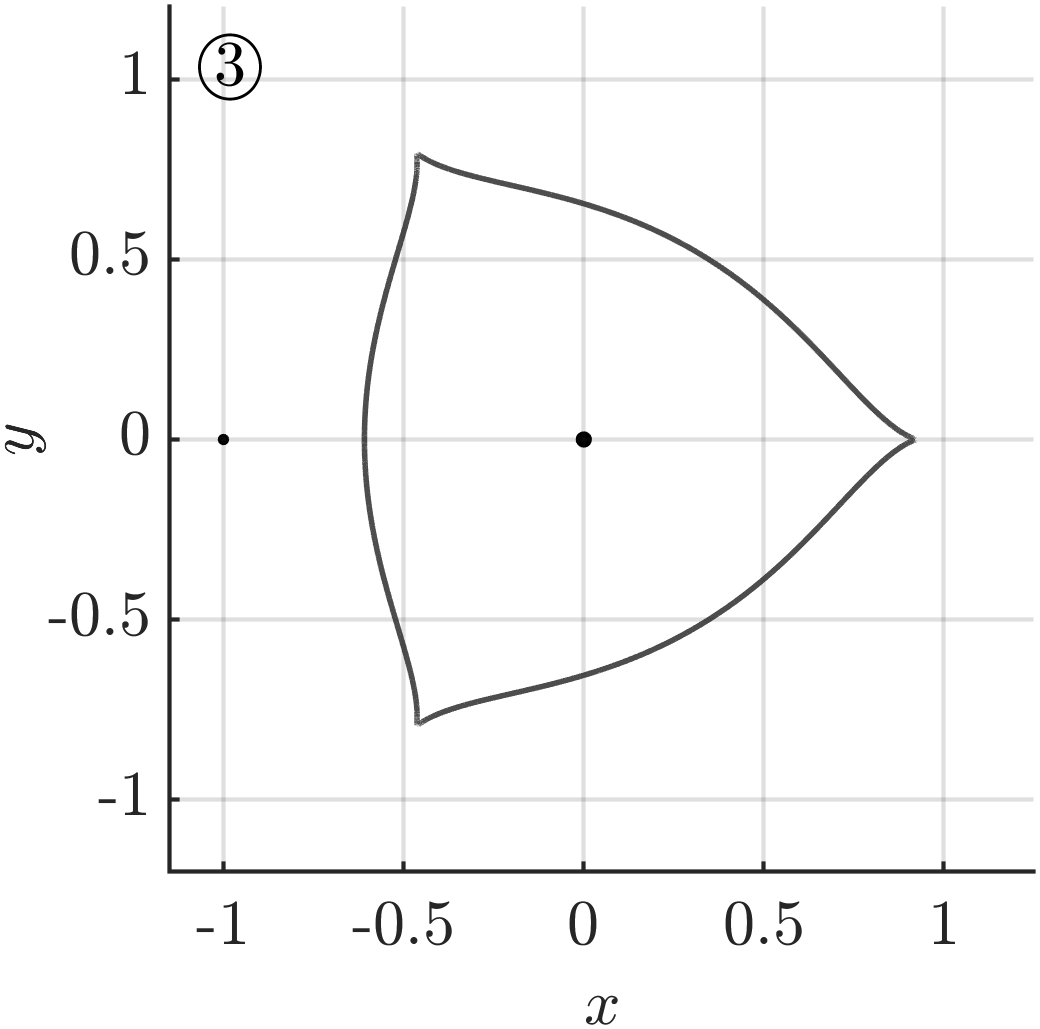}
    \includegraphics[width=0.48\columnwidth]{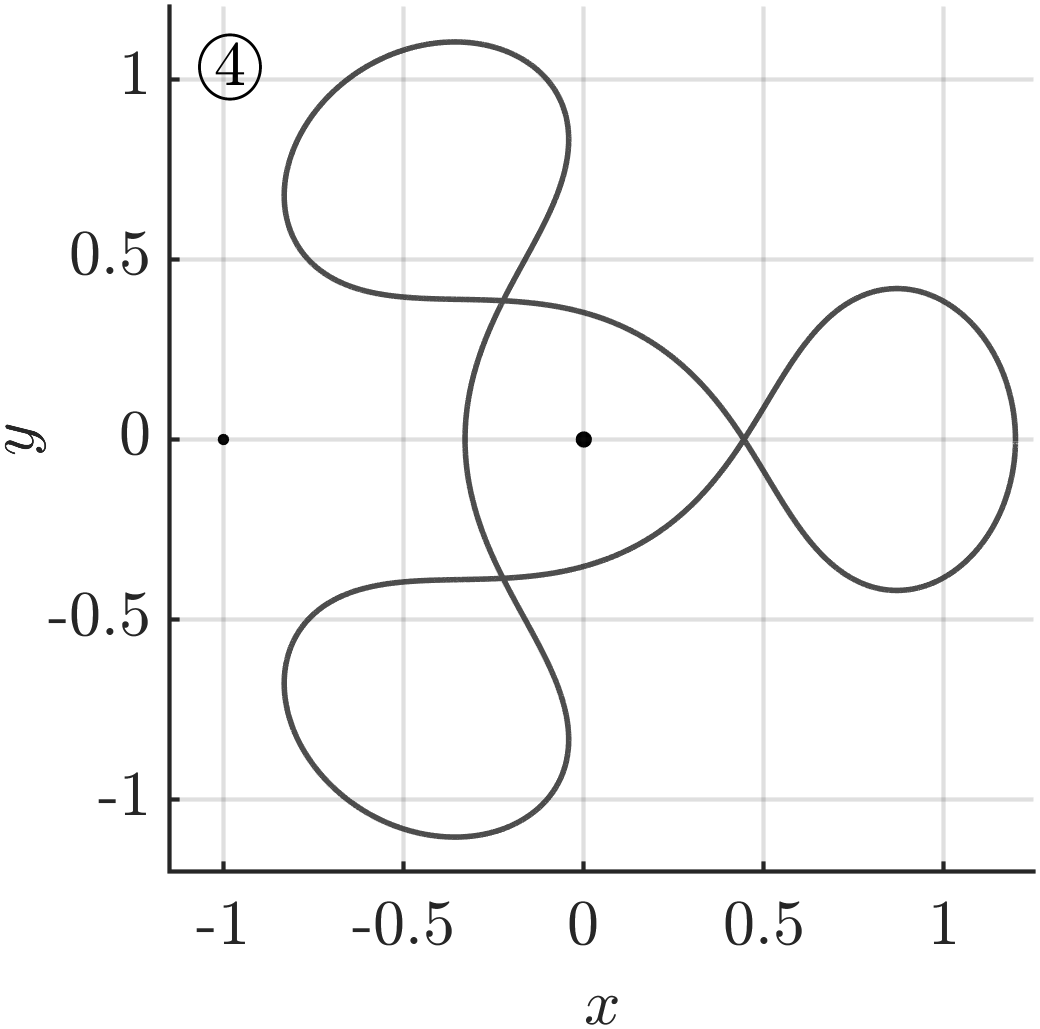}
    \caption{Up, variation of the $C$ value according to the period $T$ of the periodic orbits in the family around the first primary in the CRTBP. Two bottom rows, periodic orbit in the $xy$-plane, for each pair of values ($T$,$C$) numbered in the figure above.}
    \label{fig:OPs}
\end{figure}

In the CRTBP a family of periodic orbits around the first primary, the Sun for us, is known to exist, composed by periodic orbits that sweeps a large range of $C$ values, depending on its period $T$. This variation is shown in Figure~\ref{fig:OPs}, along with the projection in the $xy$-plane of four periodic orbits of the family, to illustrate the different behaviour of the orbits as the values of $C$ and $T$ vary. Most of the periodic orbits in this family are stable (coloured in blue in the figure) and some of them are unstable (coloured in red). For example, when the period approaches $3\pi$ these periodic orbits are unstable. Moreover, for periods slightly larger than $4\pi$ there is a small window of instability and, for larger values of $T$, the family becomes unstable.

\begin{figure}
    \centering
    \includegraphics[width=0.48\columnwidth]{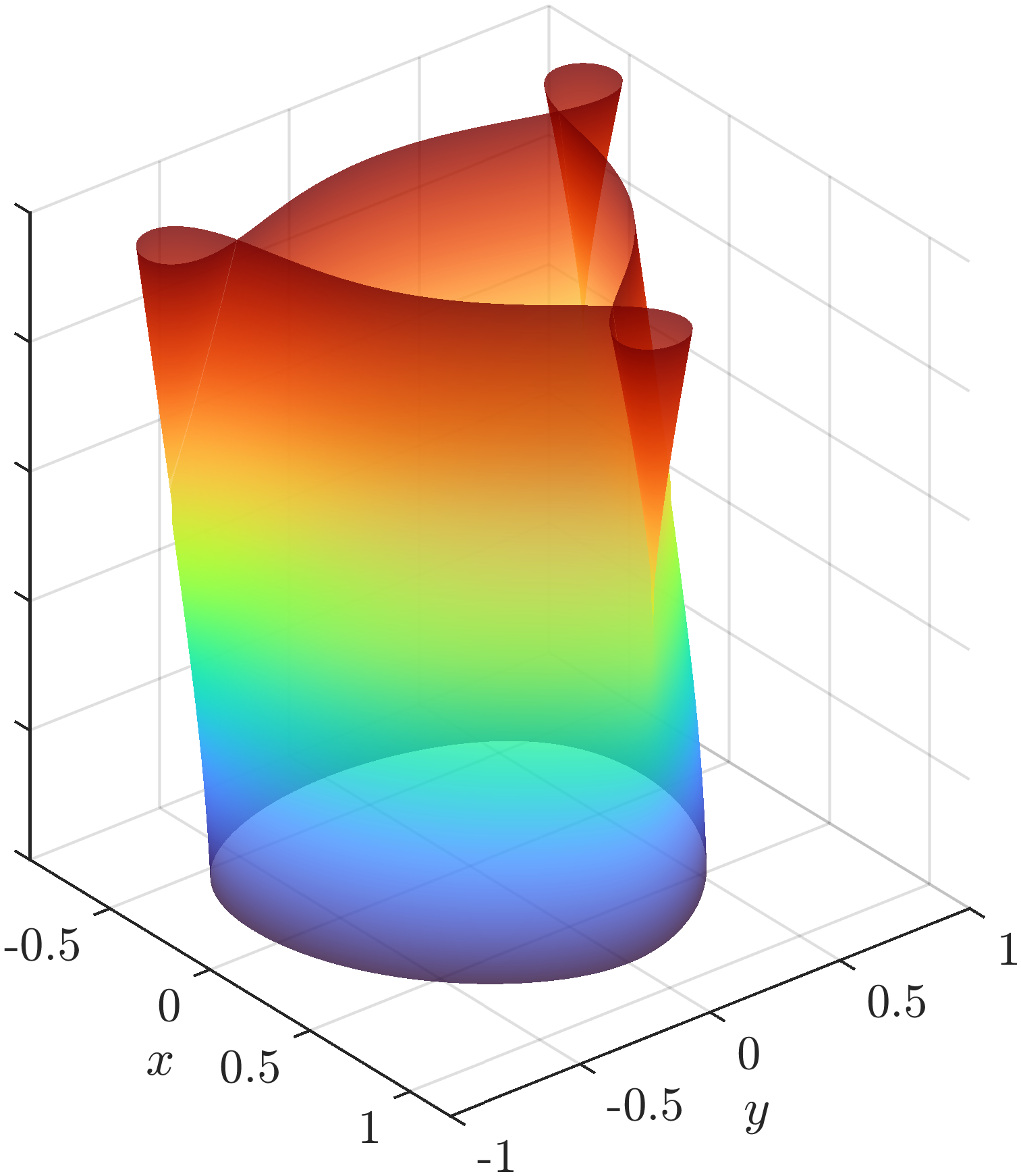}
    \hspace{2mm}
    \includegraphics[width=0.48\columnwidth]{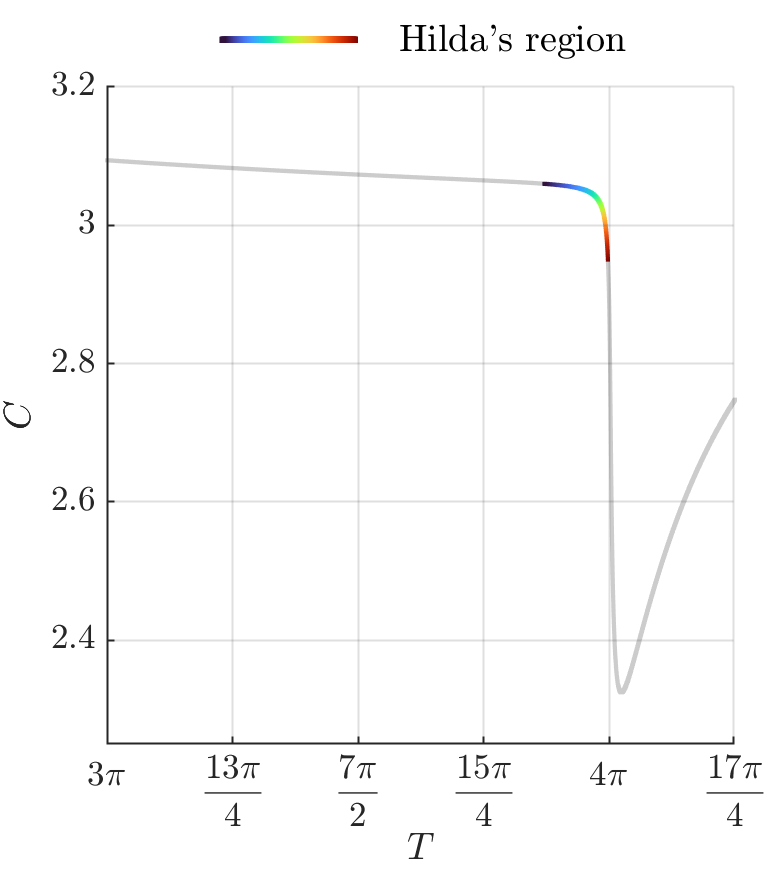}
    \caption{Left, surface composed by a continuum of periodic orbits projected in the $xy-$plane for the energy range of the Hilda group. Colour is used to identify each periodic orbit with the corresponding point in the $(T,C)$ plot at the right, where the variation of the $C$ value according to the period $T$ of the periodic orbits is shown.}
    \label{fig:OPsHildas}
\end{figure}

In spite of the wideness of this family of periodic orbits, our interest is focused on a narrow range of $C$, since, as mentioned before, the values related to the Hilda group are approximately between $2.98$ and $3.06$. Notice that this piece of the family is stable as shown in Figure~\ref{fig:OPs}. In Figure~\ref{fig:OPsHildas}, this range is marked and coloured in a ($T$,$C$) plot, such that the colour is used to identify the corresponding periodic orbit shown in the left of that figure. Notice that the plot on the left of Figure~\ref{fig:OPsHildas} shows a surface composed by a continuum of curves, that corresponds to the periodic orbit projected in the $xy-$plane for each value of the Jacobi constant considered. In this way, a transversal cut of said surface gives the periodic orbit for a given value of $C$. It is observed that for higher values of $C$, in blue colour, the shape of the periodic orbit is soft and rounded, meanwhile as $C$ decreases three peaks approaching the locations of $L_3$, $L_4$ and $L_5$ appear until the peaks become loops for the lowest values of the Jacobi constant considered, in red colour.

\begin{figure}[!ht]
    \centering
    \includegraphics[width=0.48\columnwidth]{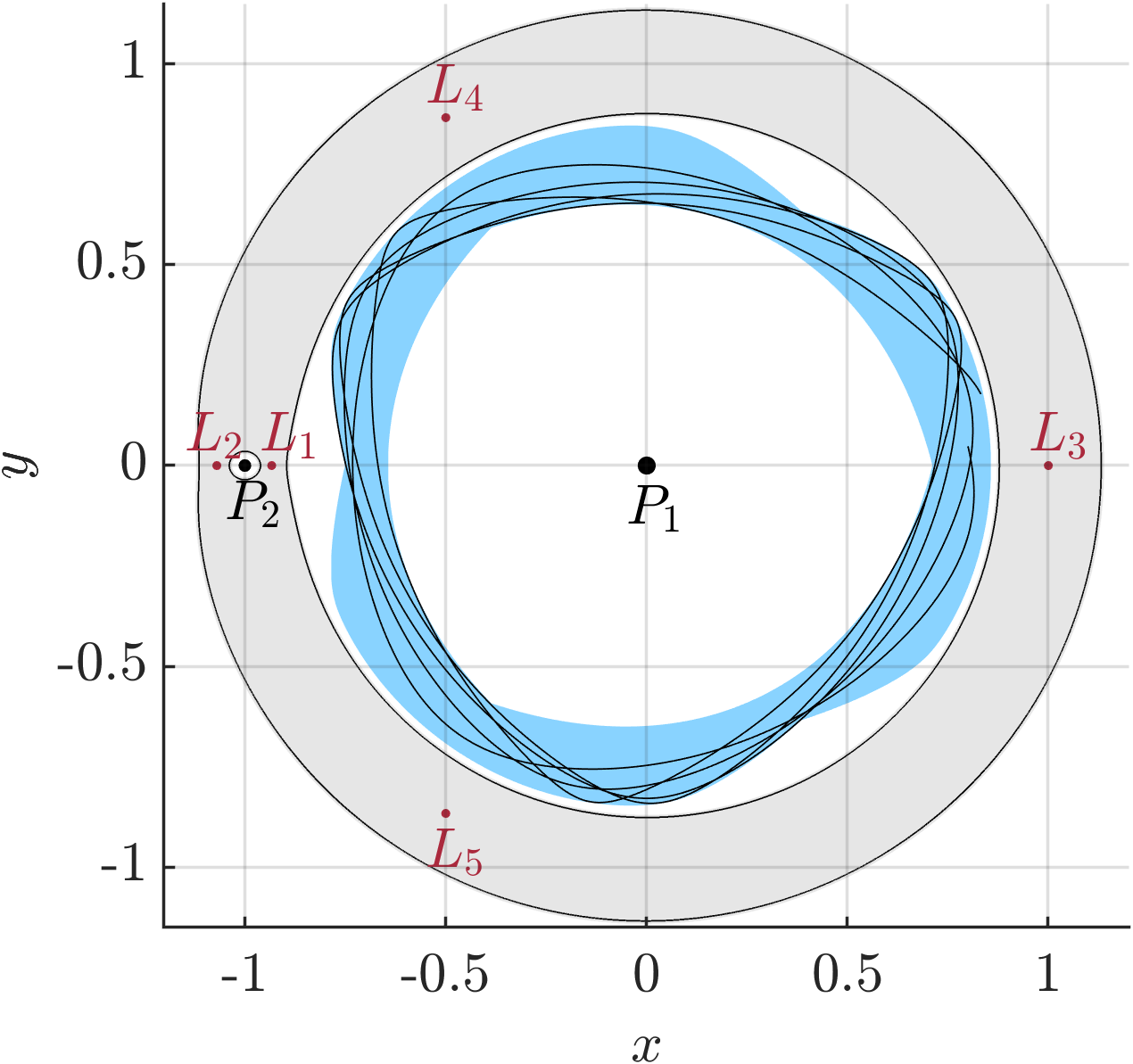}
    \includegraphics[width=0.48\columnwidth]{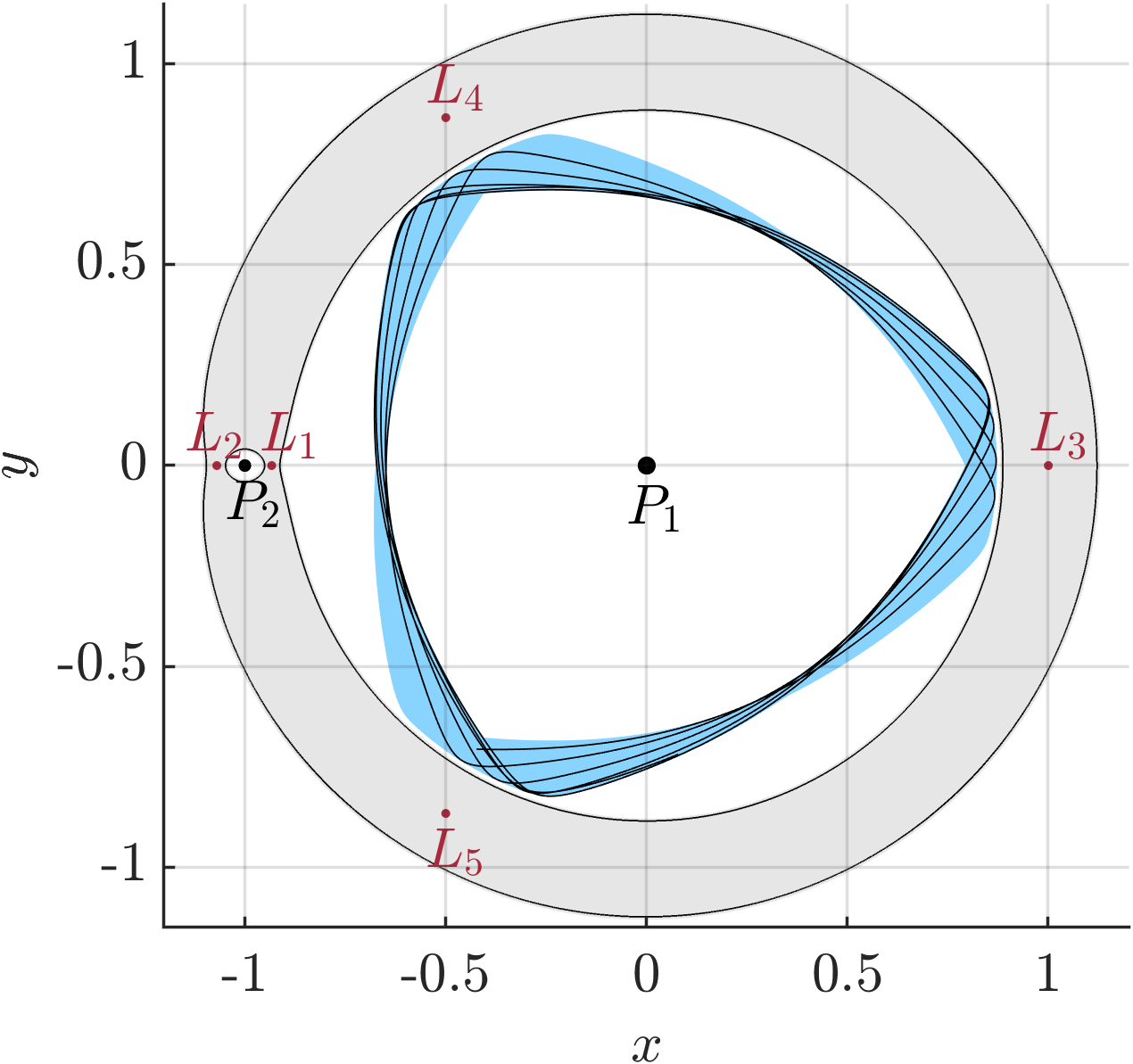}
    \includegraphics[width=0.48\columnwidth]{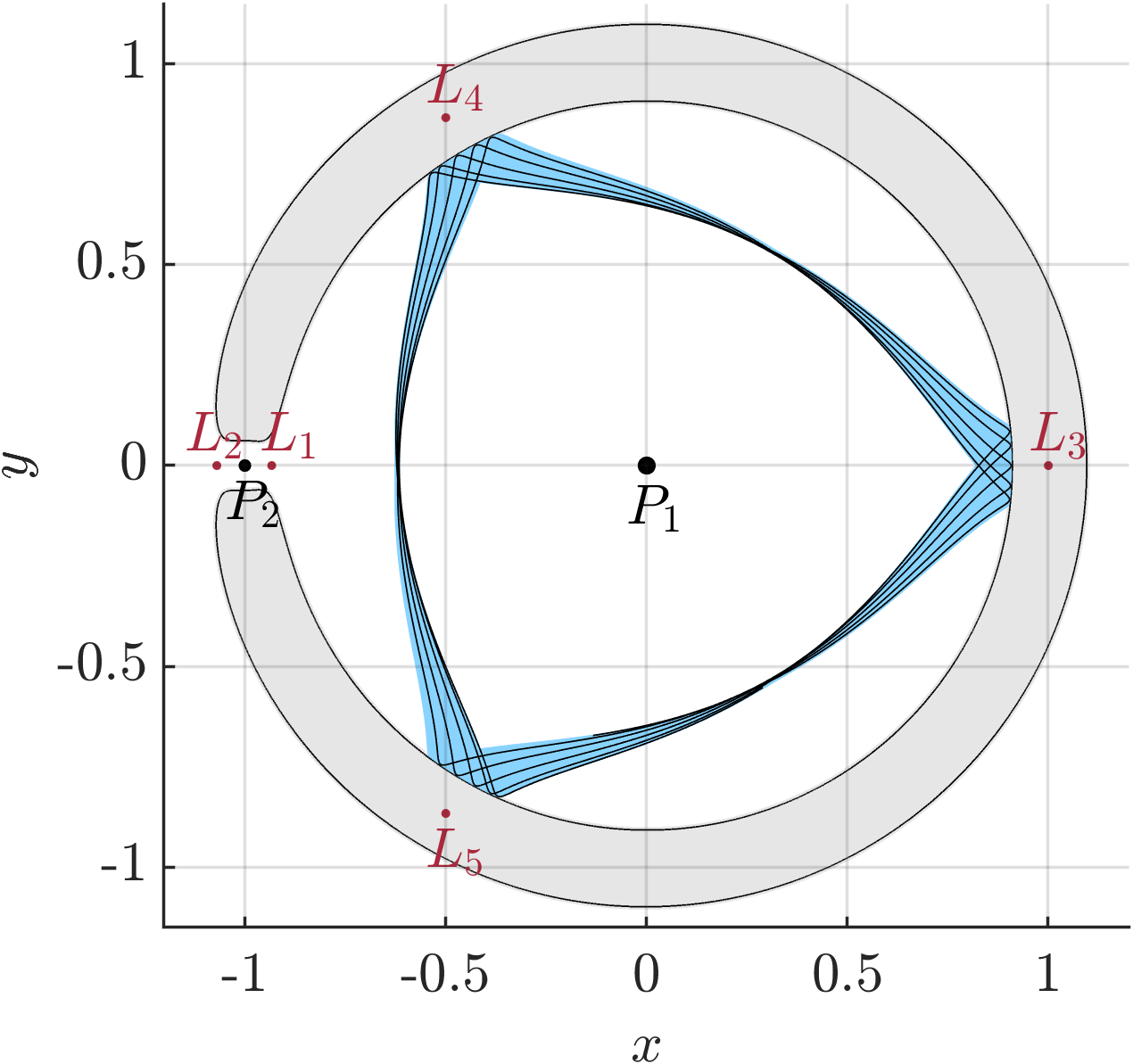}
    \includegraphics[width=0.48\columnwidth]{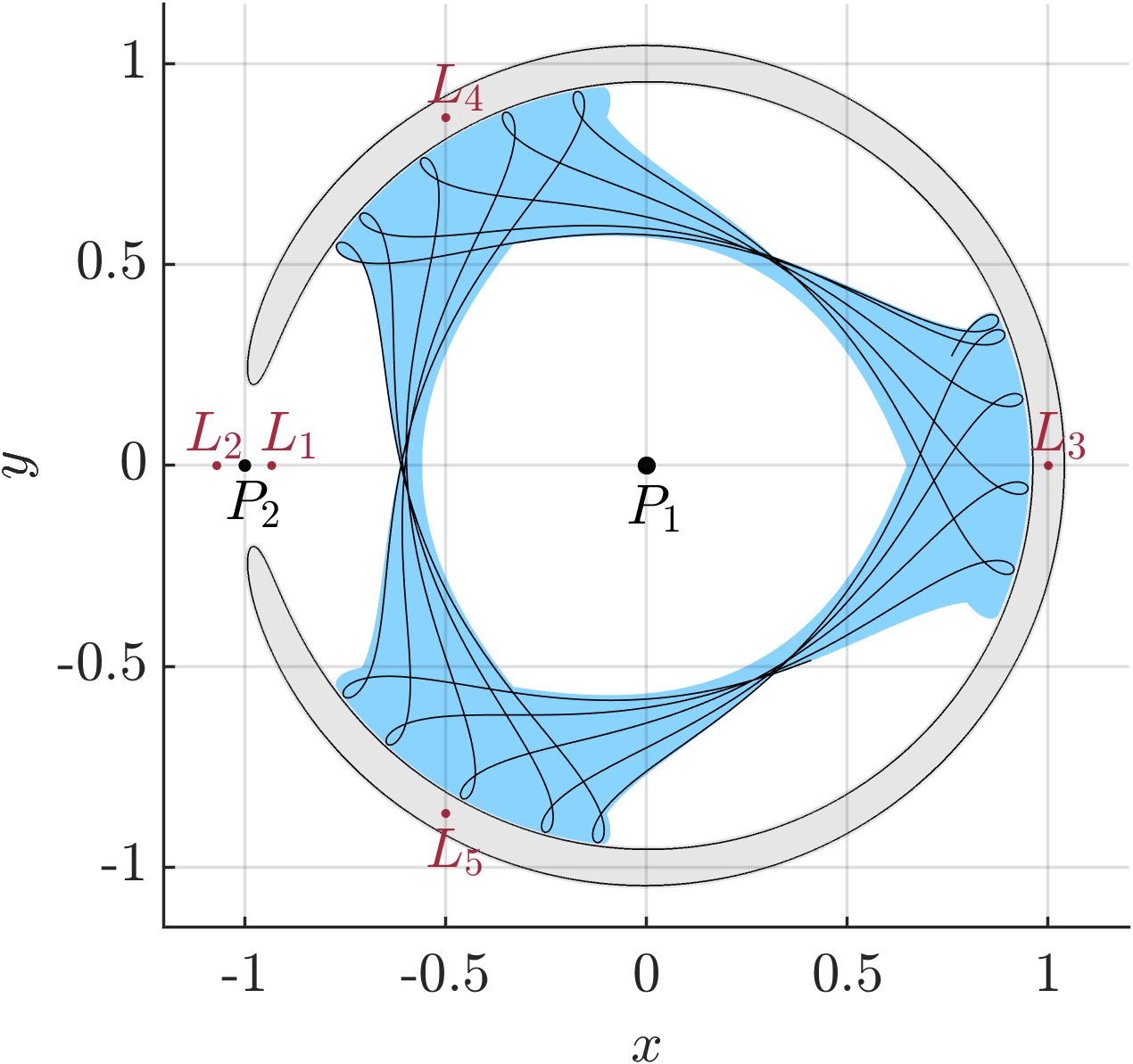}
    \includegraphics[width=0.48\columnwidth]{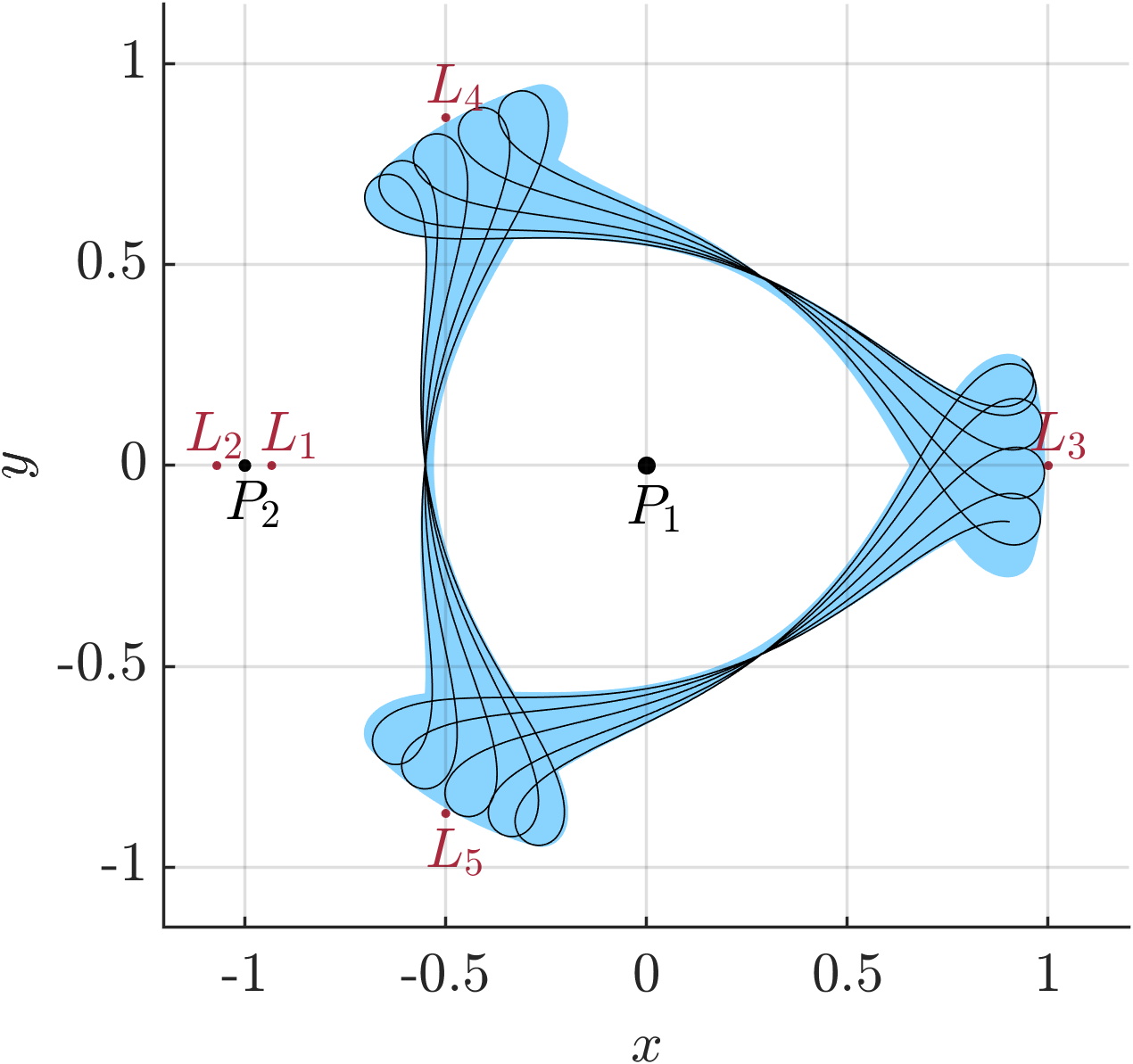}
    \includegraphics[width=0.48\columnwidth]{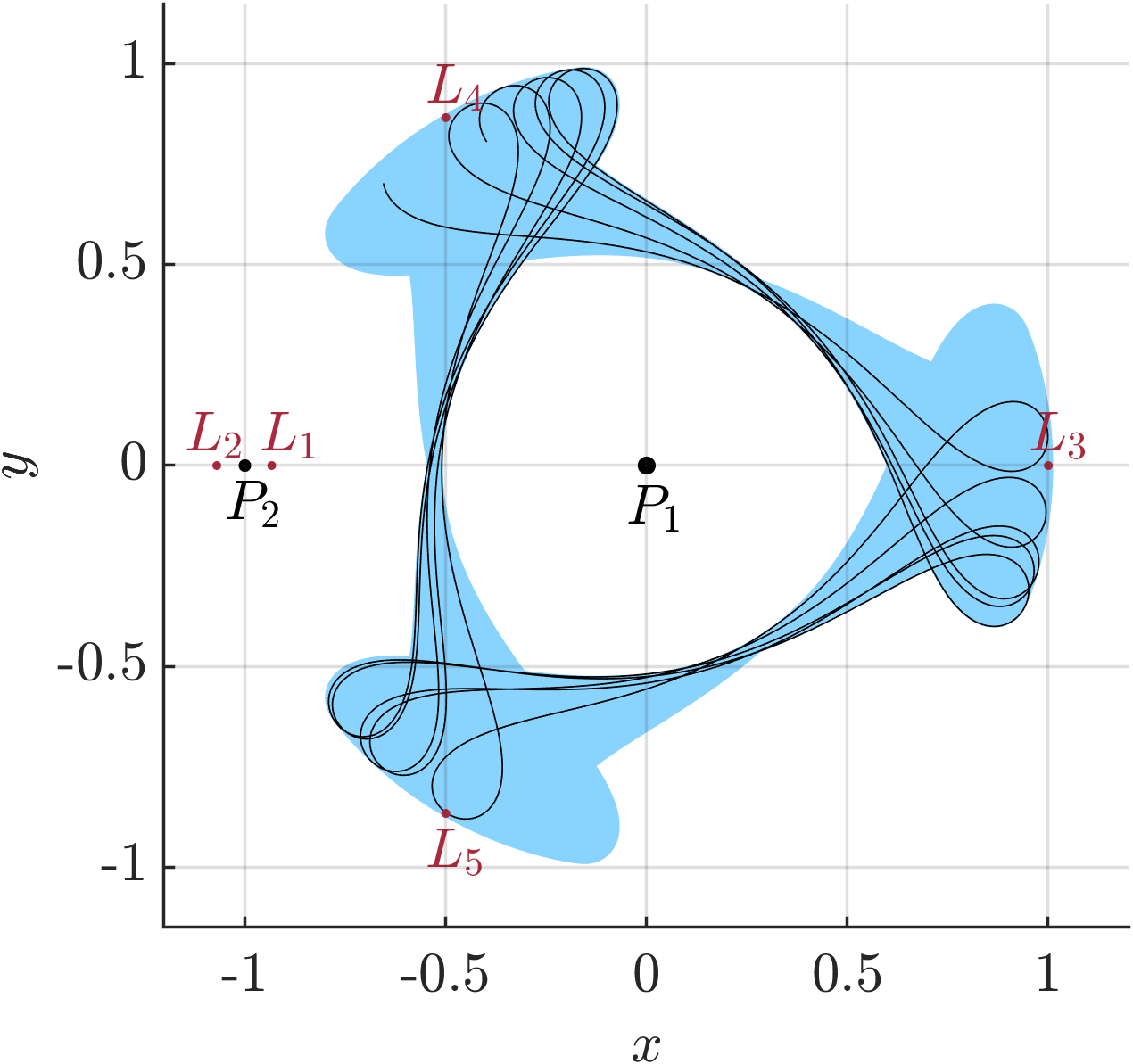}
    
    \caption{
    Trajectories of six different asteroids. From left to right: first row, (153) Hilda, $C\approx3.05021$, and (1911) Schubart, $C\approx3.04311$, second row, (499) Venusia, $C\approx3.02800$, and (1038) Tuckia, $C\approx3.00662$, and third row, (4446) Carolyn, $C\approx2.98210$, and (2483) Guinevere, $C\approx2.96911$. See text for further details.  
    }
    \label{fig:trasHildas}
\end{figure}

We also observe this variation of the behaviour of the orbits according to the value of $C$ when we integrate in the planar CRTBP the Hilda-type asteroids collected in Figure~\ref{fig:Cs}. In Figure~\ref{fig:trasHildas} the trajectories of six asteroids at different energy levels (different values of $C$) are shown.  The orbits are plotted in black in the $xy$-plane together with the Lagrangian points, the primaries (being $P_1$ the Sun and $P_2$ Jupiter) and the forbidden region for each value of $C$, marked in grey. The blue shadowed regions correspond to the areas swept by the asteroid for long times (the integrations are long enough so that longer integrations sweep the same regions). There, we observe how (153) Hilda, with the largest $C$ value of the six considered asteroids, describes a rounded trajectory around the Sun, meanwhile (2483) Guinevere, with the lowest $C$ value of the six, clearly approaches and loops around the locations of $L_3$, $L_4$ and $L_5$.

These different behaviours according to the different values of $C$ resembles to those exhibited by the periodic orbits of the family previously described and shown in Figures~\ref{fig:OPs} and \ref{fig:OPsHildas}.
The blue shadowed areas swept by the asteroids for long times, suggest that these asteroids orbit around the periodic orbits shown in Figure~\ref{fig:OPsHildas}. Then, the movement of Hilda asteroids in the planar CRTBP seems to be governed by the family of quasi-periodic orbits (invariant tori of dimension higher than 1) that exists around the family of periodic orbits previously described.

\subsection{Quasi-periodic orbits}\label{subsec:Quasi_C}

To compute the families of two-dimensional invariant tori we use a Poincaré section so that the intersection of each torus is a one-dimensional invariant curve. Concretely, we choose $y=0$ as a plane of section and, as the periodic orbits (Figure~\ref{fig:trasHildas}) cut at least twice this plane, we choose the part of the section that is between the primaries, so that $\dot{y}<0$. Note that with this choice we avoid the problem of tangencies and loops cutting more than once that happens if the section is chosen near $L_3$. Moreover, the system has a first integral (the Jacobi constant) so we use it as a parameter of the Poincar\'e map to lower the dimension of the section. Therefore, the Poincar\'e map $\bm{P}_C:\Sigma\rightarrow\Sigma$ is defined as follows. For each $C$ and each $(x_0,\dot{x}_0)\in\Sigma$, we obtain the only value $\dot{y}_0<0$ such that $(x_0,0,\dot{x}_0,\dot{y}_0)$ has Jacobi constant $C$. Then, we integrate the orbit that starts at this point till it cuts again $y=0$ for an $x_1$ value between that of the primaries (it happens that then $\dot{y}_1<0$). Finally we take the values $(x_1,\dot{x}_1)\in\Sigma$ as the image of the Poincar\'e map. An extra advantage of this approach is that, as we will see, we can use the symmetry \eqref{eq:sym} of the problem.

Given a fixed value $C$ for the Jacobi constant,
each invariant curve in $\Sigma$, denoted by $\bm{\varphi}: \mathbb{T} \rightarrow \mathbb{R}^2$, is characterized by its rotation number $\rho$. These curves must satisfy the invariance condition,
\begin{equation}\label{eq:inveq}
\bm{P}_C(\bm{\varphi}(\theta)) = \bm{\varphi}(\theta + \rho) \quad \text{for all} \quad \theta \in [0, 2\pi).
\end{equation}
Each curve is approximated as a truncated real Fourier series,
\[
\bm{\varphi}(\theta) \approx \bm{a}_0 + \sum_{k=1}^N \bm{a}_k \cos(k \theta) + \bm{b}_k \sin(k \theta),
\]
where $\bm{a}_0$, $\bm{a}_k$, and $\bm{b}_k \in \mathbb{R}^2$ are the Fourier coefficients with $k = 1, \ldots, N$ and $\theta \in [0, 2\pi)$.
At this point, we can impose the symmetry \eqref{eq:sym} of the problem, which implies that the $x$ coordinate is an even function and the $\dot{x}$ is an odd function. Therefore, we can reduce the number of unknowns,
\[
\begin{aligned}
    x(\theta) &= \varphi_1(\theta) = a_0 + \sum_{k=1}^N a_k \cos(k \theta), \\
    \dot{x}(\theta) &= \varphi_2(\theta) = \sum_{k=1}^N b_k \sin(k \theta).
\end{aligned}
\]

The final number of unknowns is $2N + 2$, since $\rho$ is not known. These unknowns can be solved by imposing the invariance condition \eqref{eq:inveq} using a Newton method that involves a linear system.

For tori close to the periodic orbit this procedure is extremely fast, because the number of Fourier nodes required (the value of $N$) is small and this allows to find the invariant curve along with its rotation number $\rho$. Therefore, this is what we have done for some of the tori close to the periodic orbit. However, as we move further away, the number of nodes increases, and resonances appear, which are not easy to handle.
To avoid these problems, a common procedure to  analyse families of invariant objects in a CRTBP is to perform Poincaré Section Plots (PSP) for a fixed value of the energy and a frequency analysis to find the rotation number. As we will see, the computation of tori will be unavoidable for the Elliptic RTBP.

This section is devoted to analyse the existence of a family of quasi-periodic solutions around each of the six selected asteroids, shown in Figure~\ref{fig:trasHildas}, by studying a PSP at each of their energy levels.
The spatial Poincaré section that we have used is the one ($\Sigma$) introduced above.
Then, for a given asteroid, we propagate long time its initial conditions and plot all the intersections to section $\Sigma$. 
Next, we modify its coordinates such that the energy level, or $C$ value, remains constant and again look for the crosses of the new trajectory to the defined section $\Sigma$. 

\begin{figure}[!ht]
    \centering
    \includegraphics[width=0.48\columnwidth]{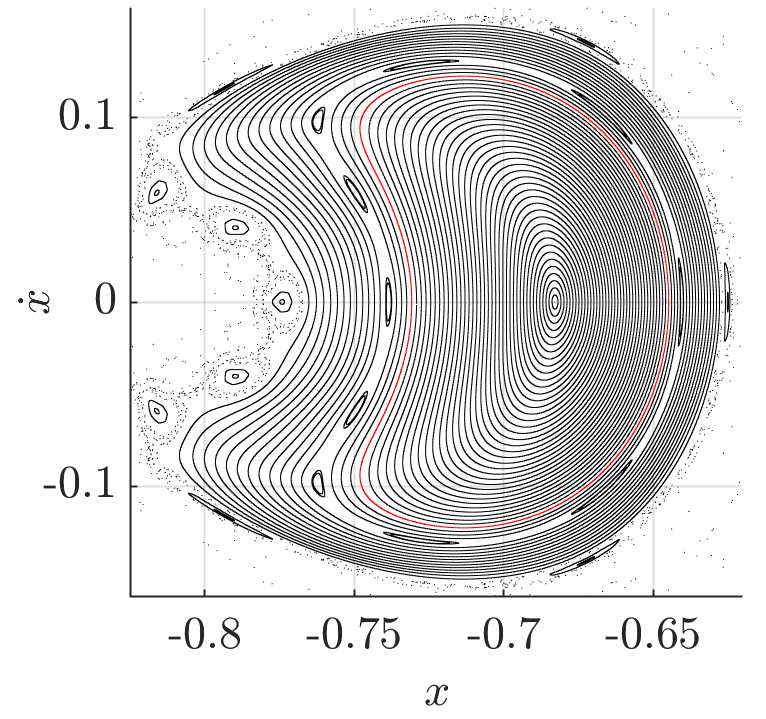} 
    \includegraphics[width=0.48\columnwidth]{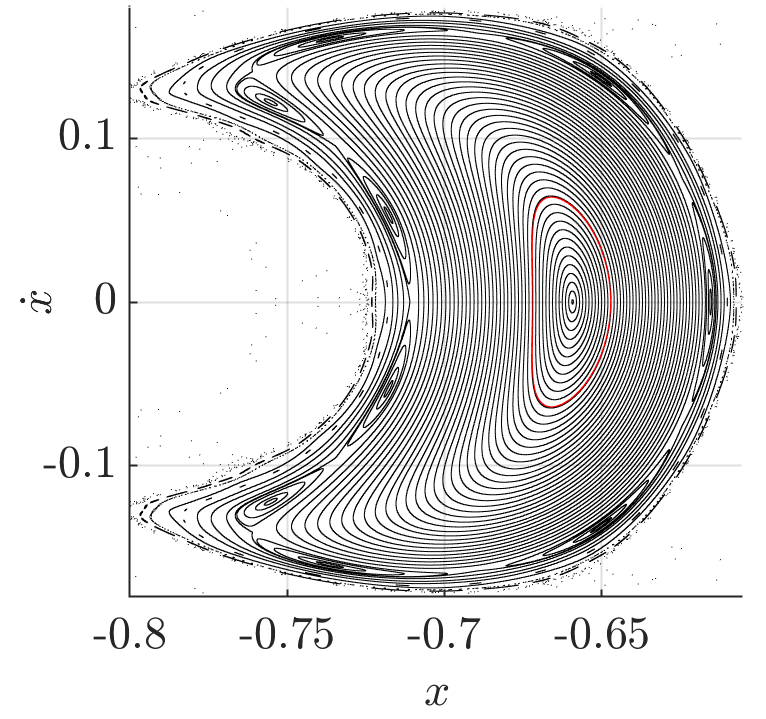} 
    \includegraphics[width=0.48\columnwidth]{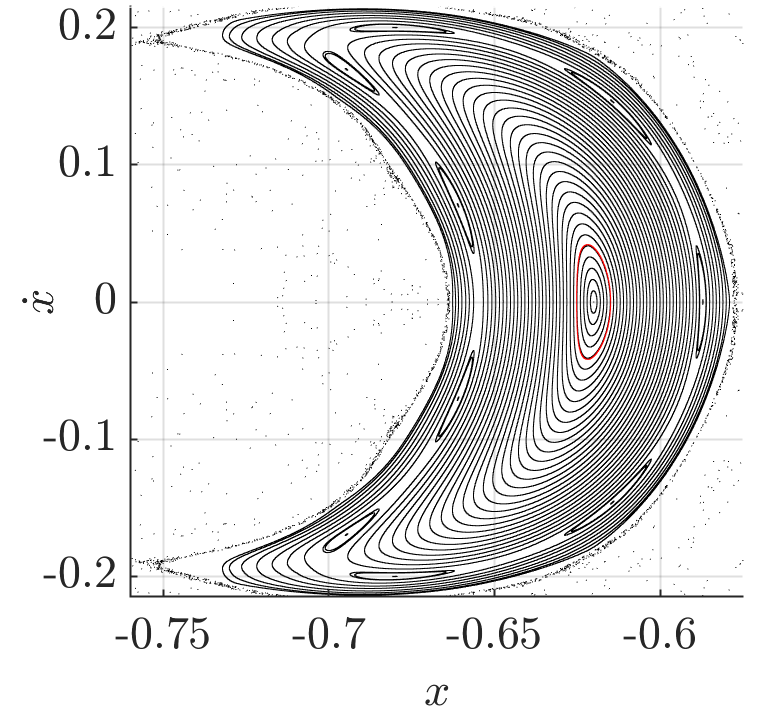} 
    \includegraphics[width=0.48\columnwidth]{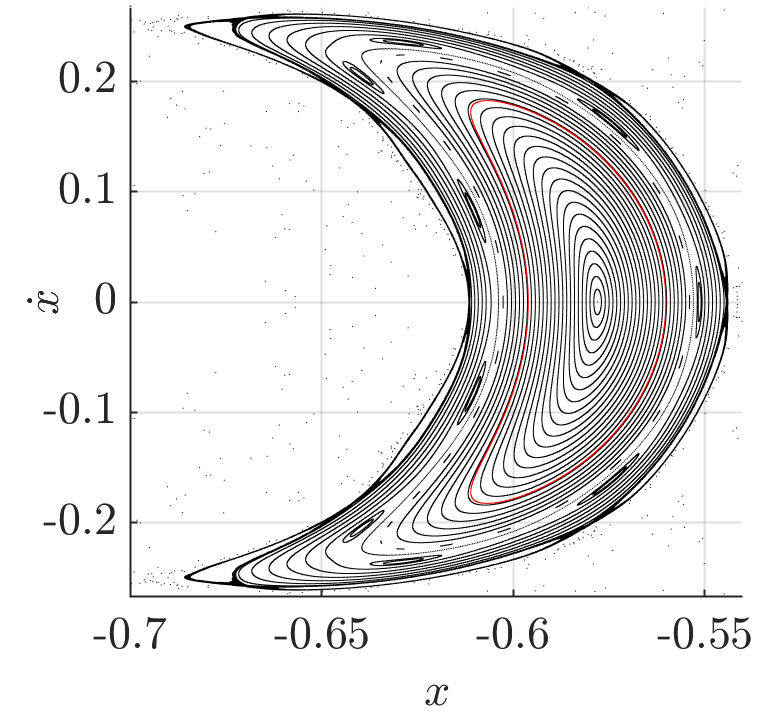} 
    \includegraphics[width=0.48\columnwidth]{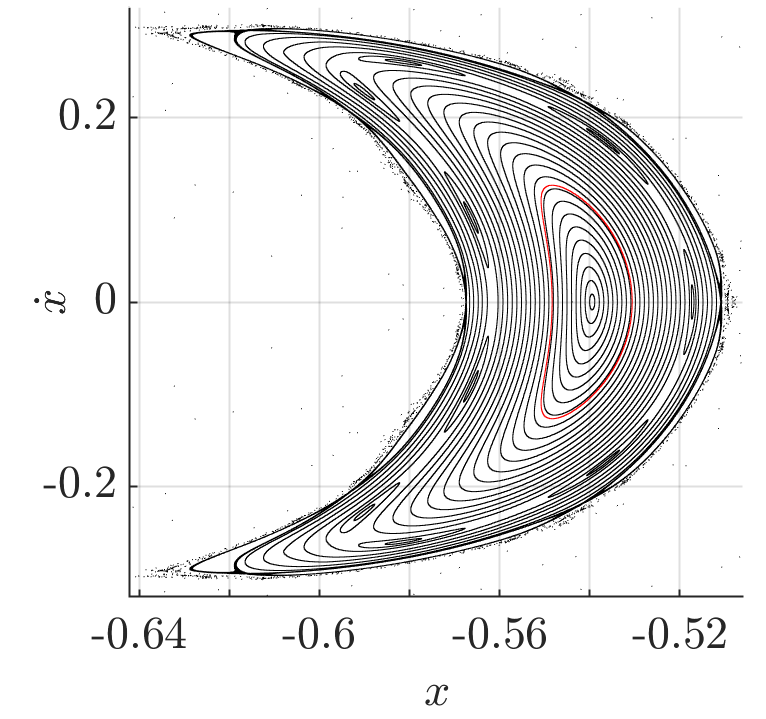} 
    \includegraphics[width=0.48\columnwidth]{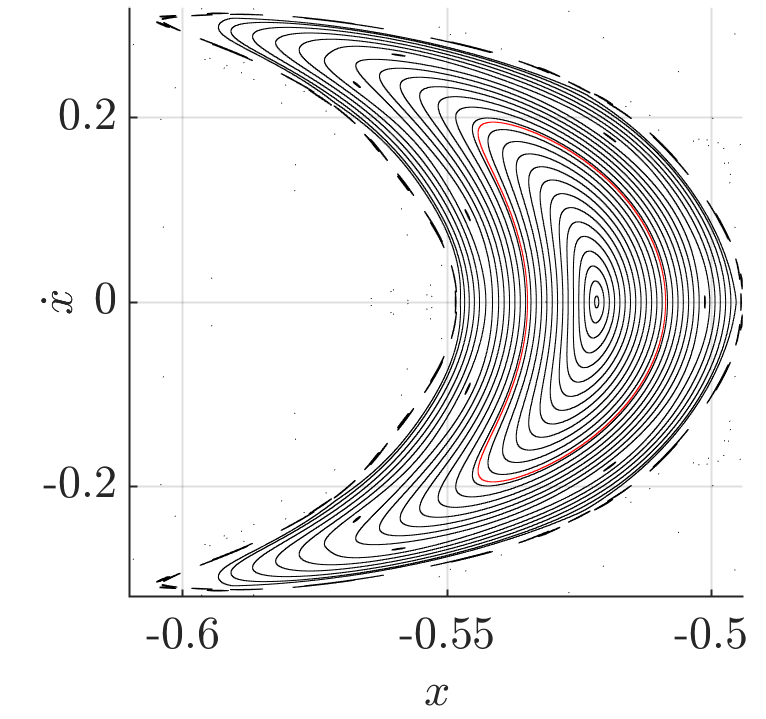} 
    \caption{Poincaré Section Plots (PSP) at the energy level of the six selected Hilda-type asteroids: (153) Hilda and (1911) Schubart in the first row, (499) Venusia and (1038) Tuckia in the second row, and (4446) Carolyn and (2483) Guinevere in the third row. In red the resulting closed curve for the initial conditions of each of these asteroids.}
    \label{fig:PSPast}
\end{figure}

In Figure~\ref{fig:PSPast}, the PSP for section $\Sigma$
at the energy levels of the six selected Hilda-type asteroids are shown in descending order of $C$ value. The curves in red correspond to intersections of the trajectories of the asteroids with section~$\Sigma$ and the black ones to intersections of the trajectories in the same energy level. We observe how each asteroid describes one closed curve in the ($x$,$\dot{x}$) plane, what indicates that their motion is trapped in a two-dimensional quasi-periodic solution of the planar CRTBP. Also, we find that each of those solutions is surrounded by other concentric closed curves, or chain of closed curves in the same PSP. This is sometimes known as islands that are surrounded by a ``chaotic sea''.
This means that the quasi-periodic solution found for each asteroid belongs to a family of two-dimensional quasi-periodic solutions that grow from a stable periodic orbit, that in the PSP is seen as a point in the ``centre'' of these families.

In a PSP, resonances are usually seen as chains of islands, where the number of islands on each chain gives the order of the resonance. For example, in the PSP of (153) Hilda, the first one in Figure~\ref{fig:PSPast}, we find some sets of 10 islands, which means resonances of order 10. Meanwhile, in the next PSPs in the same figure, for the other five selected asteroids, we find resonances of order 9.

Our goal is to make a general analysis for the Hilda group of asteroids and not only for the six selected above. However, since the number of asteroids in the Hilda group, described in Section~\ref{sec:selec}, is pretty large and each one has its specific energy level, it is not be feasible to perform a PSP for each of them. Therefore, we aim to make an analysis grouping them by narrow ranges of energy, or similarly, by narrow ranges of the Jacobi constant. 

In particular, we take the ranges in values of $C$ centered at the values corresponding to the same six Hilda asteroids shown in Figure~\ref{fig:PSPast}, and select for each of them the 40-50 asteroids closest in energy. All those asteroids in each group are analysed together in a plot that is similar to a Poincaré Section Plot (PSP), but for a narrow range of values of $C$. We refer to this plot as Thick-Poincaré Section Plot (T-PSP).

\begin{figure}[!ht]
    \centering
    \includegraphics[width=0.48\columnwidth]{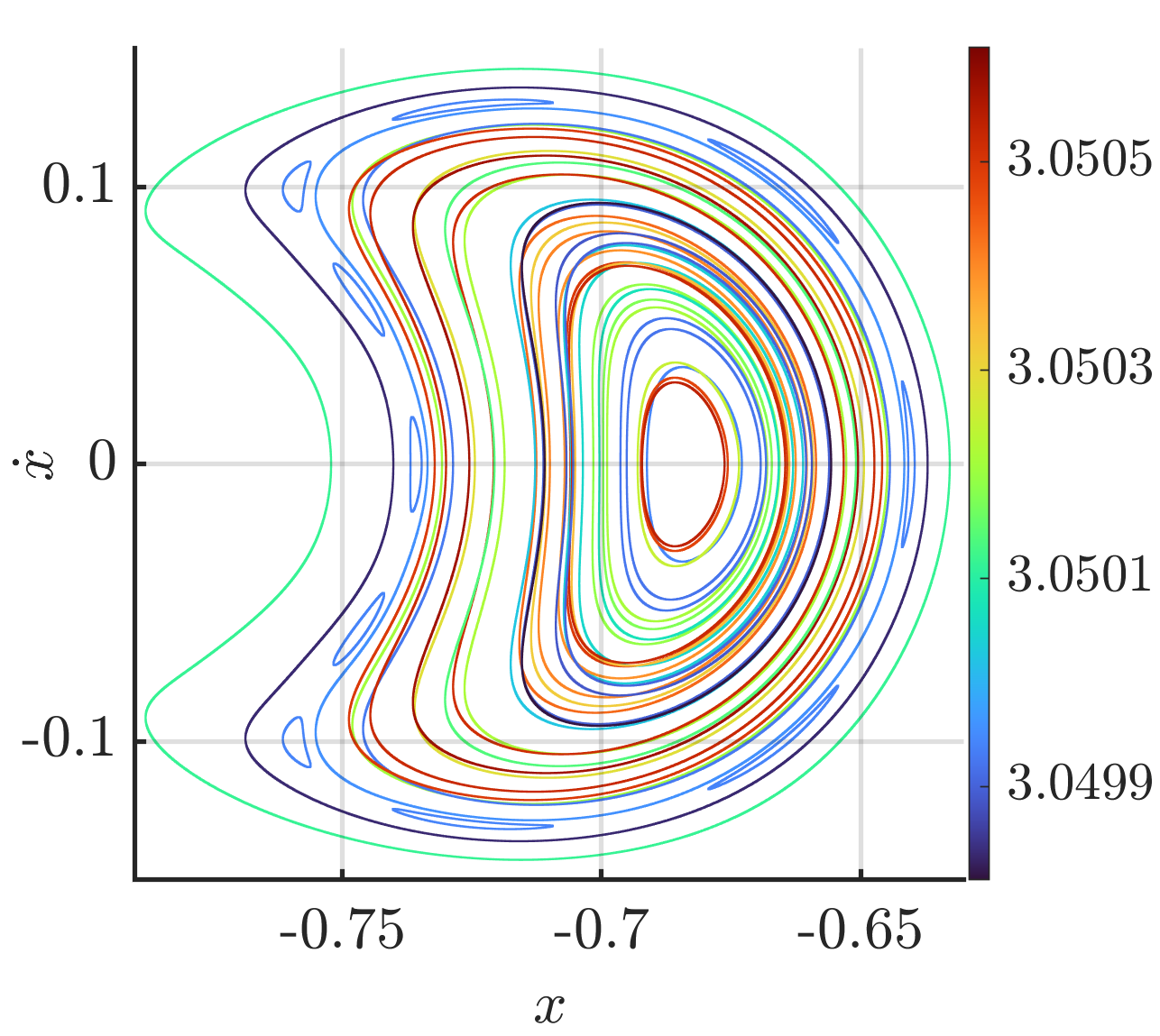} 
    \includegraphics[width=0.48\columnwidth]{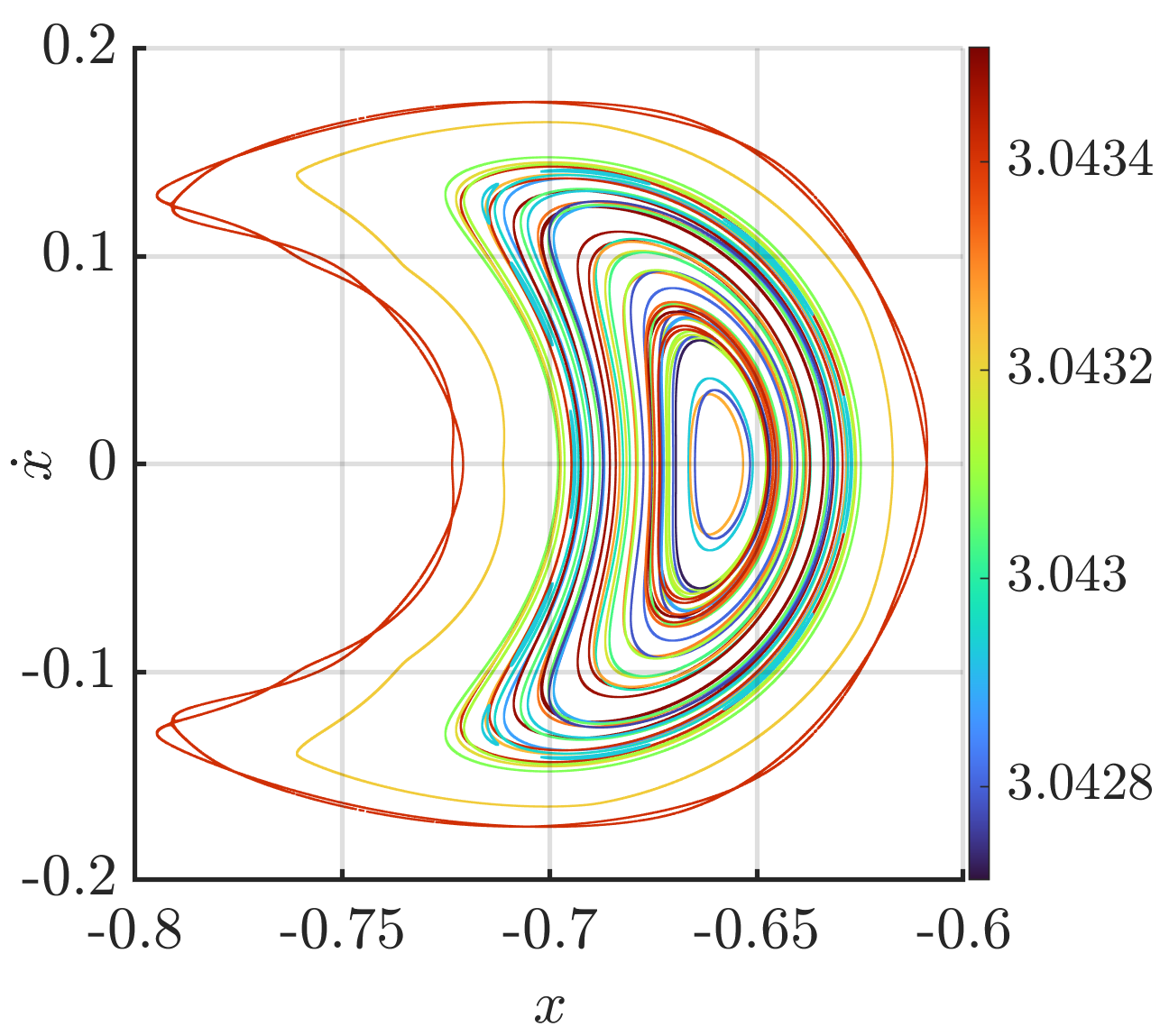} 
    \includegraphics[width=0.48\columnwidth]{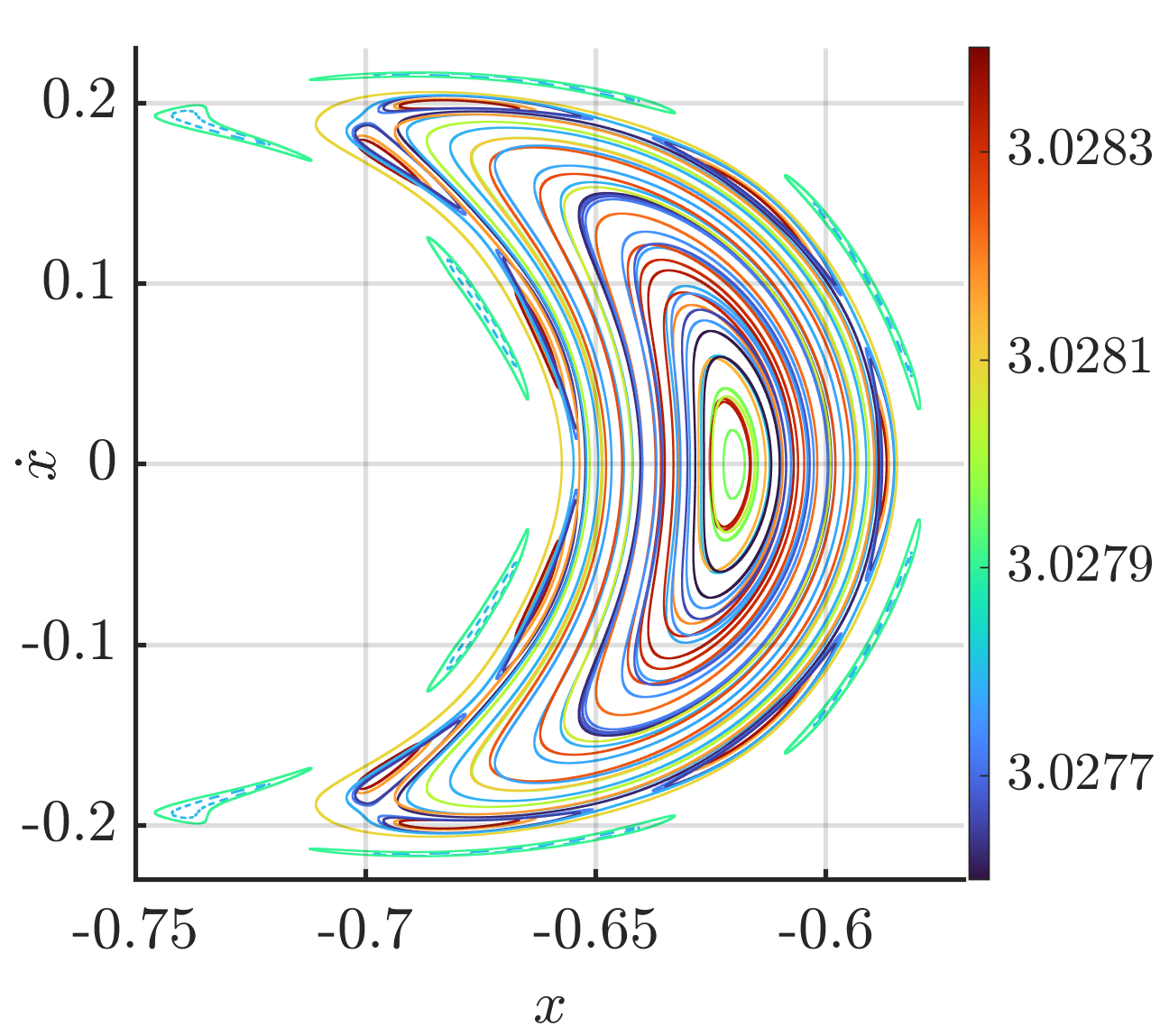} 
    \includegraphics[width=0.48\columnwidth]{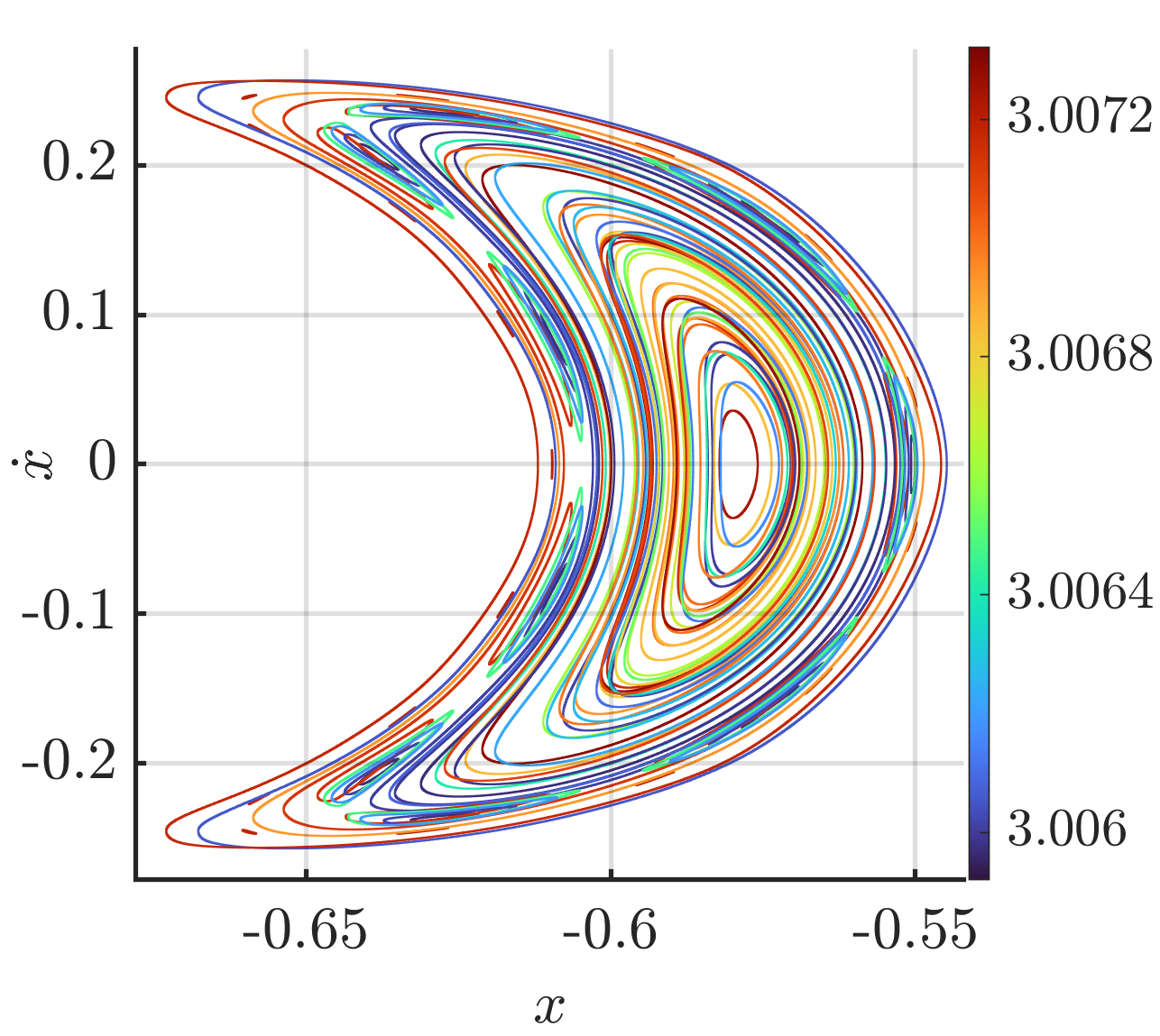} 
    \includegraphics[width=0.48\columnwidth]{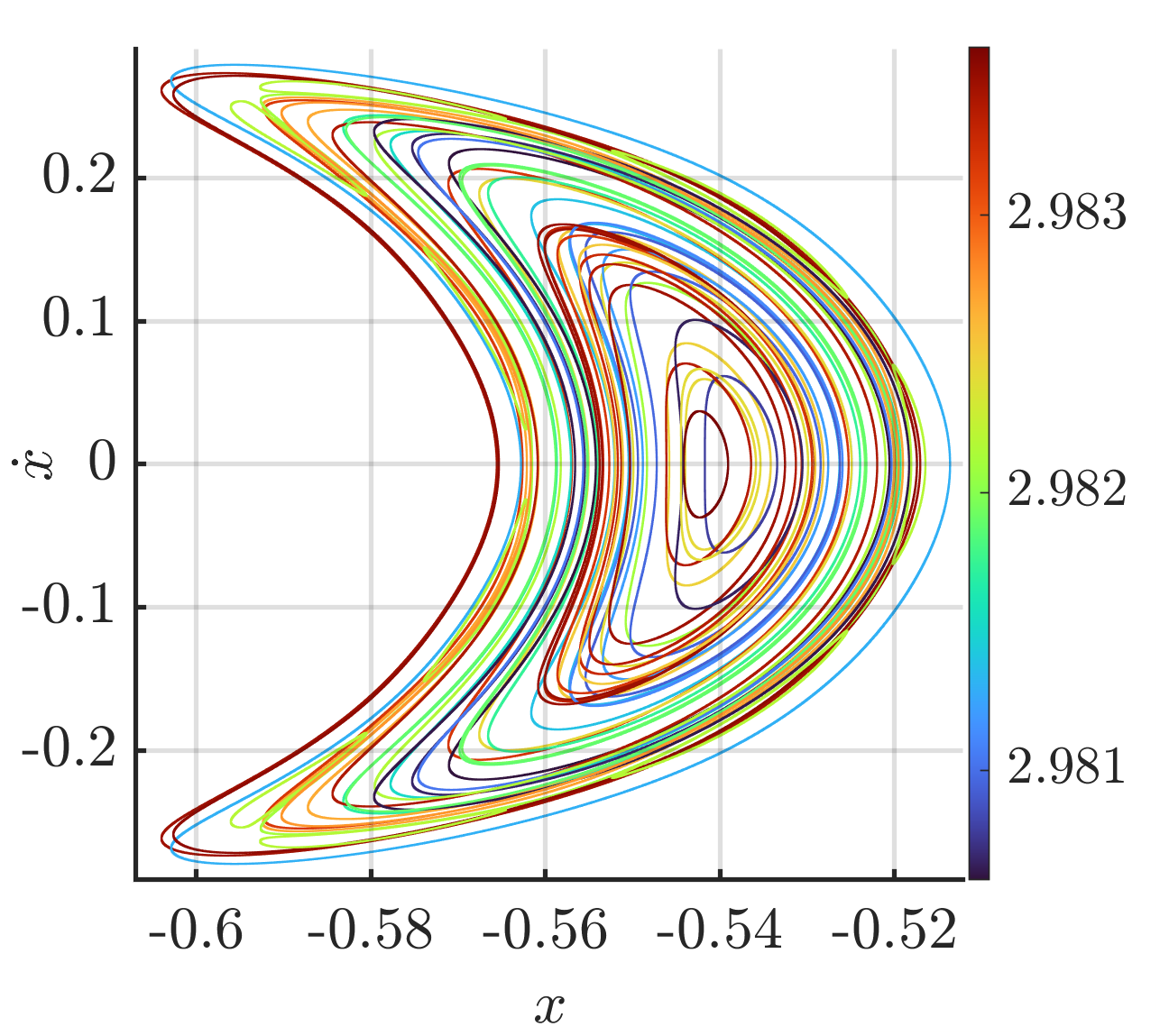} 
    \includegraphics[width=0.48\columnwidth]{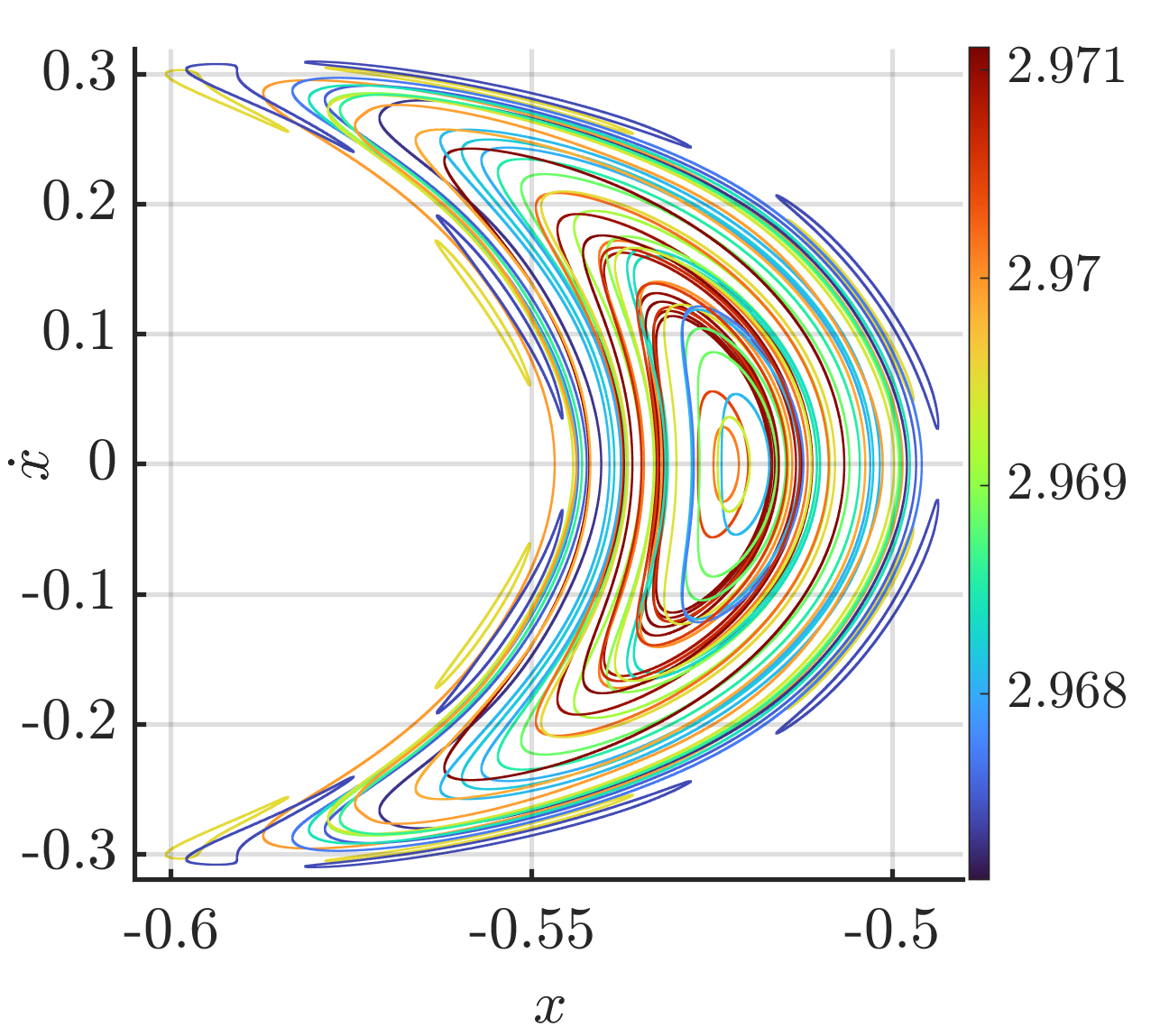} 
    \caption{Thick-Poincaré Section Plots (T-PSP) for Hilda-type asteroids grouped by narrow ranges of the Jacobi constant. Colour corresponds to $C$ value. See text for details.}
    \label{fig:TPSPast}
\end{figure}

Again, we use the section given by $\Sigma$ as the spatial Poincaré section, such that for each asteroid, we plot all the crosses of their trajectories to that section in a ($x$,$\dot{x}$) plane and colour them according to their $C$ value. Therefore, each colour corresponds to a different asteroid. In Figure~\ref{fig:TPSPast}, six of these T-PSP are presented. The groups of asteroids shown in these T-PSP have $C$ values centred at the $C$ value of the asteroids (153) Hilda and (1911) Schubart, in the first row, (499) Venusia and (1038) Tuckia, in the second row, and (4446) Carolyn and (2483) Guinevere, in the third row.

 Analysing the T-PSPs, we observe that most of the asteroids describe one closed curve on the ($x$,$\dot{x}$) plane, more or less concentric. Notice that now each curve in the T-PSP corresponds to a different asteroid, with slightly different energies among them, then it is expected the non-exactly concentric aspect in the T-PSPs when we plot them together. In addition to the asteroids that only describe one closed curve, we find that some of them describe sets of closed curves, denoting the presence of resonances. 

For example, in the first T-PSP of Figure~\ref{fig:TPSPast}, the one centred at the energy level of (153) Hilda, we observe a chain of 10 islands for one asteroid with $C$ value in blue. Note that, for the first PSP in Figure~\ref{fig:PSPast}, the one generated from the energy level of (153) Hilda, we found the same resonance. 
Likewise, looking at the other PSP and T-PSP we find the same correlations, with some small differences due to the slightly different energy levels of the asteroids in the T-PSP. For the case of (1911) Schubart the order of the resonance in the PSP was 9, meanwhile in the T-PSP we find an asteroid in a resonance of order 10, in turquoise colour. Both, PSP and T-PSP performed for (499) Venusia, show the resonance of order 9, and a resonance of 16 islands inside the one of order 9 is also found for the case of the T-PSP. Again, the T-PSP centered at (1038) Tuckia and (4446) Carolyn energy levels, show asteroids in resonances of order 9, as the ones observed in the PSP at their exact energy levels. Finally, in the T-PSP for (2483) Guinevere a resonance of order 8 is found, meanwhile in the PSP we found one of order 9.

The fact that, when performing the Thick-Poincaré Section Plots of the trajectories of a vast number of Hilda asteroids we see closed curves supports the hypothesis that their movement is indeed governed by two-dimensional invariant tori. Also, the T-PSP's have allowed us to visualise that the behaviour of this vast number of Hilda asteroids is well represented by the PSP's of the six selected asteroids.

\begin{figure}[!ht]
    \centering
    \includegraphics[width=0.48\linewidth]{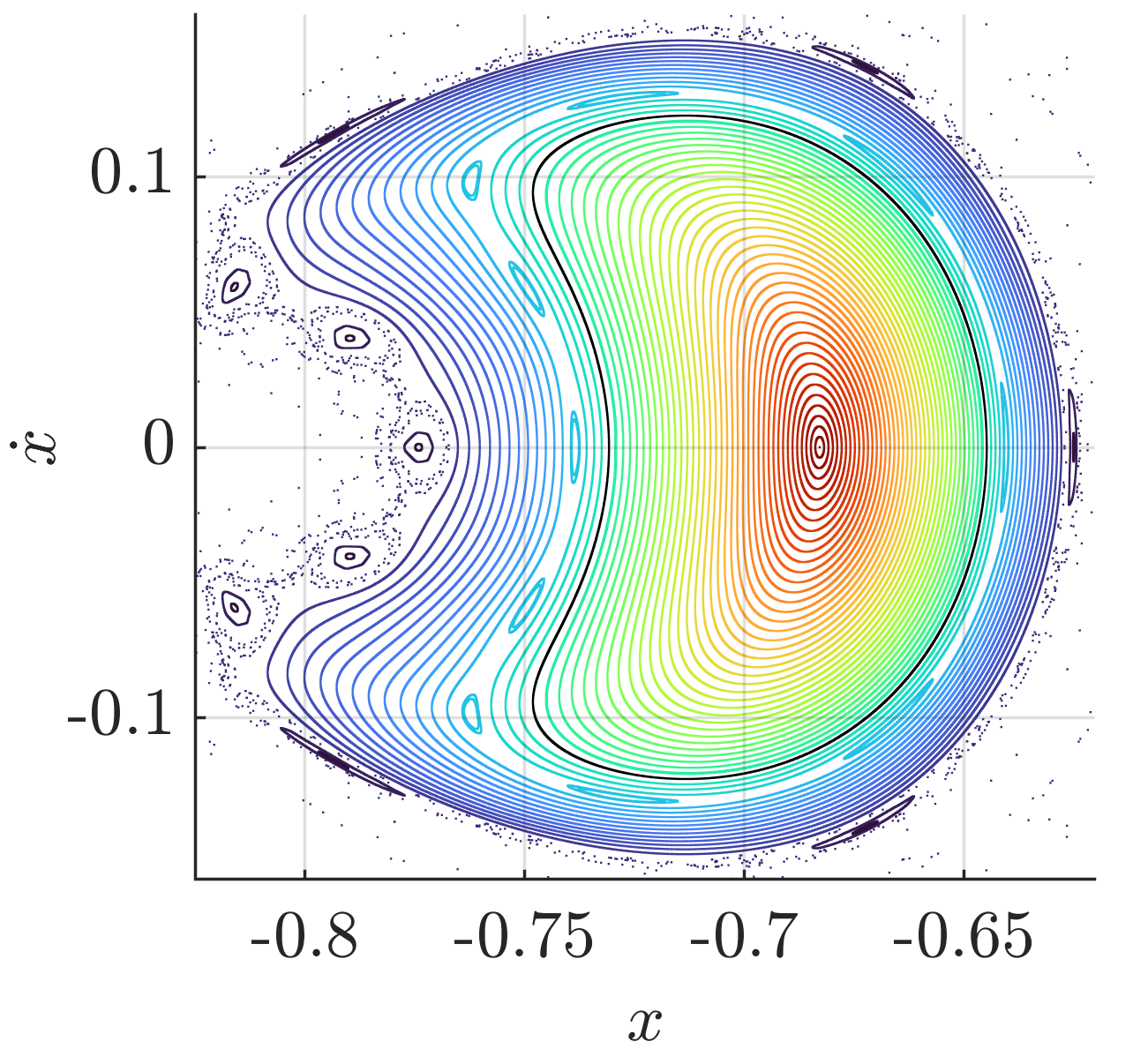}
    \includegraphics[width=0.48\linewidth]{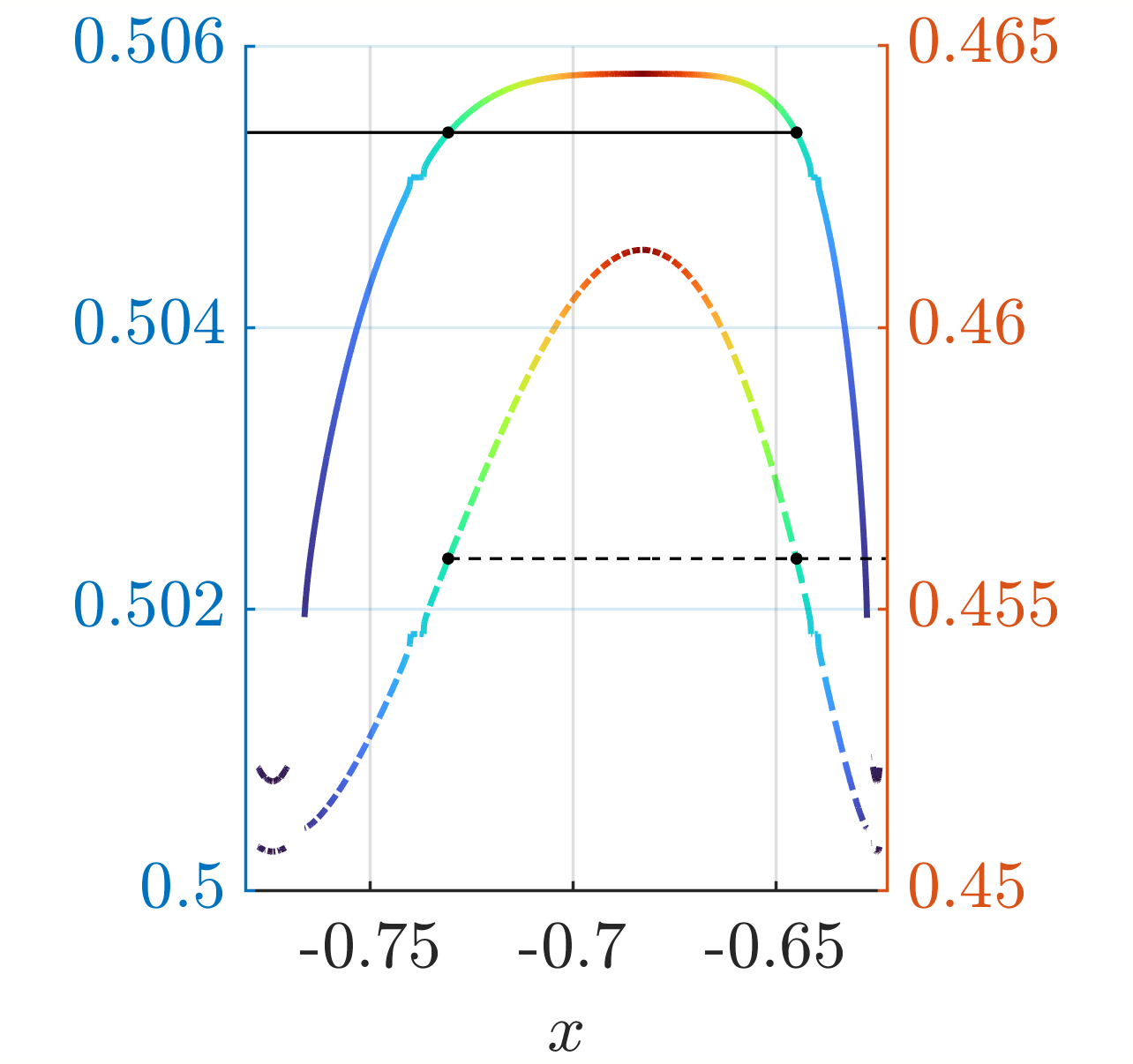}
    \includegraphics[width=0.48\linewidth]{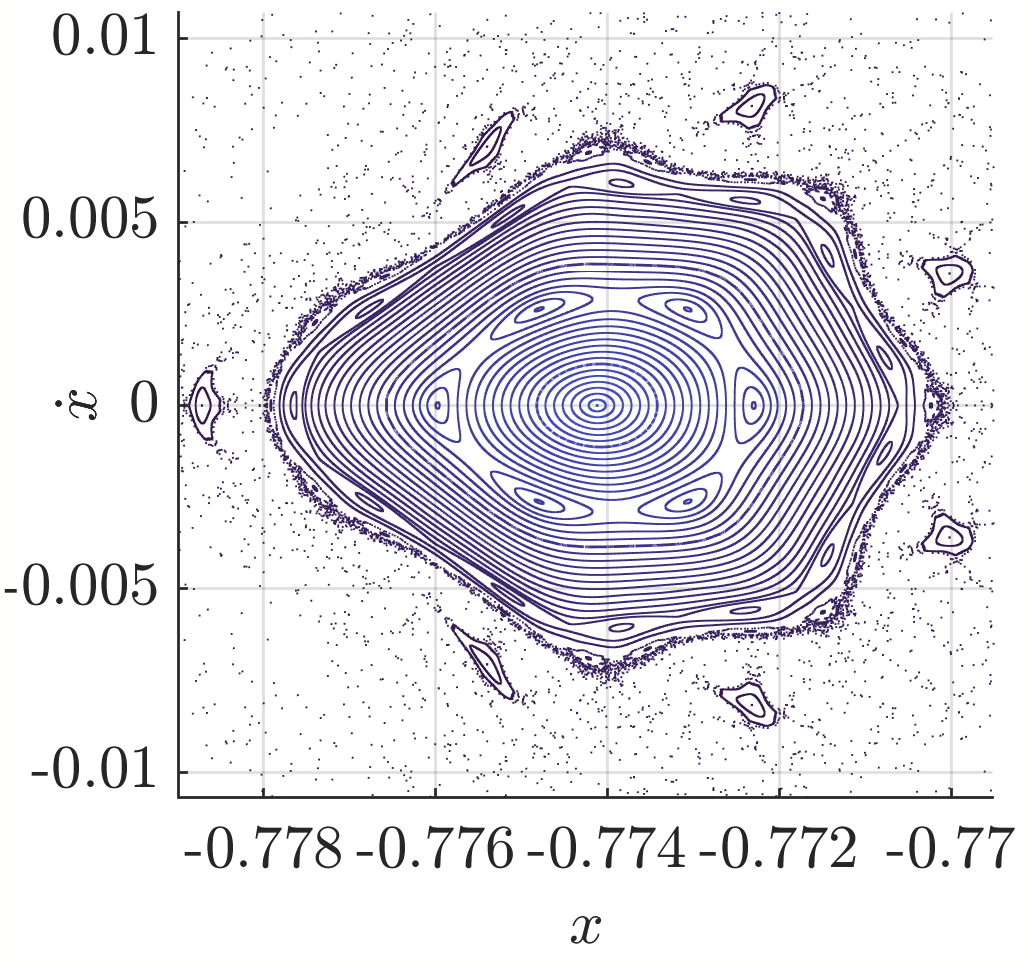}
    \includegraphics[width=0.48\linewidth]{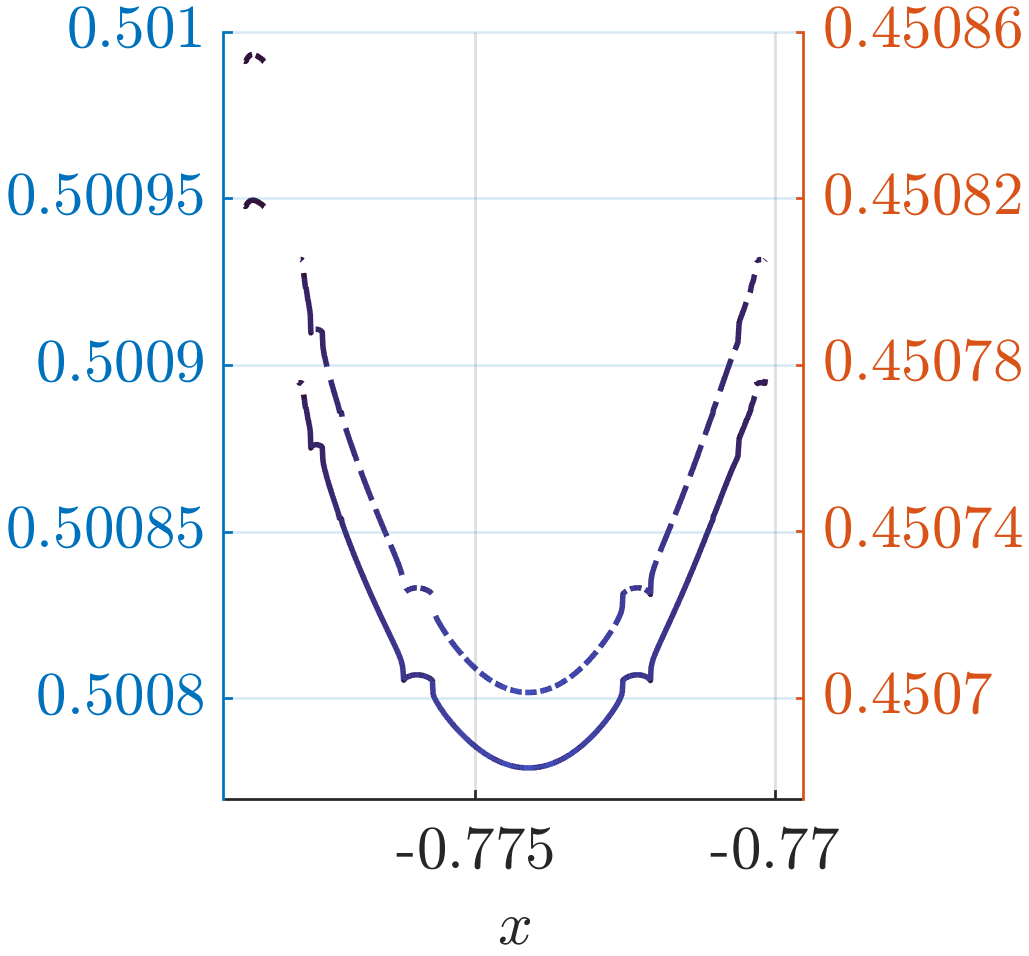}
    \caption{First row left, PSP at the energy level of asteroid (153) Hilda. Black curve corresponds the asteroid, while colour is used to identify each curve in the PSP (2D invariant tori of the flow) with their two main frequencies, included on the right. The value of these two main frequencies are shown as a function of $x$ when $\dot{x} = 0$.
    Second row, a zoom around one of the islands of the PSP at left, and the frequency analysis at right. See text for further details.}
    \label{fig:PSPfreqHilda}
\end{figure}

Therefore, the next step is to perform a frequency analysis of the two-dimensional quasi-periodic solutions shown in Figure~\ref{fig:PSPast}, since they were proved to be good representatives of the whole group of asteroids.
In Figure~\ref{fig:PSPfreqHilda} the family of 2D quasi-periodic solutions for (153) Hilda asteroid energy level is shown along with a zoom performed around one its islands, that contains other resonances itself. Now the curve corresponding to the asteroid is coloured in black and the others of the family have different colours to identify them with the two main frequencies shown on the right of the same figure. This plot shows the variation of the two main frequencies of the quasi-periodic solutions as they cut the $x$-axis, that is when the invariant curves intersect with $\dot{x} = 0$. Left vertical axis corresponds to the first main frequency, plotted in solid line, and right vertical axis corresponds to the second main frequency, plotted in dash line. The horizontal black line points the two frequencies of (153) Hilda asteroid. In order to compute these two frequencies for each solution, we have implemented the frequency analysis described in Appendix~\ref{sec:fran}. 

Analysing the upper right plot of Figure~\ref{fig:PSPfreqHilda}, we observe that both frequencies have values close to $0.5$, consistent with the fact that these asteroids are close to a $4\pi$ periodic orbit (Figure~\ref{fig:OPsHildas}) whose normal frequency is also close to $0.5$. It is also remarkable that the variation of these two frequencies along the family is limited to a narrow range of values, as shown in the two vertical axis of the frequencies plot. Resonances are also clearly identified in this plot, since the variation of the frequencies along the $x$-axis suffers some discontinuities and jumps. 

For example, a order 10 resonance, shown in turquoise colour in the PSP at the left of Figure~\ref{fig:PSPfreqHilda}, has two islands cutting the $x$-axis at $x\approx-0.74$ and $x\approx -0.64$, what produces two small horizontal segments, at the same $x$ values in the frequency analysis shown on the right of that figure. If we compute the ratio between the two frequencies of these flat segments in the plot, we see it corresponds to a rational number, as it is expected for resonances. Similarly, at $x\approx-0.774$ (and $x\approx -0.6$) we observe a jump in the frequency analysis due to another order 10 resonance, the one coloured in dark blue. Second row of Figure~\ref{fig:PSPfreqHilda}, shows a zoom of an island of the last mentioned resonance on the left, along with the corresponding frequency analysis on the right. This island contains other resonances itself, that are also clearly identified in the frequency analysis. At the ``center'' of this island, $x\approx-0.774$, a point corresponding to a stable periodic orbit is found, surrounded by 2D quasi-periodic solutions until a six order resonance (6 islands) appears, producing two discontinuities in the frequency analysis of the right. If we continue the analysis in the outside direction from the center of the island we find another resonances. Three of them, of orders 13, 20 and 7 are easily to identify, as well as the discontinuities they produce in the frequency graph. However, as we refine the PSP in this small area we can find more resonances inside the island, that displays a really rich dynamics.

\begin{figure}[!ht]
    \centering
    \includegraphics[width=0.48\linewidth]{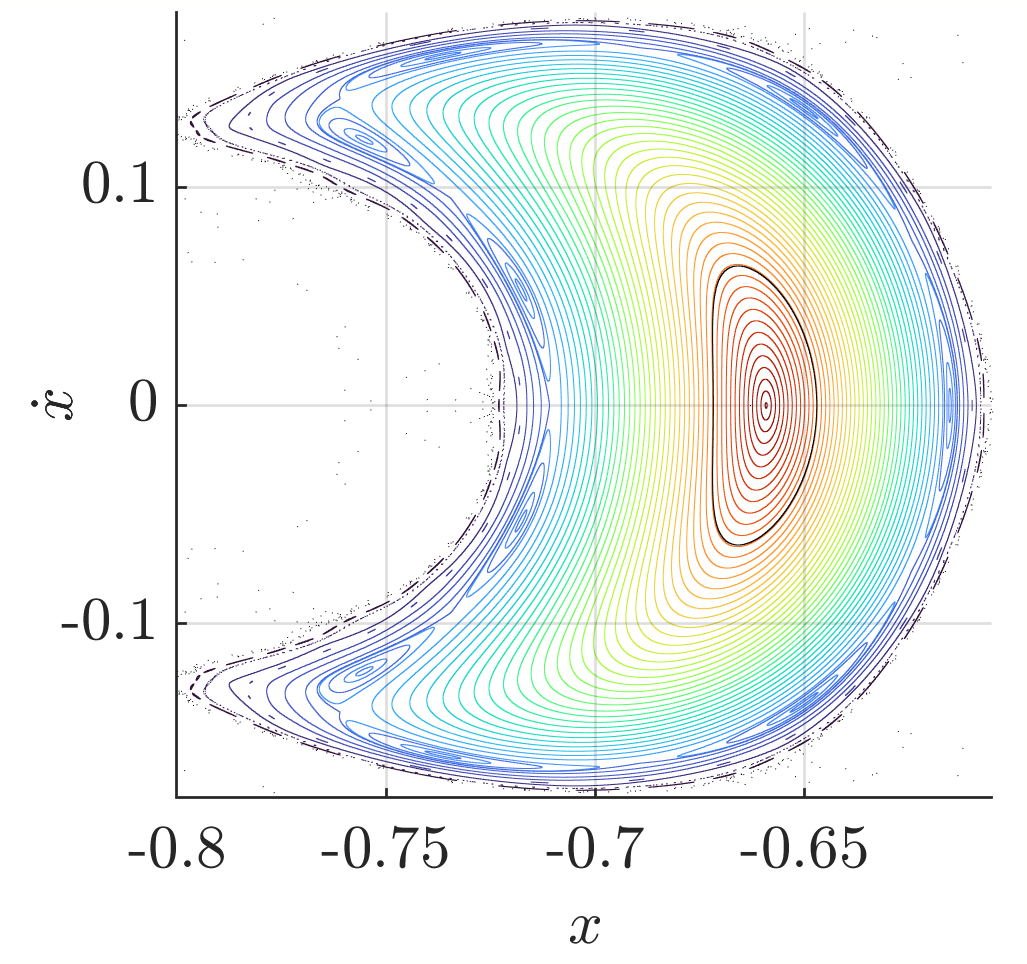} 
    \includegraphics[width=0.48\linewidth]{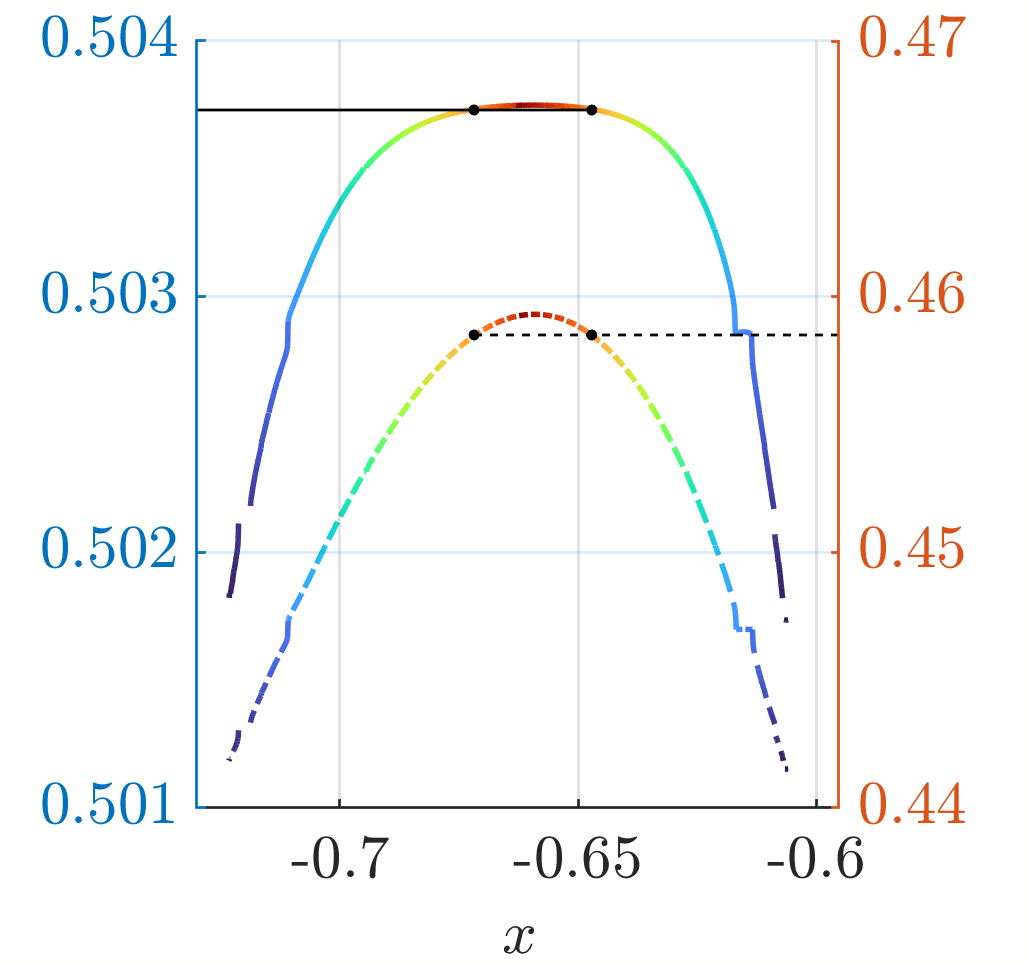}
    \includegraphics[width=0.48\linewidth]{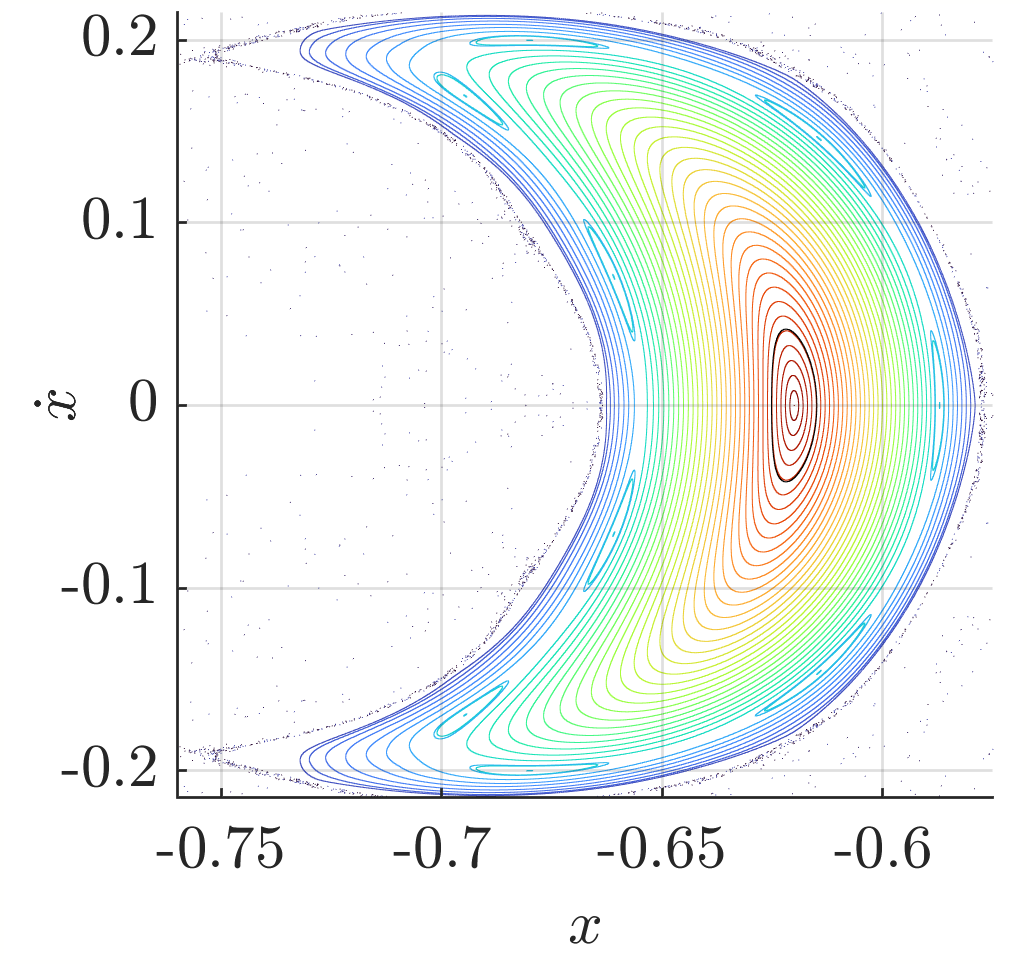} 
    \includegraphics[width=0.48\linewidth]{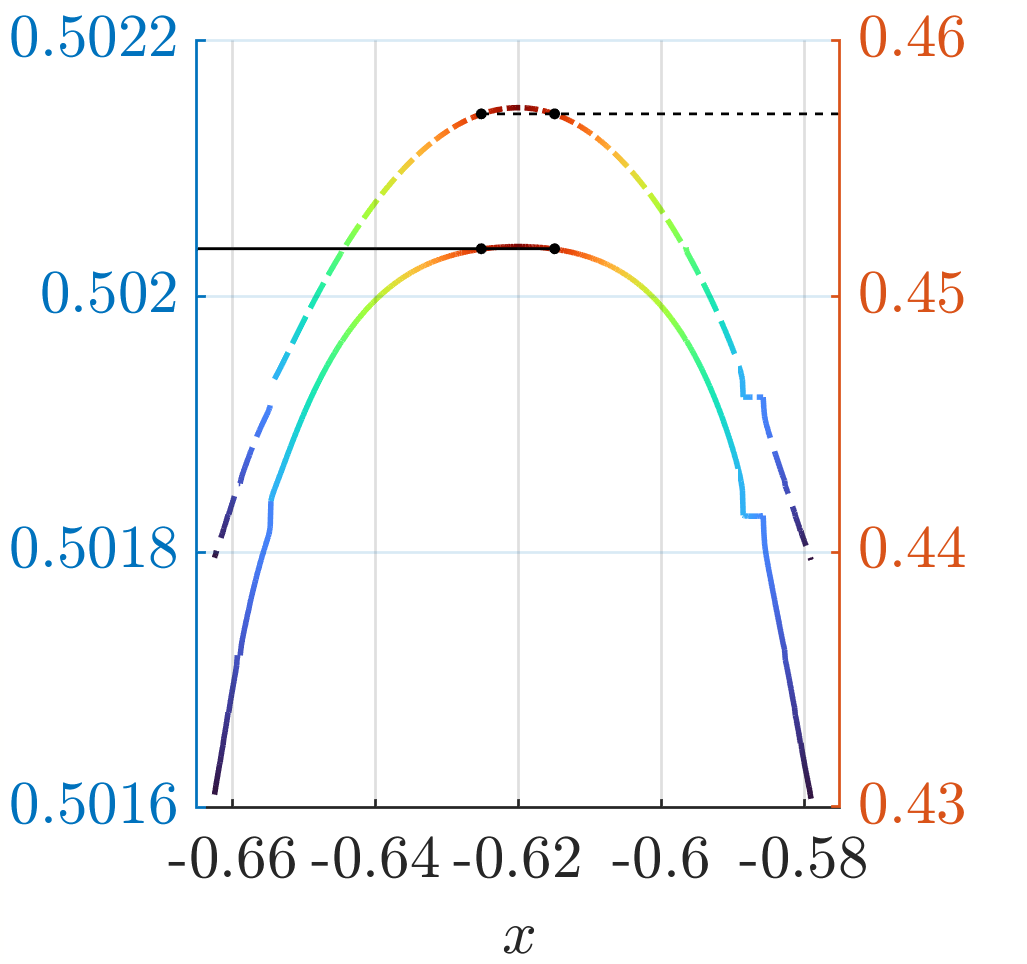}
    \includegraphics[width=0.48\linewidth]{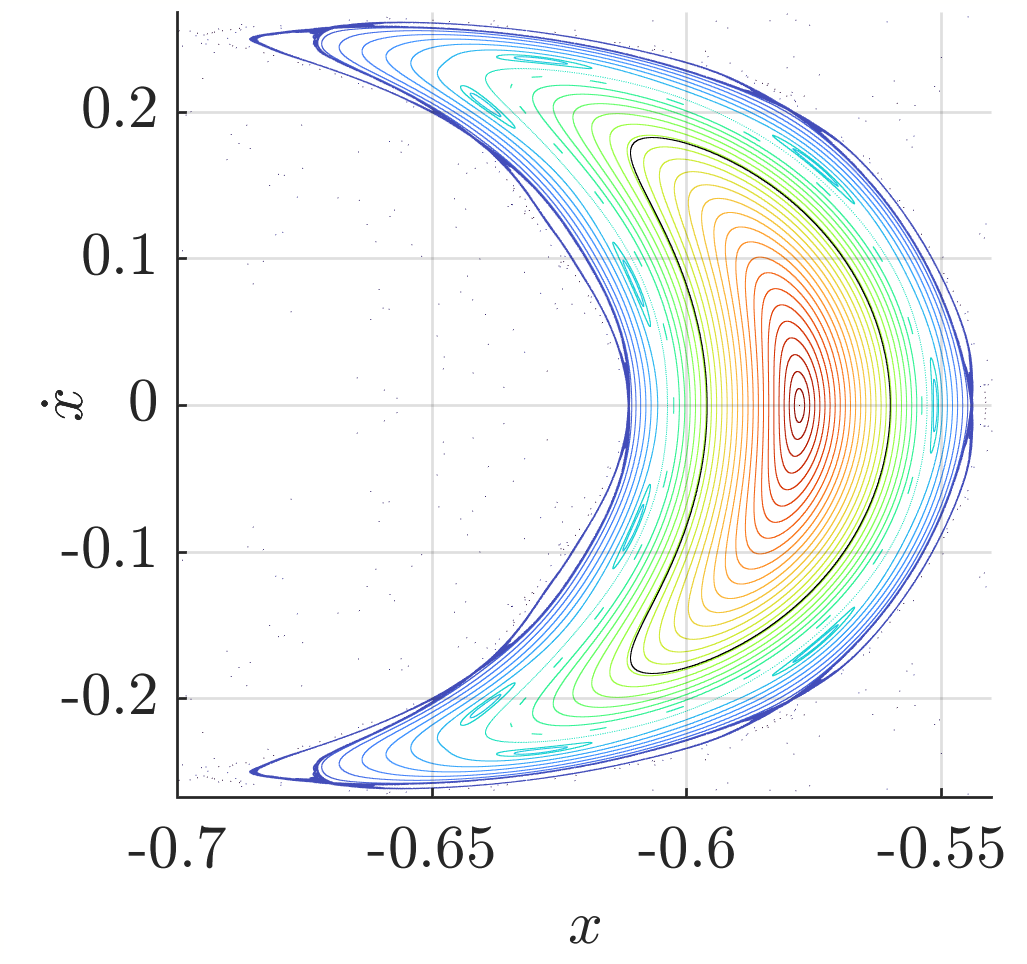} 
    \includegraphics[width=0.48\linewidth]{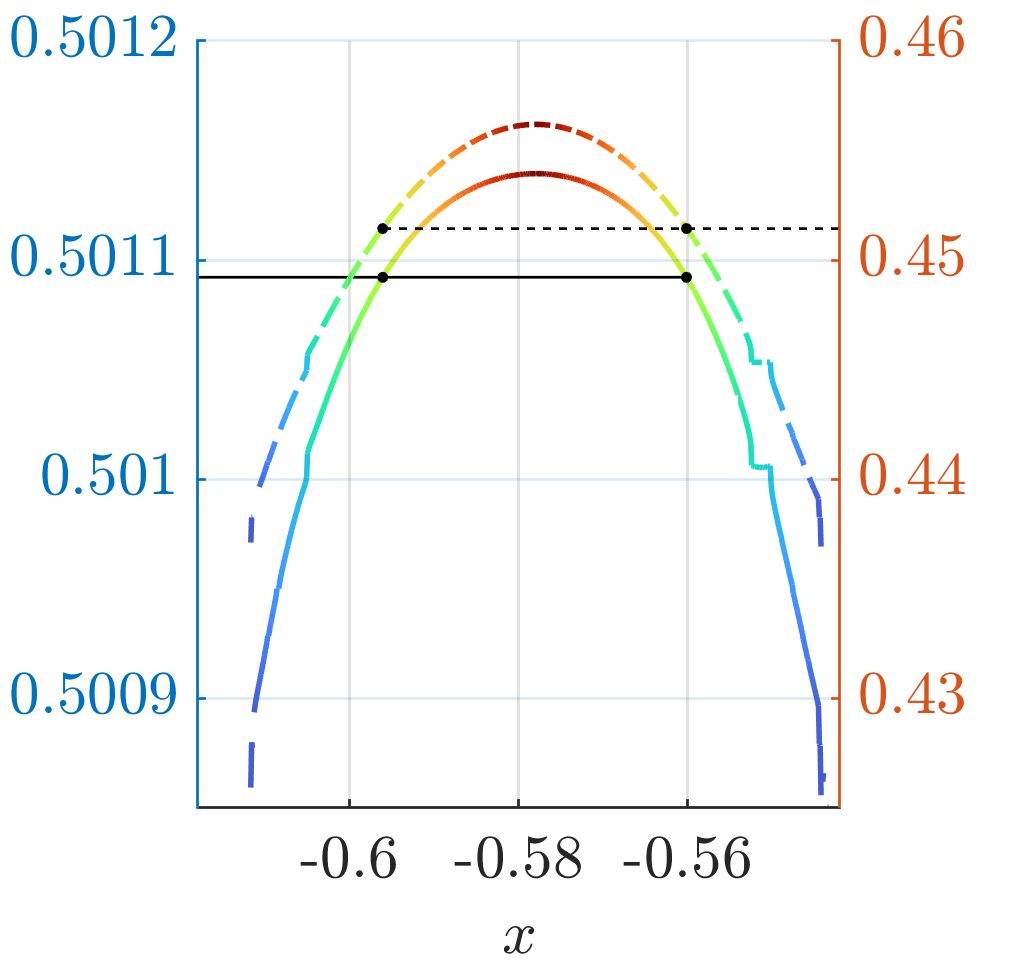}
    \includegraphics[width=0.48\linewidth]{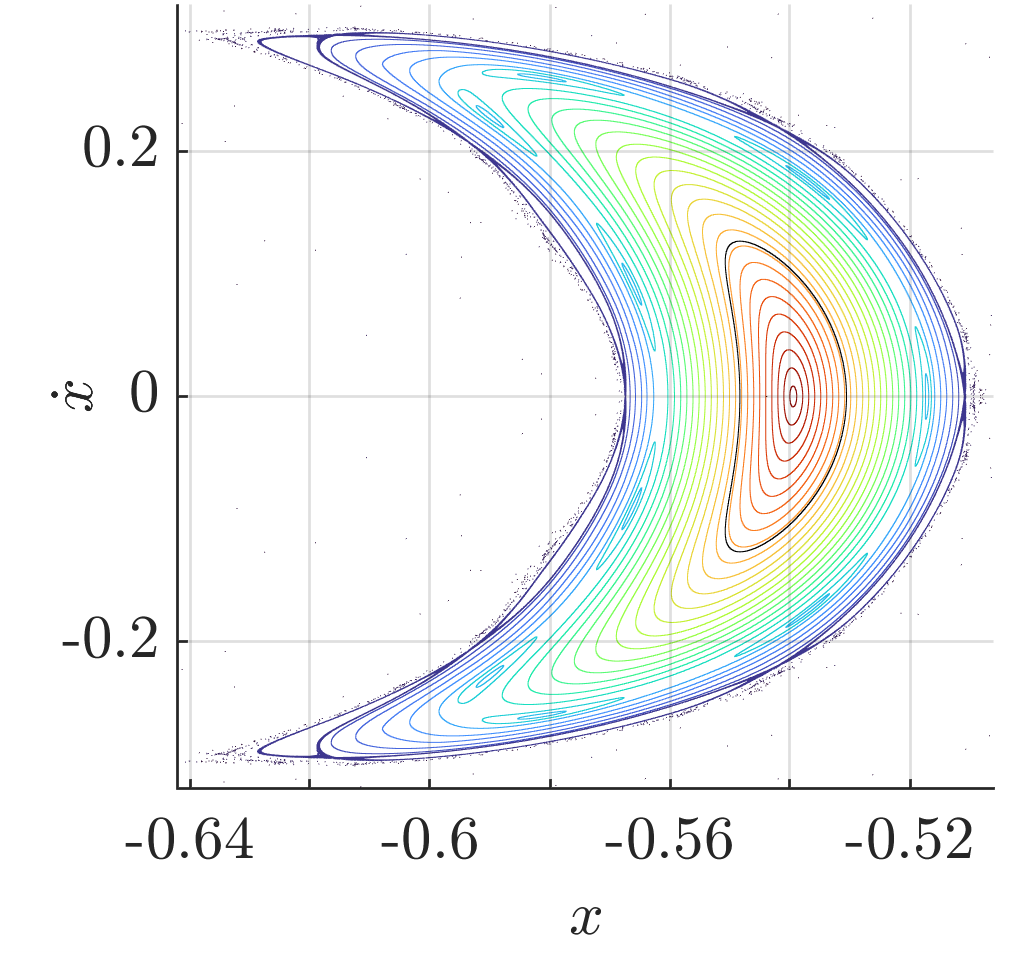} 
    \includegraphics[width=0.48\linewidth]{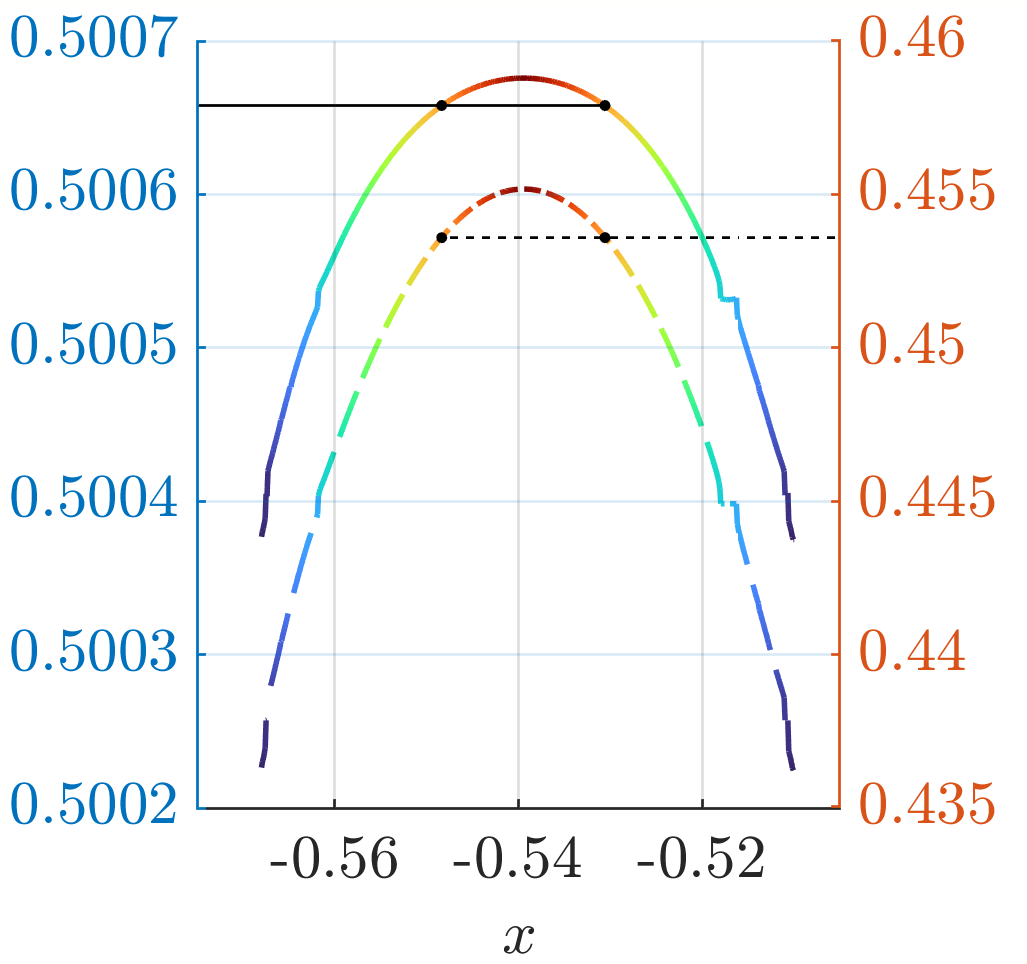}
    \includegraphics[width=0.48\linewidth]{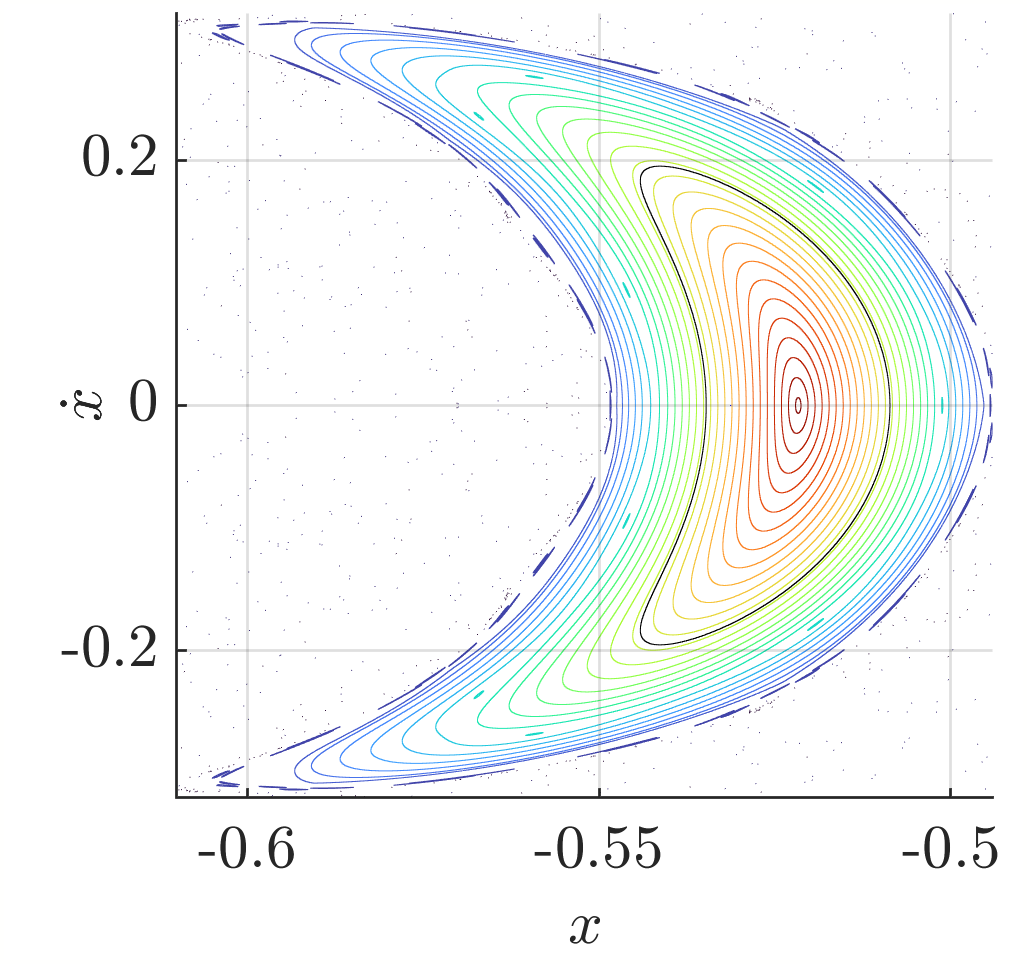} 
    \includegraphics[width=0.48\linewidth]{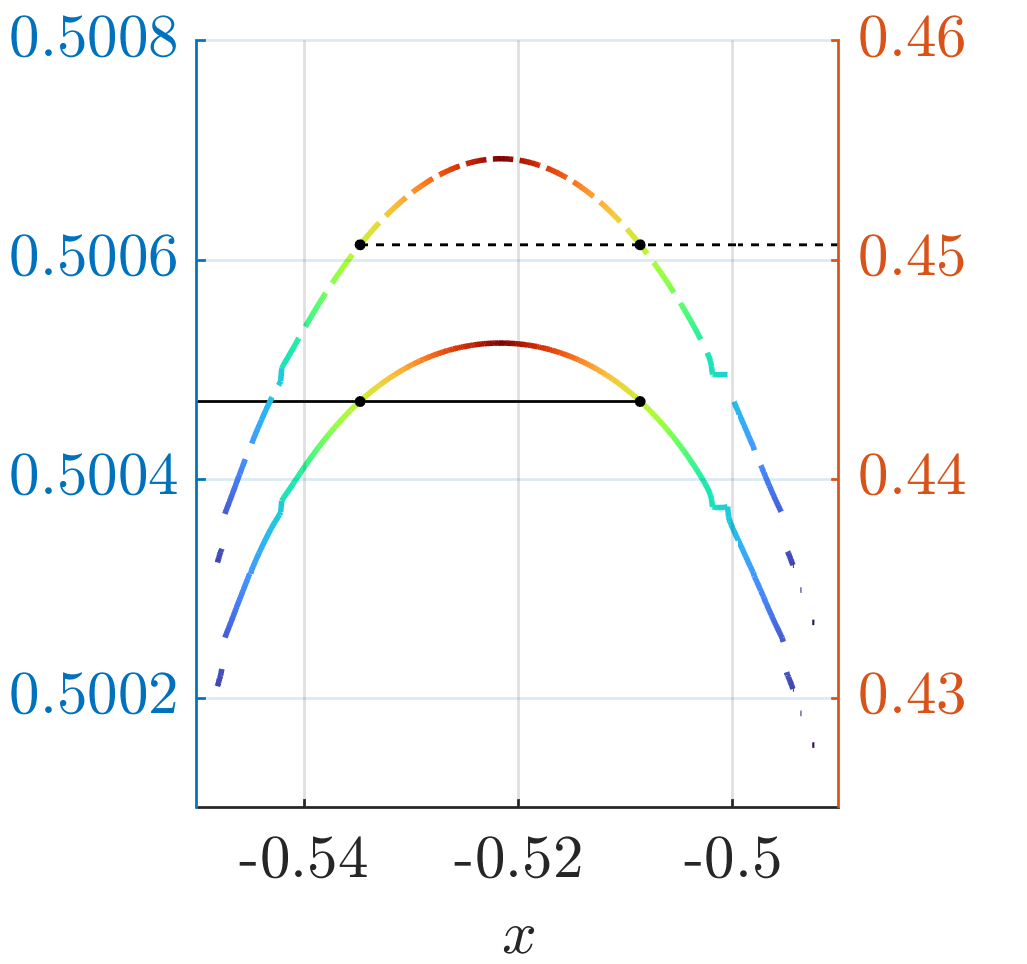}
    \caption{Analougous representation as in the first row of Figure \ref{fig:PSPfreqHilda} for the energy levels of asteroids (1911) Schubart, (499) Venusia, (1038) Tuckia, (4446) Carolyn and (2483) Guinevere (from top to bottom).}
    \label{fig:PSPfreqOther}
\end{figure}

Figure~\ref{fig:PSPfreqOther} shows the same analysis previously described for the other five selected asteroids. First column shows in black the intersections with section $\Sigma$ defined in Section~\ref{subsec:Quasi_C} of the trajectories of asteroids
(1911) Schubart, (499) Venusia, (1038) Tuckia, (4446) Carolyn and (2483) Guinevere, from up to down.
The other coloured points correspond to the families of two dimensional quasi-periodic solutions at the same energy level as each asteroid. At the right of each PSP, second column, the variation of the two main frequencies of the quasi-periodic solutions is shown as they cut the $x$-axis. Again, left vertical axis corresponds to the first main frequency, plotted in solid line, right vertical axis corresponds to the second main frequency, plotted in dash line, and the horizontal black lines points the frequencies of the asteroids.

The same comments as before about the small ranges of the frequencies and their closeness to value $0.5$ are applicable in these new plots for the other 5 representative asteroids. Again, resonances produce some horizontal segments and discontinuities in the frequency analysis, as it was expected.

Therefore, the main conclusion of this section devoted to the analysis in the planar CRTBP is that the motion of the Hilda group of asteroids is confined in an island of two-dimensional quasi-periodic solutions.

\section{Analysis in the planar ERTBP}\label{sec:ERTBP}
Here we focus on the dynamics of the Hilda-type asteroids in the planar Elliptical Restricted-Three Body Problem,
described in Section~\ref{sec:elliptic}.
We can apply directly the change of coordinates described in Section~\ref{sec:CC} to our selection of Hilda-type asteroids, and project them in the horizontal plane $z=z'=0$.

The eccentricity of the ERTBP acts as a $2\pi$ time-periodic perturbation of the CRTBP. Then, under generic non-degeneracy and non-resonance conditions, the family of periodic orbits of Section~\ref{sec:periodic} (Figures~\ref{fig:OPs} and \ref{fig:OPsHildas}) becomes a family of two dimensional invariant tori, adding the frequency of the eccentricity (i.e., 1) to the one of the family. Due to the resonances, this family has a Cantor structure where the holes that form the Cantor structure are exponentially small with the order of the resonance \cite{JorbaVillanueva1997}. For this reason, only low order resonances can be observed in a double precision calculation.
We compute these quasi-periodic solutions in Section~\ref{subsec:E_2D}. Analogously, the two dimensional invariant tori studied in Section~\ref{subsec:Quasi_C} for the CRTBP (Figures~\ref{fig:PSPast} and \ref{fig:TPSPast}) become three dimensional invariant tori in the ERTBP, that are studied in Section~\ref{subsec:E_3D}, along with the frequency analysis as we did in the circular case.

\subsection{Two-dimensional quasi-periodic solutions}\label{subsec:E_2D}

The first quasi-periodic solutions mentioned, 2D invariant tori, are defined by two frequencies, one that comes from the periodic orbit of the CRTBP and other that corresponds to the time-periodic perturbation. Therefore, we can use the period associated to one of these frequencies to define a temporal Poincaré map in which the 2D invariant tori of the flow is seen as a 1D invariant curve of the map. Typically, the period of the perturbation is used to define the temporal Poincaré map, since all the quasi-periodic solutions of the family share this frequency and hence, it can be used to study the whole family in the same stroboscopic map. 

We recall that in the ERTBP, the temporal magnitude is the true anomaly, denoted by $f$, and the perturbation is given by the eccentricity of the orbits of the primaries acting with periodicity in $f$ of $2\pi$. Therefore the stroboscopic map that we use, denoted by $P$, is defined as a temporal map corresponding to values of the true anomaly equal to integer multiples of $2\pi$.

The computation and stability analysis of this kind of solutions is a covered topic and there is not a unique way of doing so. In particular, we follow the procedures presented in \cite{CastellaJorba2000, Jorba2001, GabernJorba2004}, since these do not require any particular property of the invariant curve. 

Similarly as the computation of invariant curves in spatial maps summarised in Section~\ref{subsec:Quasi_C}, the computation of the invariant curves of the temporal map is based on looking for a representation of the invariant curve, $\bm{\varphi}: \mathbb{T} \rightarrow \mathbb{R}^4$, that satisfies the invariance condition
\begin{equation}
\bm{P}(\bm{\varphi}(\theta))=\bm{\varphi}(\theta+ \rho) \ \ \theta \in \mathbb{T},    
\label{eq:invcond}
\end{equation}
where $\rho\in\mathbb{T}$ is the rotation number of the curve. In this case, the rotation number is known, since it corresponds to the relation between the period of the periodic orbit in the CRTBP and the period of the stroboscopic map, 
\begin{equation}
    \rho = 2\pi \frac{\omega}{\omega_{e}} =\frac{4 \pi^2}{T},
    \label{eq:rho}
\end{equation}
being $\omega$ and $T$ the frequency and period, respectively, of the periodic orbit in the CRTBP, and $\omega_{e}$ the natural frequency of the ERTBP, that in normalised units of the model is $\omega_{e}=1$.

The representation that we choose for the invariant curves is an approximation through a real truncated Fourier series, as we did in Section~\ref{subsec:Quasi_C}. Then we apply a Newton method to Equation~(\ref{eq:invcond}) in order to find the Fourier coefficients for the four dimensions of the phase space. Notice that, the symmetry \eqref{eq:symE} with respect to the $x$-axis could be used again 
to reduce the number of Fourier coefficients in a half.

Once we have the invariant curves, we analyse their stability by looking for pairs of eigenvalue and
eigenfunction $(\lambda, \bm{\psi}) \in \mathbb{C} \times (C(\mathbb{T},\mathbb{C}^4)/\{0\})$ that satisfies the following generalised eigenvalue problem,
\begin{equation}
    D_x \bm{P}(\bm{\varphi}(\theta)) \bm{\psi} (\theta) = \lambda \Gamma_\rho \bm{\psi}(\theta),
    \label{eq:GEV}
\end{equation}
where $\Gamma_\rho$ denotes the operator $\Gamma_\rho:C(\mathbb{T},\mathbb{C}^4) \rightarrow C(\mathbb{T},\mathbb{C}^4)$ such that $\Gamma_\rho \bm{\psi}(\theta) = \bm{\psi}(\theta + \rho)$. Again, Equation~(\ref{eq:GEV}) is solved by Newton method in terms of the Fourier coefficients.

In practice, we have performed a continuation of the periodic orbits in the CRTBP shown in Figure~\ref{fig:OPsHildas}, increasing little by little the eccentricity from zero (CRTBP) to the corresponding value in the Sun-Jupiter ERTBP, $e\approx 0.04869$. Notice that the periodic orbits are also seen as invariant curves of the stroboscopic map $P$.

\begin{figure}[!ht]
    \centering
    \includegraphics[width=0.48\linewidth]{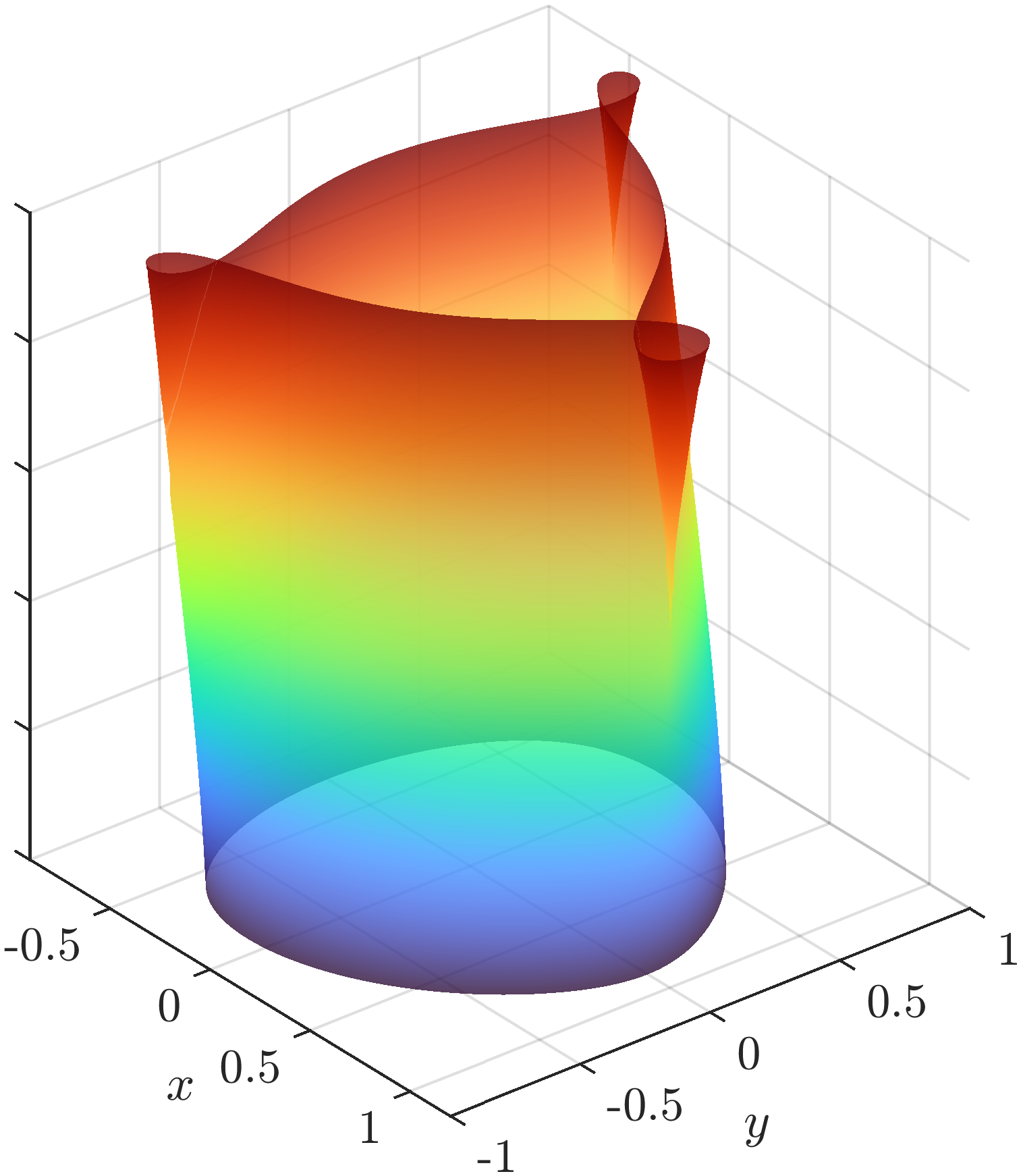}
    \includegraphics[width=0.48\linewidth]{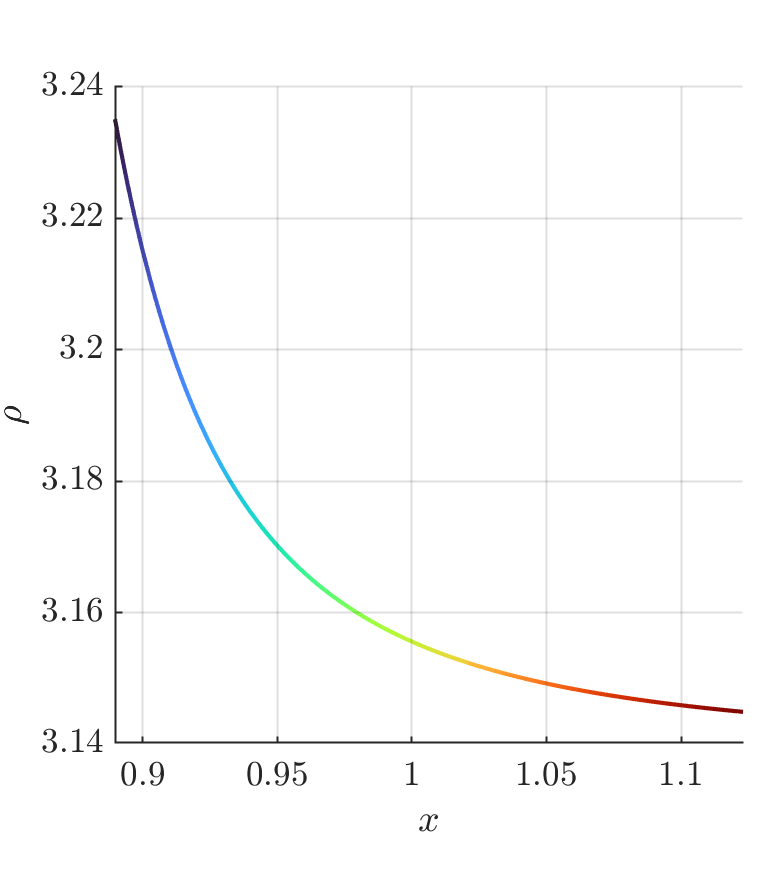}
    \caption{Left, surface composed by  an effective continuum of invariant curves of the stroboscopic map $\bm{P}$ in the ERTBP, computed from the periodic orbits of CRTBP related to Hilda behaviour.
    The colour is used to identify each invariant curve of the left with the point at the ($x$,$\rho$) plot in the right, that shows the variation of the rotation number according to the cut in $x$-axis of each invariant curve. }
    \label{fig:ICs}
\end{figure}

As mentioned in Section~\ref{sec:CRTBP}, the periodic orbits in the family related to Hilda behaviour are found for a range of values of Jacobi constant approximately among 2.98 and 3.06, or similarly, for a range of values of periods approximately $T\in(12.2,4\pi)$. Then, according to Equation~(\ref{eq:rho}), the selected periodic orbits, become invariant tori that are seen as invariant curves of the map $\bm{P}$ with rotation numbers $\rho \in (\pi, 3.2359)$ approximately, meanwhile non-resonance conditions are satisfied. Notice that the resulting family of invariant objects in the ERTBP is then a Cantorian family of 2D invariant tori. However, the gaps in the family are so small (below double precision resolution) that the family can be considered as effectively continuous \cite{JorbaVillanueva1997}. Figure~\ref{fig:ICs} left, shows a surface composed by the family of invariant curves of the map $\bm{P}$ are shown projected in the $xy-$plane, such that each transversal cut to this surface corresponds to a different invariant curve with its corresponding rotation number. Same figure right, shows the variation of the rotation number with the cut in the positive part of the $x$-axis. As the cut is given for higher values of $x$, the rotation number decreases, approaching the resonant value $\pi$.

In spite of the fact that the surface shown in Figure~\ref{fig:ICs} recalls to the surface shown in Figure~\ref{fig:OPsHildas}, it must be clear that the one in  Figure~\ref{fig:OPsHildas} corresponds to a family of periodic orbits of the flow of the CRTBP, meanwhile the one in Figure~\ref{fig:ICs} corresponds to a family of 2D invariant tori of the flow of the ERTBP, seen as a family of invariant curves of the stroboscopic map $\bm{P}$. That means that each transversal cut of the surface in Figure~\ref{fig:ICs} is an invariant curve of $\bm{P}$, i.e. a temporal section of the 2D invariant tori at $f=0$ or $f=2\pi$. 

\begin{figure}[ht!]
    \centering
    \begin{minipage}[t]{0.49\linewidth}
    \centering
        \includegraphics[width=1\textwidth,valign=t]{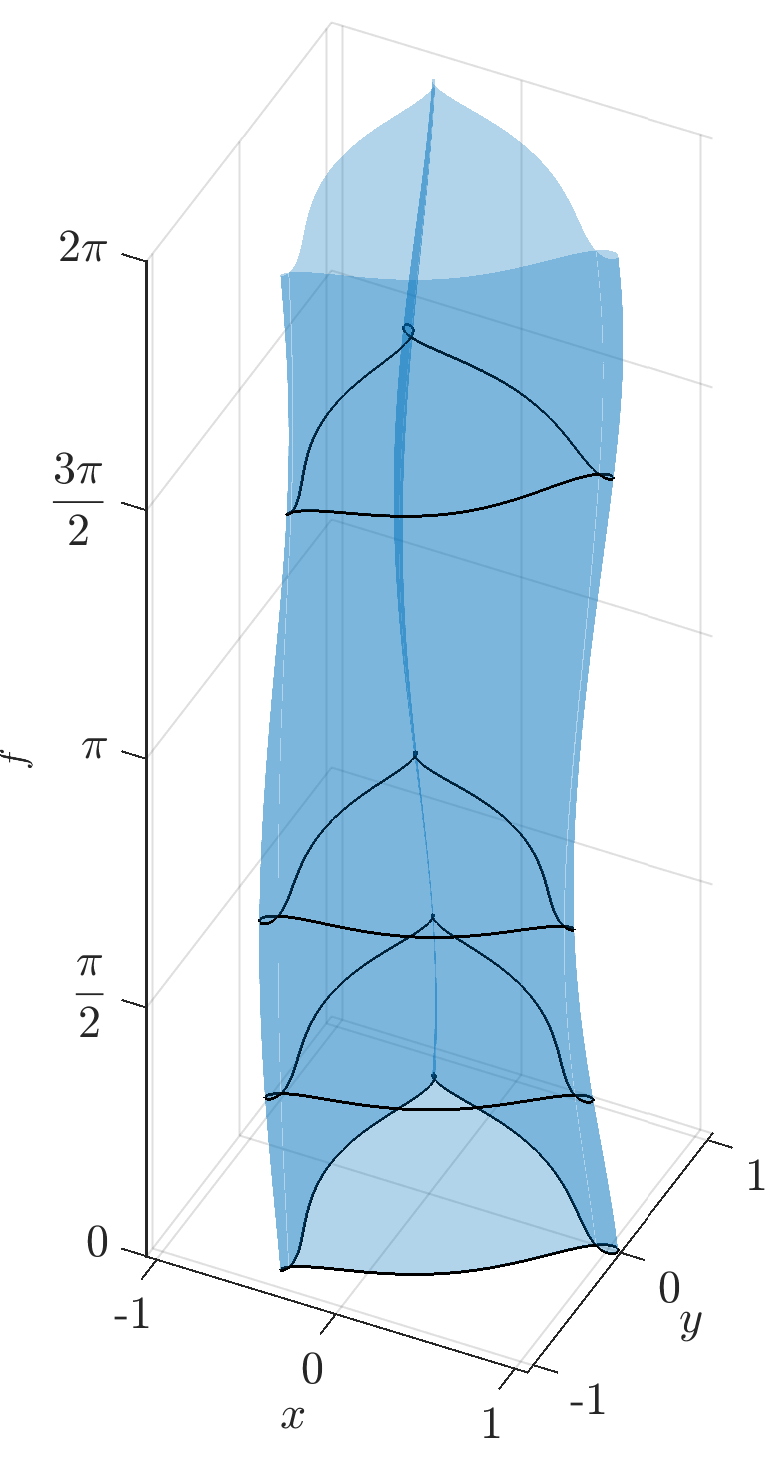}
    \end{minipage}
    \begin{minipage}[t]{0.49\linewidth}
    \centering
        \includegraphics[width=1\linewidth,valign=t]{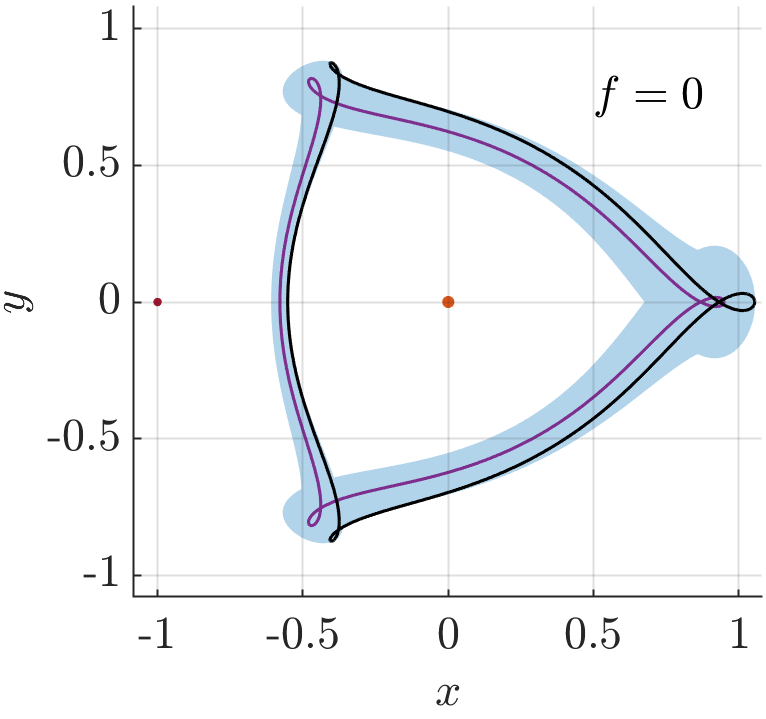}\\
        \includegraphics[width=1\linewidth,valign=t]{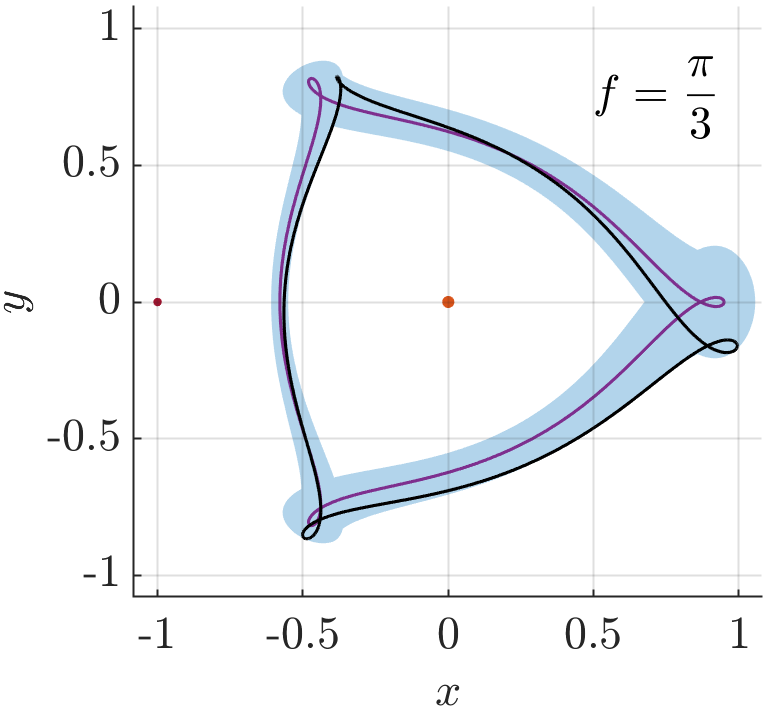}
    \end{minipage}     
        \includegraphics[width=0.49\linewidth]{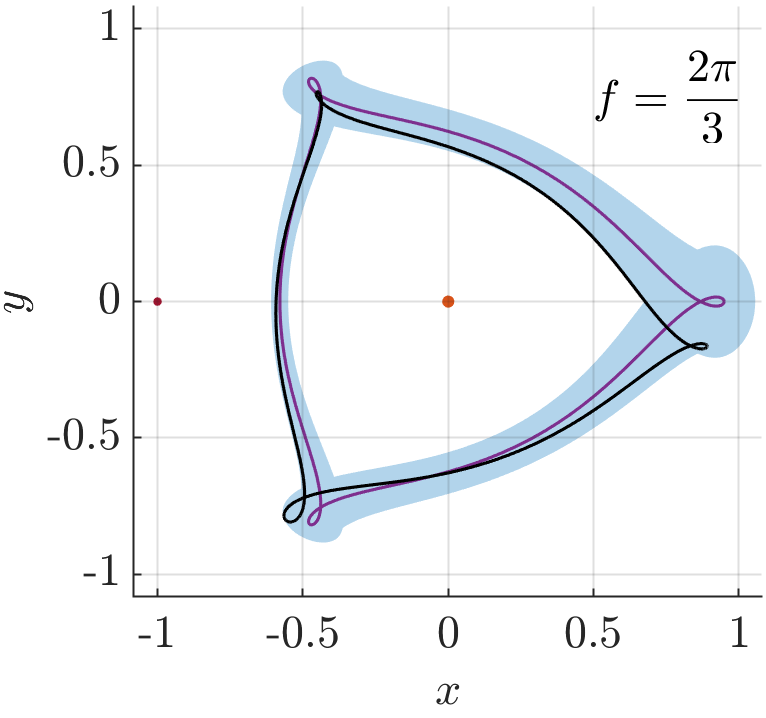}
        \includegraphics[width=0.49\linewidth]{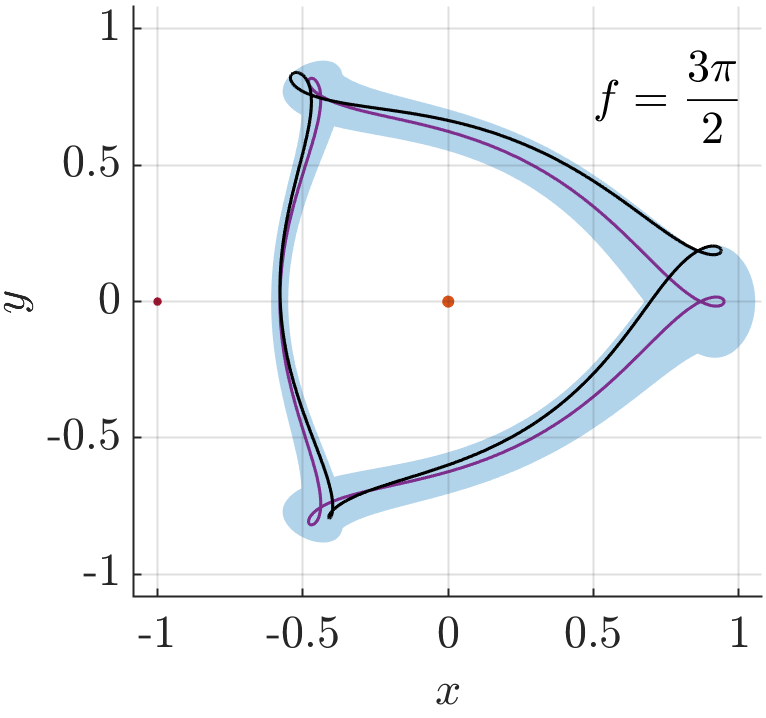}    
    \caption{Top left, 2D invariant tori of the ERTBP for $\rho=3.14878$. Right and bottom left, invariant curves of said 2D invariant tori for different values of the true anomaly, $f$. See text for further details.}
    \label{fig:invCurves}
\end{figure}

In Figure~\ref{fig:invCurves} we study just one of those 2D invariant tori. In particular, is the one with $\rho=3.14878$, coming from a periodic orbit of the CRTBP with $T=12.53796$ and $C=3.006373$. At the top left, the 2D quasi-periodic solution is shown in blue for the flow of ERTBP projected in the $xy-$plane and being $f$ the vertical axis and in black the invariant curves for specific values of $f$. At the right and bottom left, the invariant curves at those values of $f$ are shown together with the periodic orbit of the unperturbed system, coloured in purple. In shadowed blue the area in $xy-$plane swept by the 2D invariant tori for $f\in[0,2\pi]$ is shown.

The stability of these invariant tori has been studied according to Equation~(\ref{eq:GEV}), and we find that these quasi-periodic solutions are elliptic, as the periodic orbits they come from.

\subsection{Three-dimensional quasi-periodic solutions}\label{subsec:E_3D}
This section is devoted to the study of the relation of Hilda asteroids dynamics to a family of 3D invariant tori in the planar ERTBP.

Similarly as the analysis of the CRTBP performed to visualise 2D invariant tori by using Poincaré Section Plots constraining values of the Jacobi Constant, now we aim to visualise three-dimensional quasi-periodic orbits by performing a Poincaré Double Section Plot (PDSP). These PDSPs consist on study the trajectories related to Hilda asteroids in a double section; one spatial (similar to $\Sigma$ used in Section~\ref{subsec:Quasi_C}) and one temporal, using again a stroboscopic map to take advantage of the periodicity of time variable in the elliptical model. That is,
\begin{equation}
\label{eq:2PSP}
    \Sigma_{f^{*}} = \{ f=f^{*}, y=0; \ x<0 \text{ and } y'<0 \},
\end{equation}
where $f^{*}$ is a prescribed value for the true anomaly, for example with $f^{*}=0$ the temporal section would correspond to the stroboscopic map previously used, $\bm{P}$. Notice that, a 2D invariant tori of the flow dynamics is seen as a fixed point of this double section. And therefore, a closed curve found on $\Sigma_{f^{*}}$ corresponds to a three-dimensional quasi-periodic solution of the flow of the ERTBP.

To compute the slice of a torus on a double section we take an orbit on the torus and we compute each cut of the orbit with $f=f^*$ and, if its $y$ coordinate satisfies $|y|\le\delta$ for $\delta$ small enough (we have used values of $\delta$ between $10^{-5}$ and $10^{-4}$).
Then, with these points we interpolate to compute the curve in $y=0$.
These Poincar\'e  Double Section Plots have already been used in the literature\cite{Simo98,BroerHJVW03}.

\begin{figure}[!ht]
    \centering
    \includegraphics[width=0.48\linewidth]{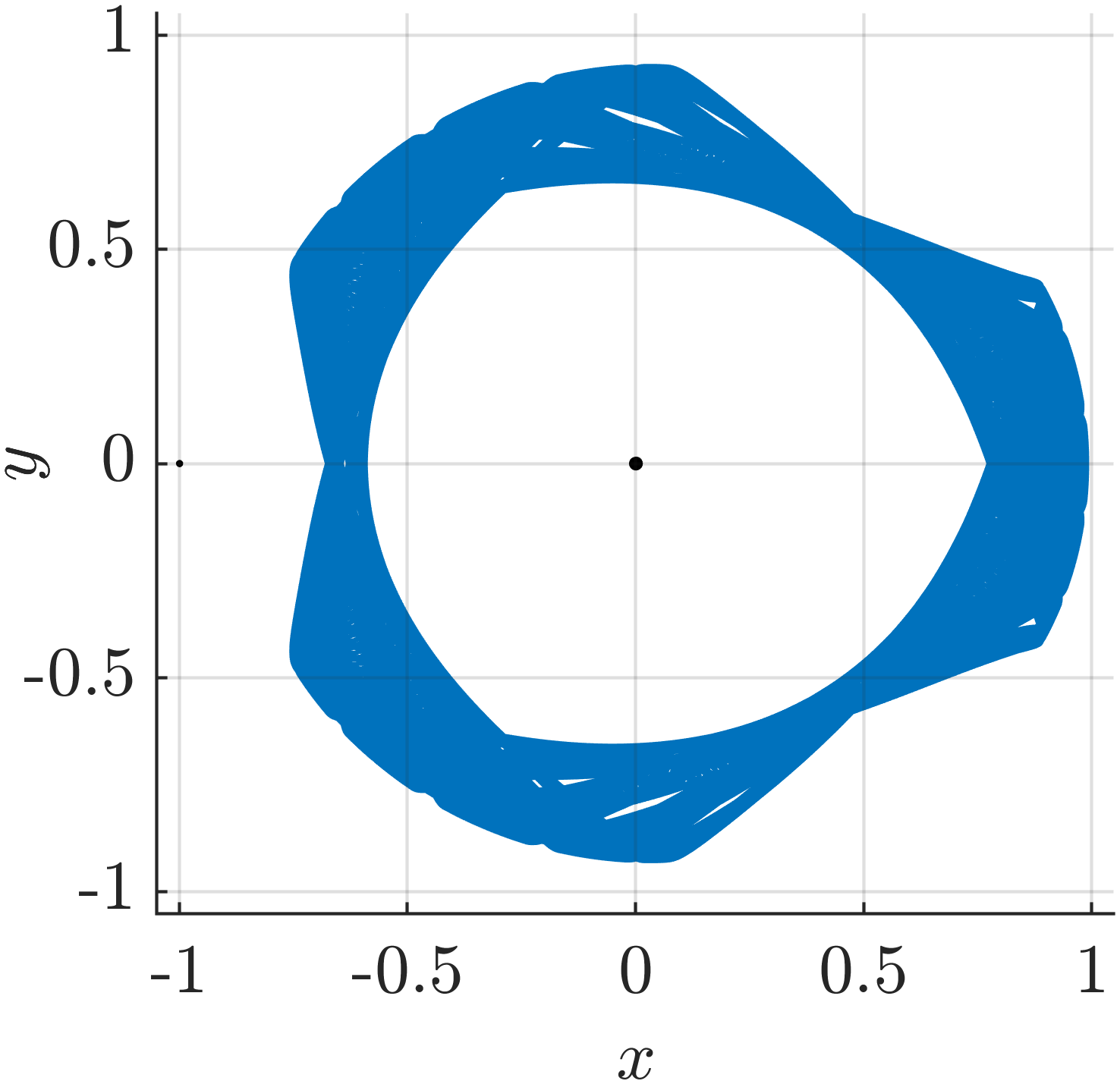}
    \includegraphics[width=0.48\linewidth]{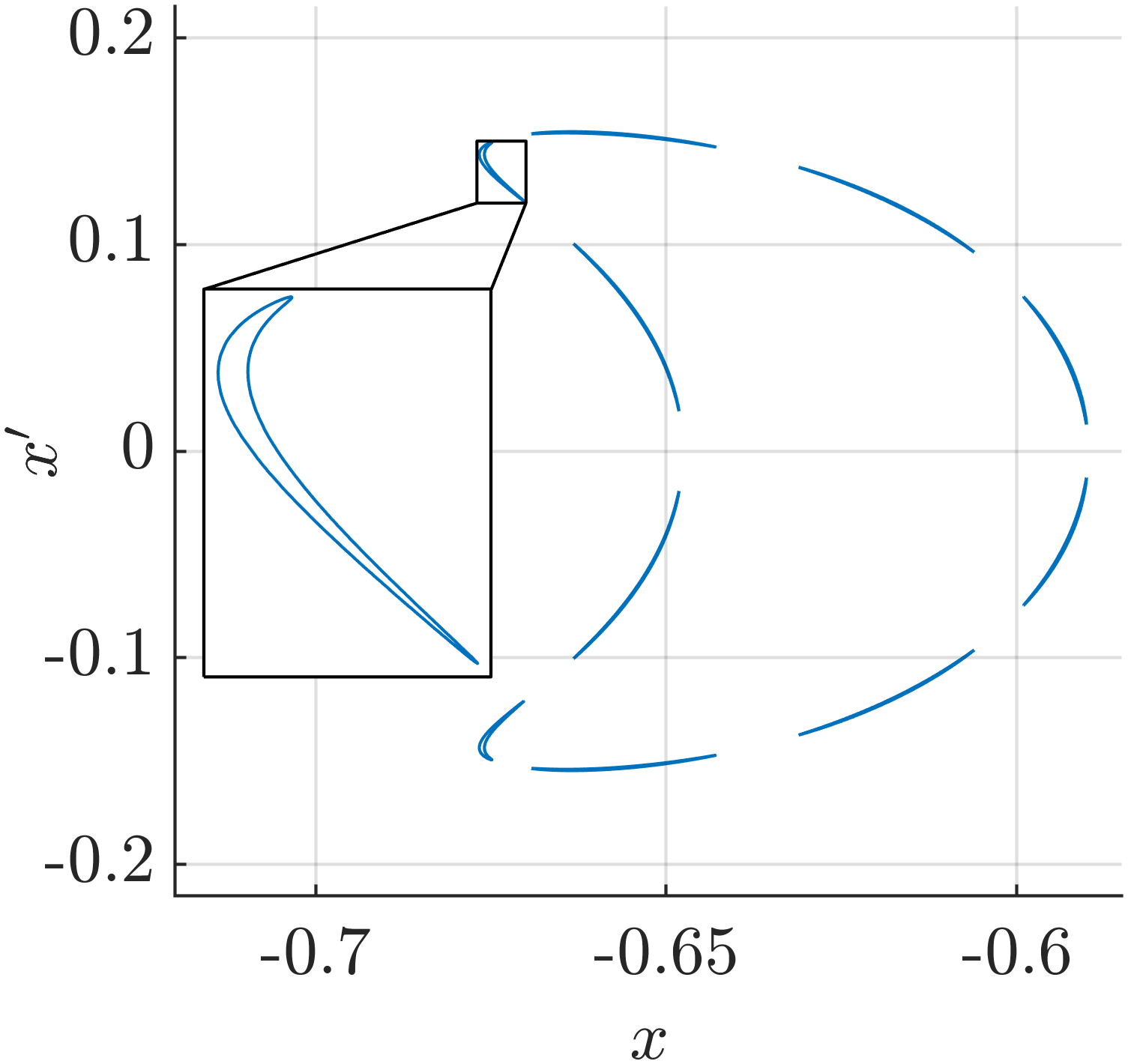}
    \includegraphics[width=0.48\linewidth]{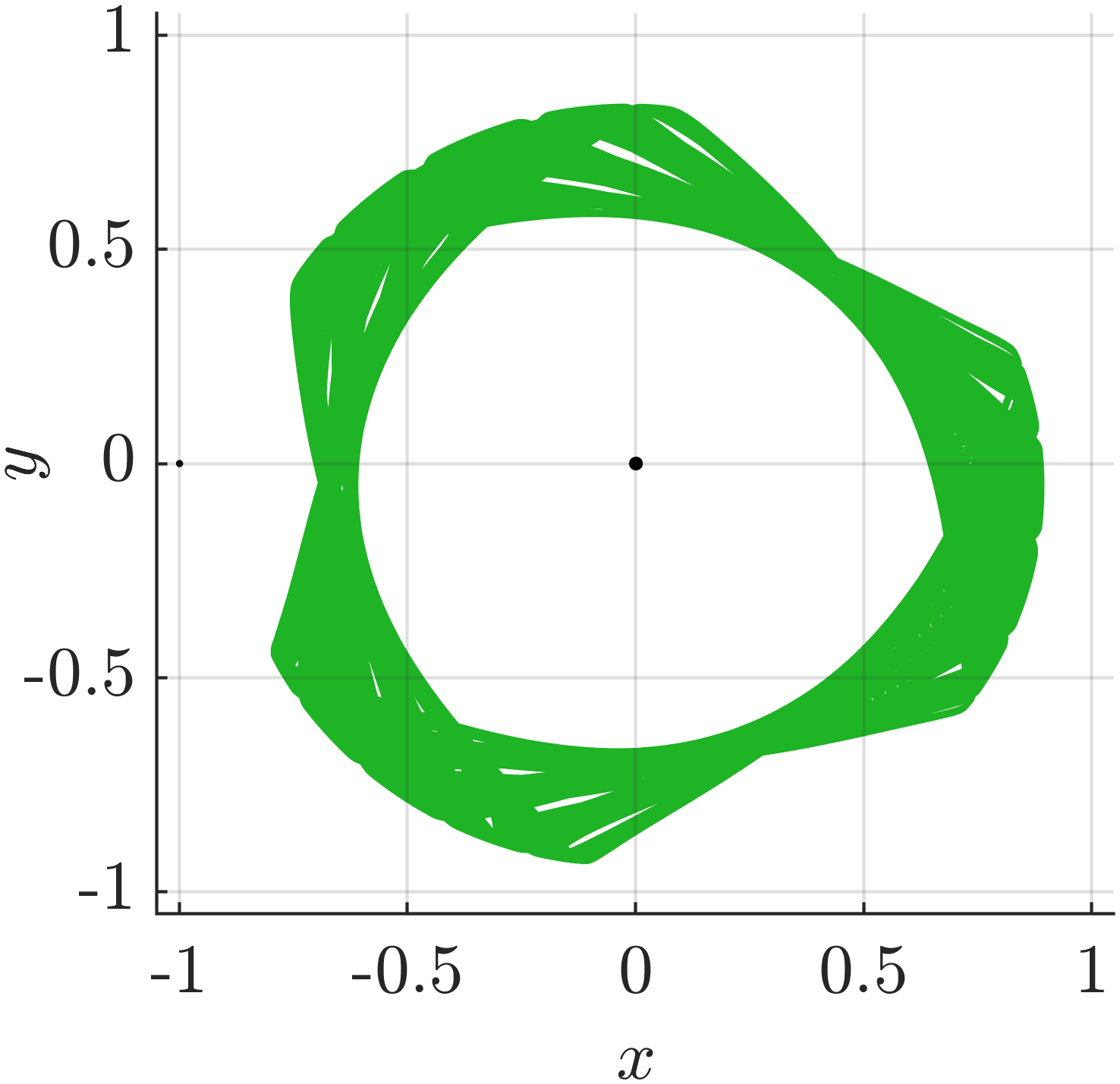}
    \includegraphics[width=0.48\linewidth]{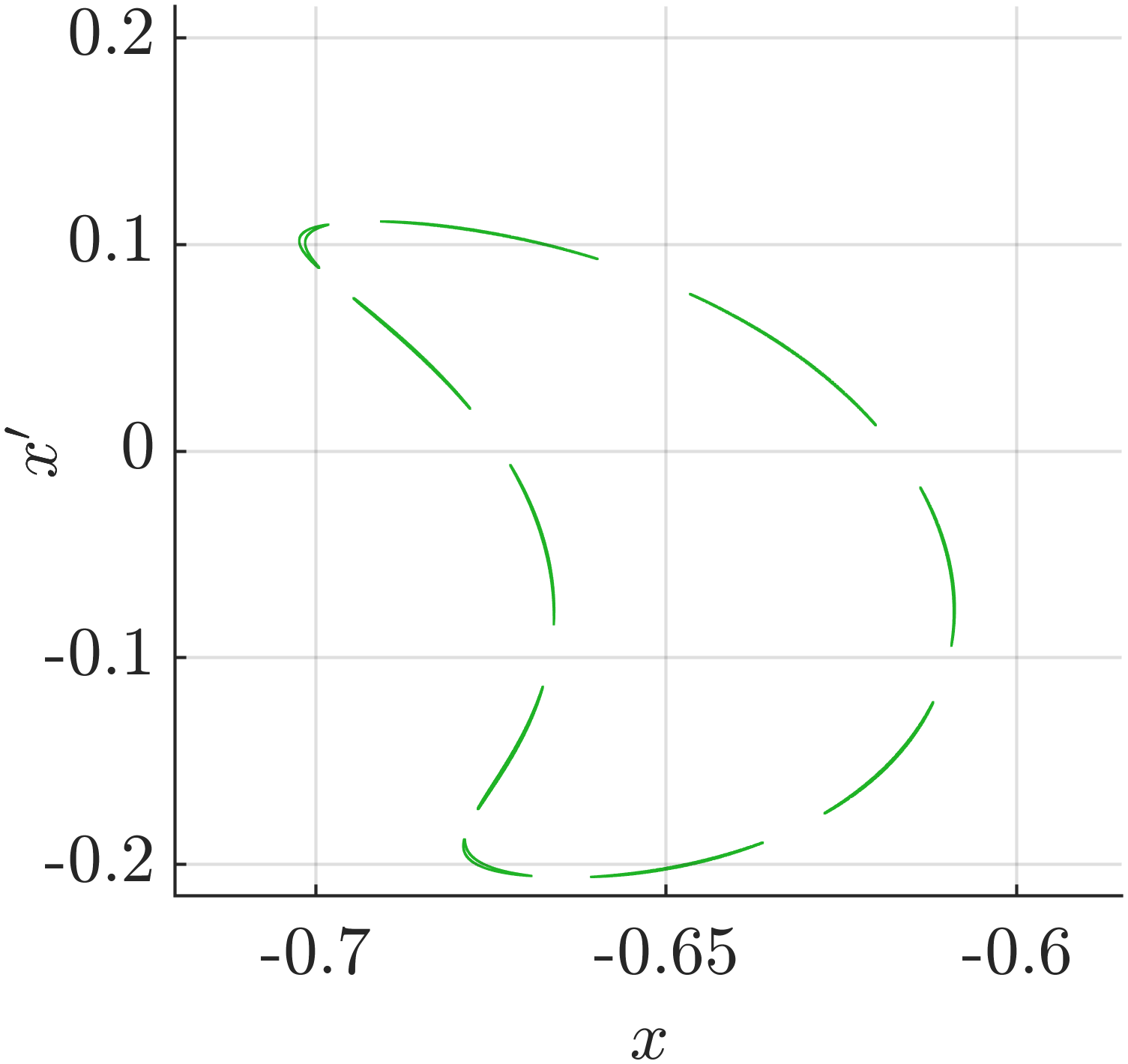}
    \includegraphics[width=0.48\linewidth]{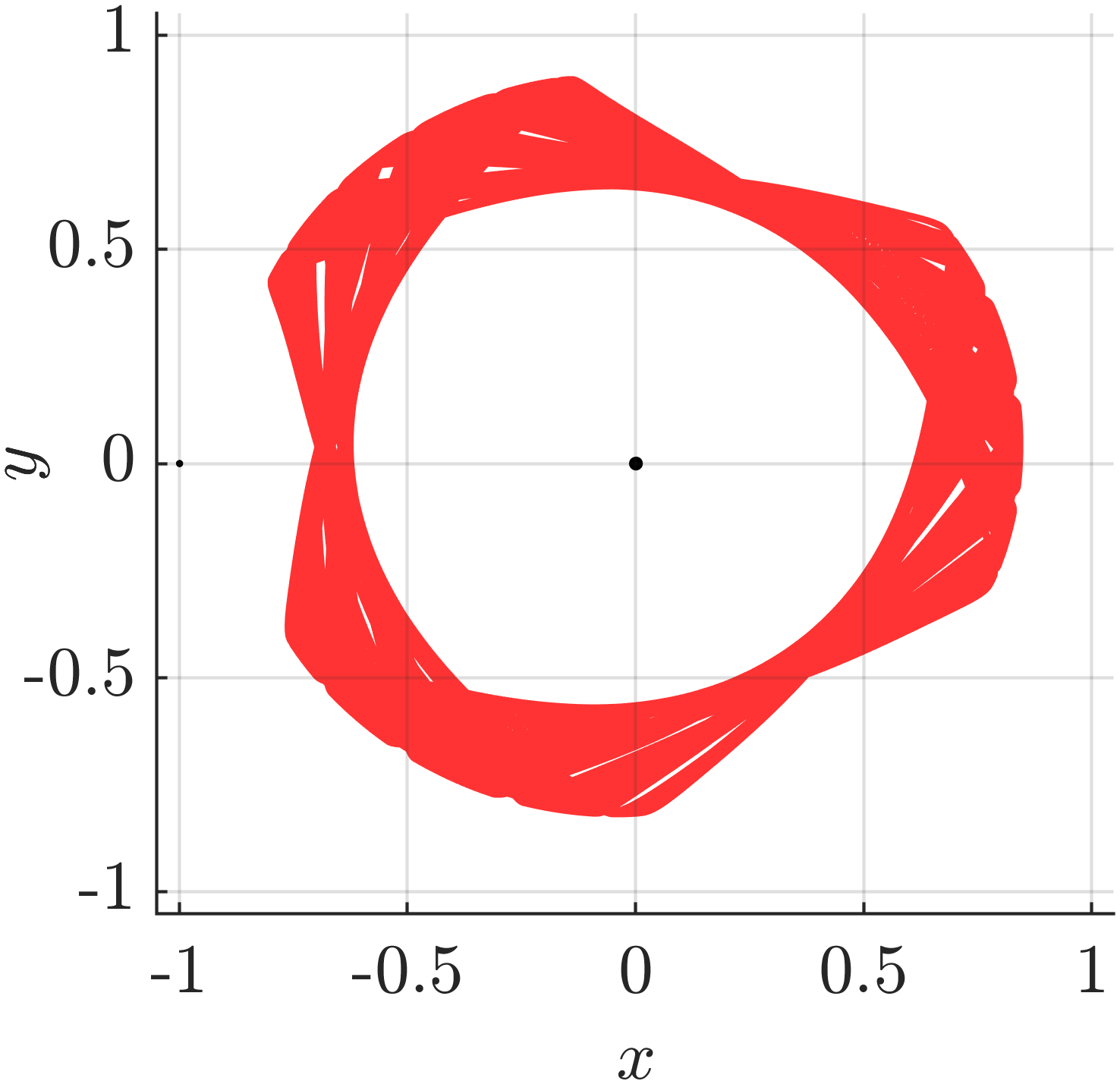}
    \includegraphics[width=0.48\linewidth]{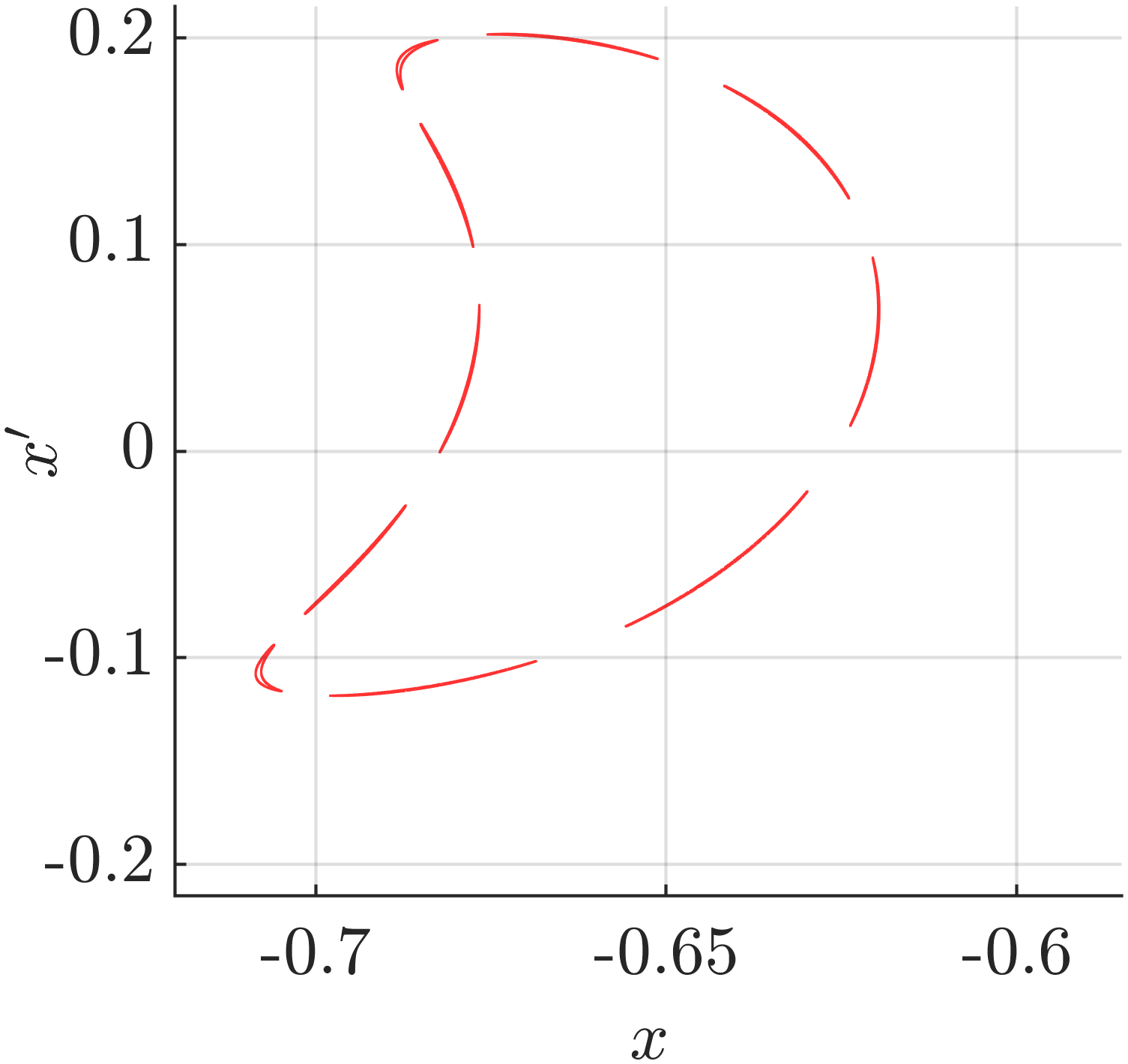}
    \caption{Left, positions of (153) Hilda asteroid trajectory in temporal sections $f^{*}=0$, $\pi/2$, $4\pi/3$, from up to down. Right, intersection of the asteroid trajectory with the double sections $\Sigma_0$, $\Sigma_{\pi/2}$ and $\Sigma_{4\pi/3}$, defined as in (\ref{eq:2PSP}).}
    \label{fig:torosdoblesecHilda}
\end{figure}

First, we start by visualising the six selected asteroids shown in Figure~\ref{fig:trasHildas} in the ($x$,$x'$) plane for this PDSP of the ERTBP. In Figure~\ref{fig:torosdoblesecHilda}, plots of the (153) Hilda asteroid are shown. In the first column, its positions are plotted at three different values of the true anomaly, $f^{*}=0$, $\pi/2$ and $4\pi/3$. At the right of each of those plots, the asteroid is shown in the double section $\Sigma_0$, $\Sigma_{\pi/2}$ and $\Sigma_{4\pi/3}$, defined as in (\ref{eq:2PSP}). The asteroid (153) Hilda describes 10 closed curves in the double PSPs, showing that its motion belongs to a resonance of order 10. A zoom of one of the islands is included for clarity. For the value of true anomaly $f^{*}=0$, the image in double section $\Sigma_0$ is symmetric with respect to the $x$ axis. Meanwhile, those in $\Sigma_{\pi/2}$ and $\Sigma_{4\pi/3}$ do not present that symmetry.

\begin{figure}[!ht]
    \centering
    \includegraphics[width=0.48\linewidth]{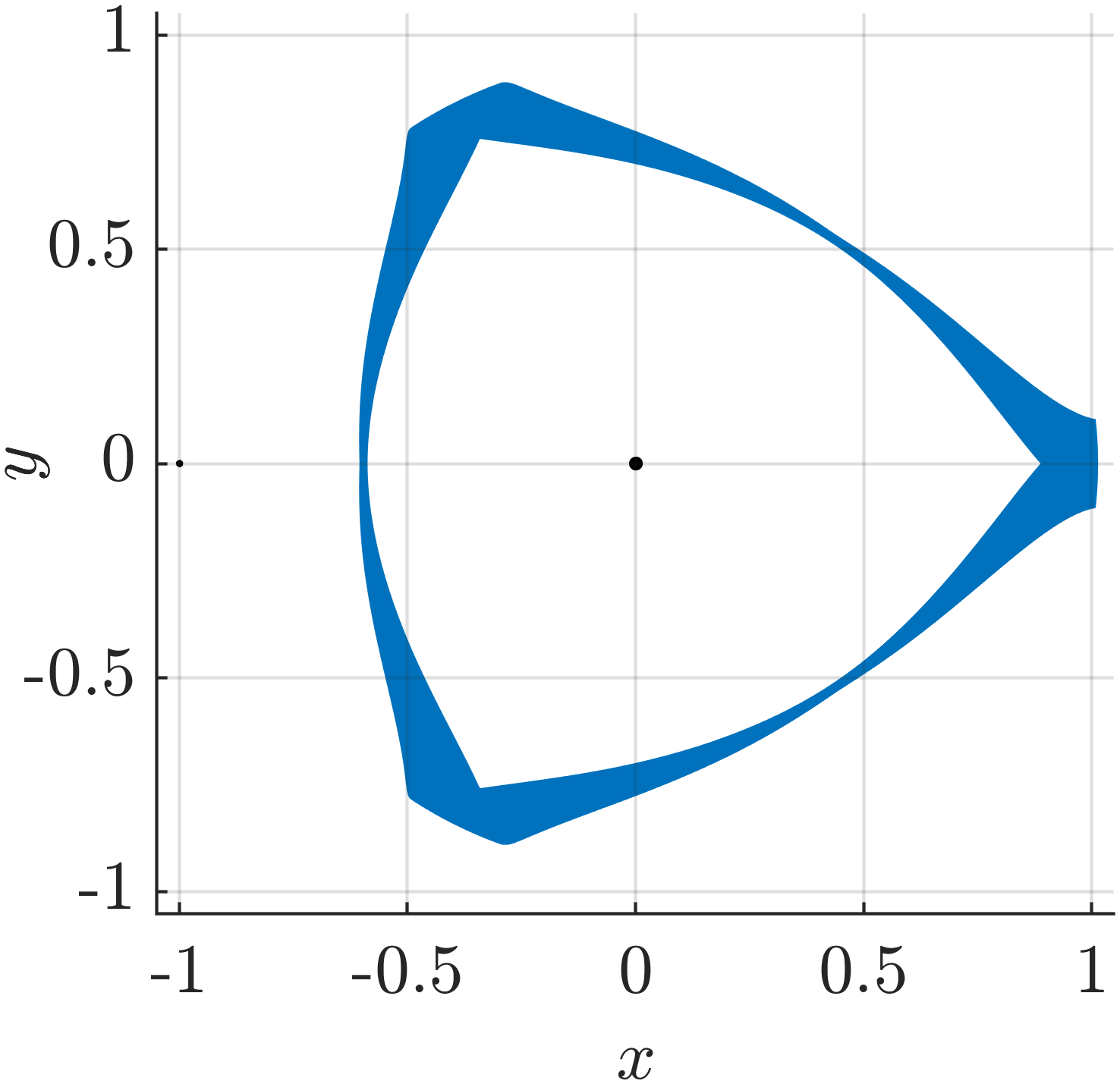}
    \includegraphics[width=0.48\linewidth]{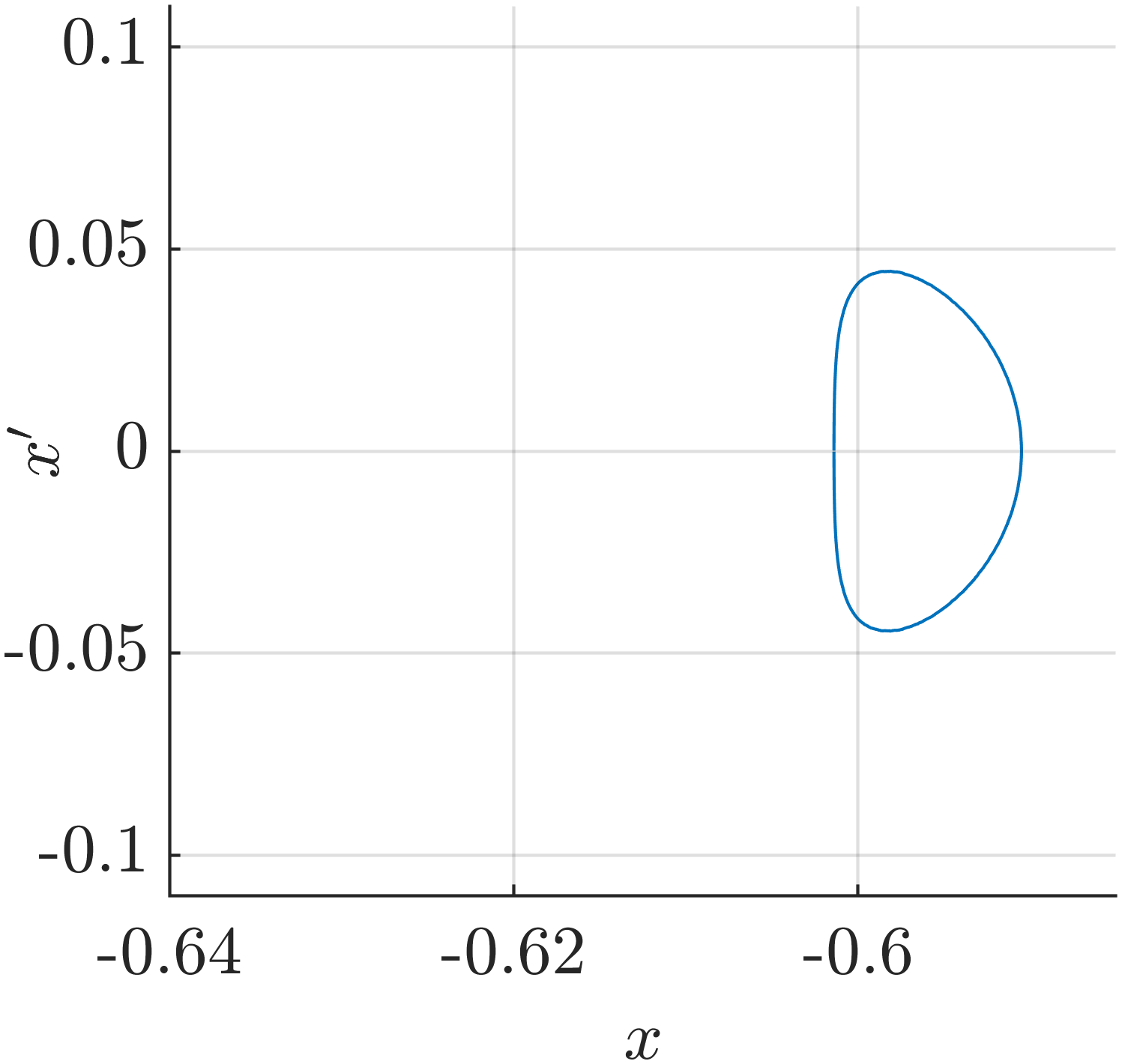}
    \includegraphics[width=0.48\linewidth]{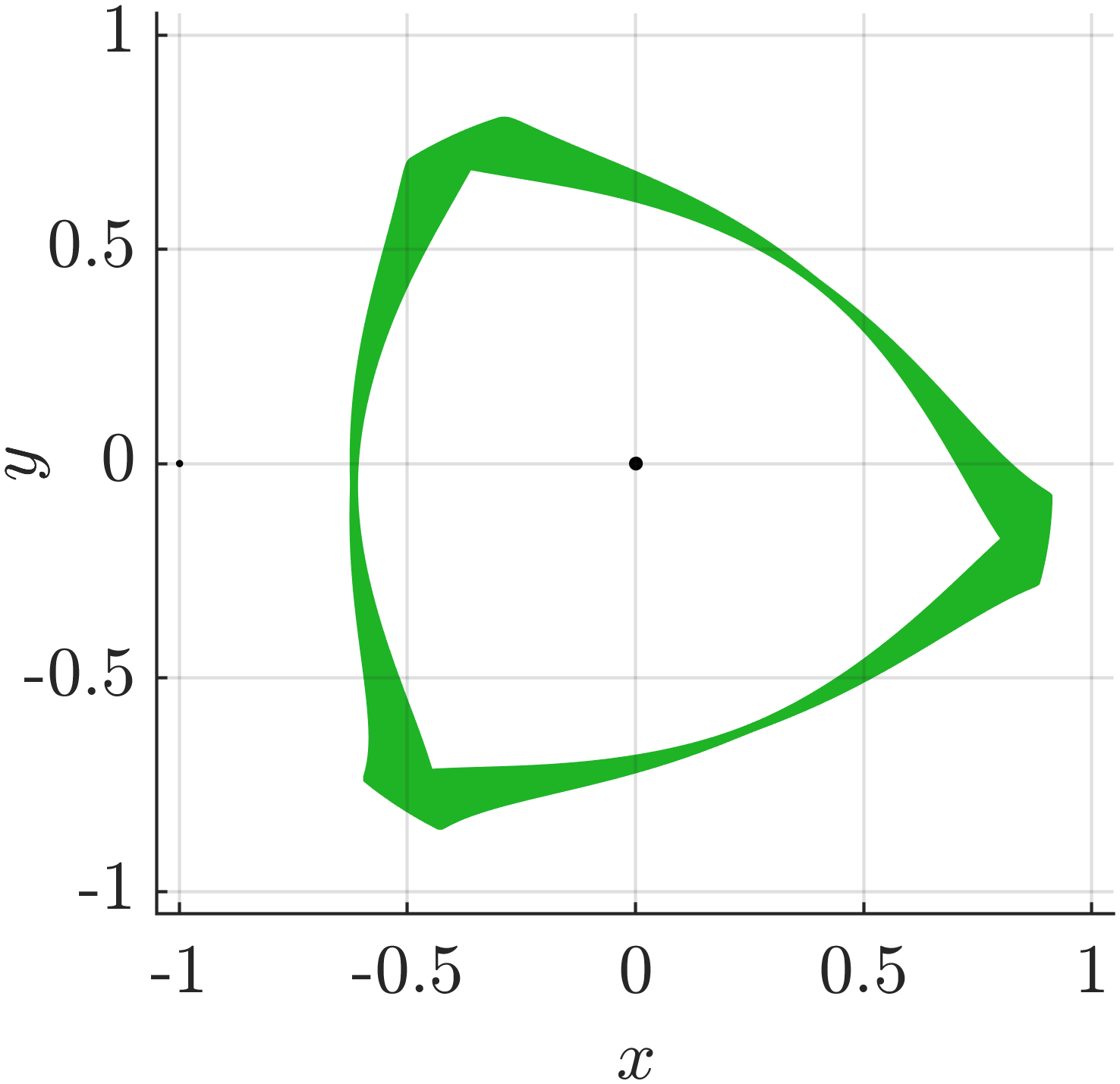}
    \includegraphics[width=0.48\linewidth]{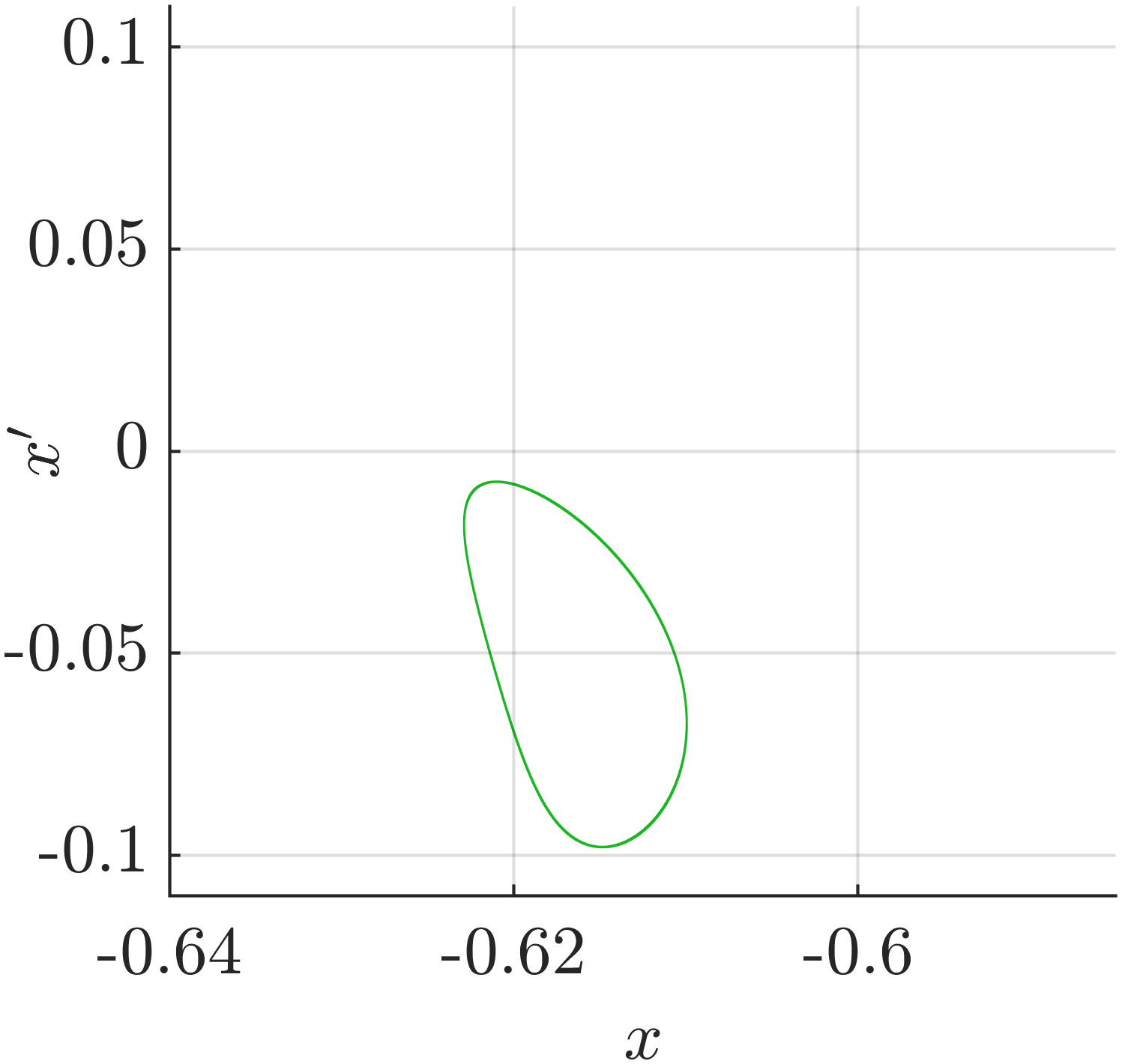}
    \includegraphics[width=0.48\linewidth]{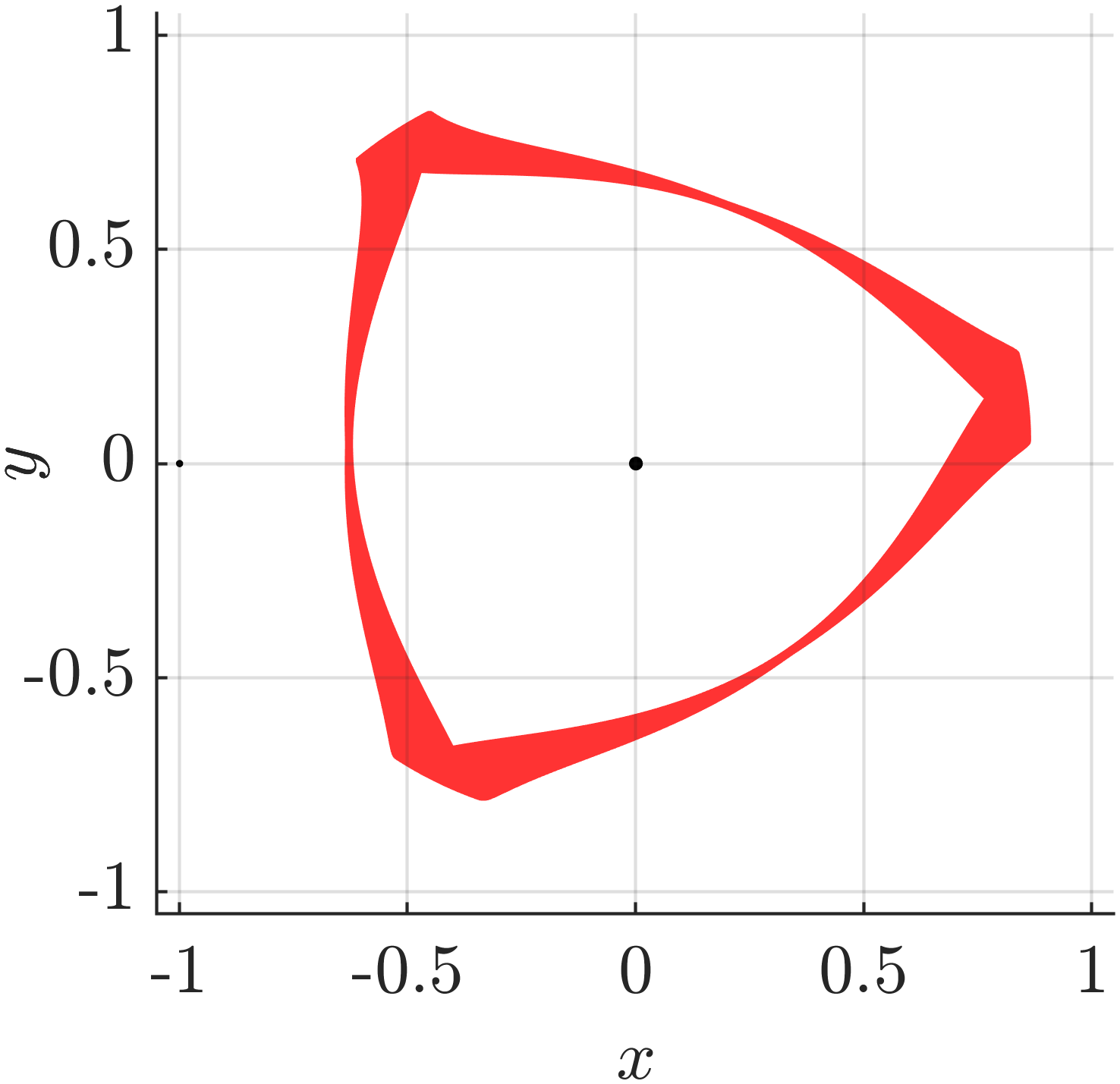}
    \includegraphics[width=0.48\linewidth]{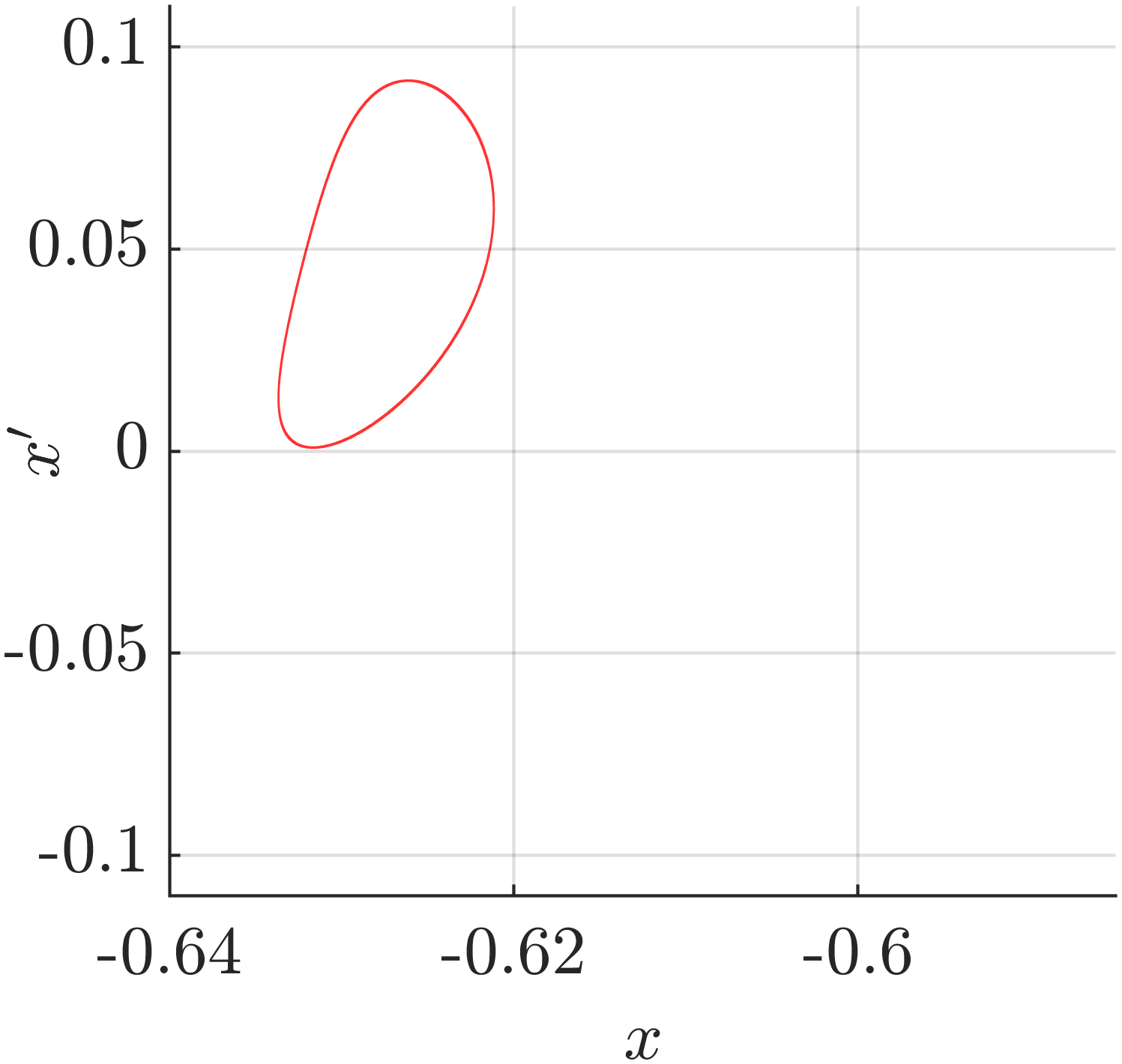}
    \caption{Left, positions of (1911) Schubart asteroid trajectory in temporal sections $f^{*}=0$, $\pi/2$, $4\pi/3$, from up to down. Right, intersection of the asteroid trajectory with the double sections $\Sigma_0$, $\Sigma_{\pi/2}$ and $\Sigma_{4\pi/3}$, defined as in (\ref{eq:2PSP}).}
    \label{fig:torosdoblesecSchubart}
\end{figure}

Figure~\ref{fig:torosdoblesecSchubart} shows the same plots for asteroid (1911) Schubart. Also for temporal sections $f^{*}=0$, $\pi/2$ and $4\pi/3$. Again, we observe the symmetry with respect to the $x$ axis in $\Sigma_0$. However, now we only observe one closed curve in the double sections. This same behaviour is found for the other four asteroids, suggesting that their motion is trapped in a three dimensional quasi-periodic solution of the flow dynamics.

In the CRTBP the points constituting the curves in the PSP share the Jacobi constant. 
In the ERTBP there is not a first integral as in the CRTBP. Therefore there is not an easy way, nor unique, to generate initial conditions dynamically related to the trajectory of a given asteroid, in order to plot them together in $\Sigma_{f^{*}}$ describing a family of 3D invariant tori to prove that the asteroid motion is part of said family.

That is our goal, and for achieving it we have designed the following procedure for the double section $\Sigma_0$ in order to take advantage of the mentioned symmetry \eqref{eq:symE}.
In spite of the fact that in the ERTBP there is not a first integral, there is an invariant relation that we can use, Equation~\eqref{eq:C_ERTBP}, 
\[
    \mathcal{C} = \frac{2\tilde{\Omega}}{1+e\cos f} - 2 e \int_{0}^{f} \frac{\tilde{\Omega} \sin f}{(1+e\cos f)^2} df - \left({x'}^2 + {y'}^2\right),
\]
taking $z=0$, since we are studing the planar case.

For each point in the PDSP of the asteroids we can evaluate this expression. Notice that, as we are interested in its evaluation in the stroboscopic map $P$, the integral term can be taken as zero, and the term $\cos f=1$. 

Given the curve for one asteroid in the double section $\Sigma_0$, $\mathcal{C}$ takes a slightly different value for each point in the curve, taking the minimum value at the right cross of the axis $x$ ($x'=0$). Then, we use that point to find the two-dimensional invariant tori satisfying Equation~\eqref{eq:invcond} with that value of $\mathcal{C}$ at $f=0$, $y=0$, $x'=0$. This 2D invariant object is seen as a fixed point of the double section $\Sigma_0$.

Once we have that point, new initial conditions are generated in the double section along the axis $x$ ($x'=0$), by increasing $x$ and modifying $y'$ such that $\mathcal{C}$ value remains. The initial conditions are propagated in the ERTBP for long spans of time saving the crosses to the PDSP in Equation~\eqref{eq:2PSP}. For each of the initial conditions we find one closed curve, or set of closed curves, in the ($x$,$x'$) plane. All together constitute a family of three-dimensional quasi-periodic solutions where the motion of the Hilda asteroids takes place.

\begin{figure}
    \centering
    \includegraphics[width=0.48\linewidth]{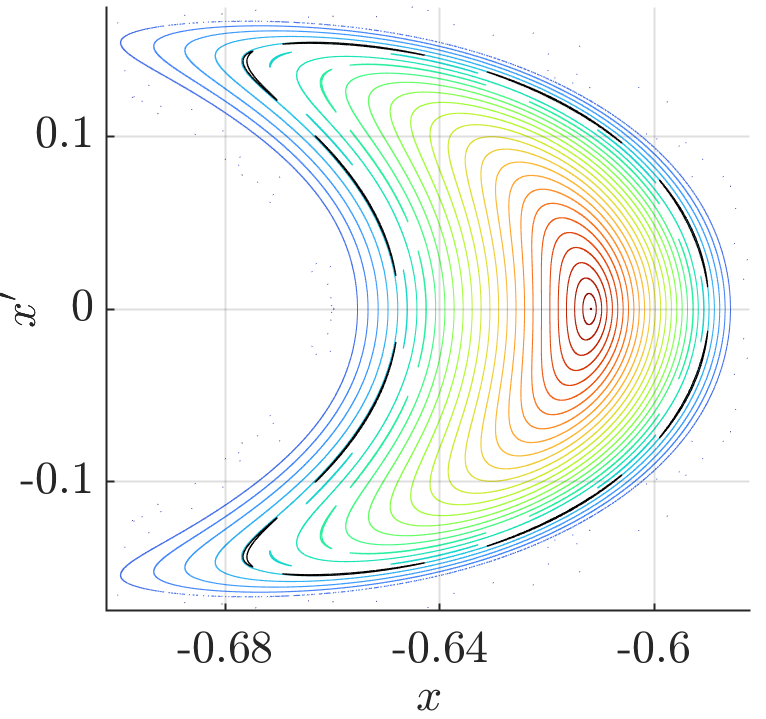}
    \includegraphics[width=0.48\linewidth]{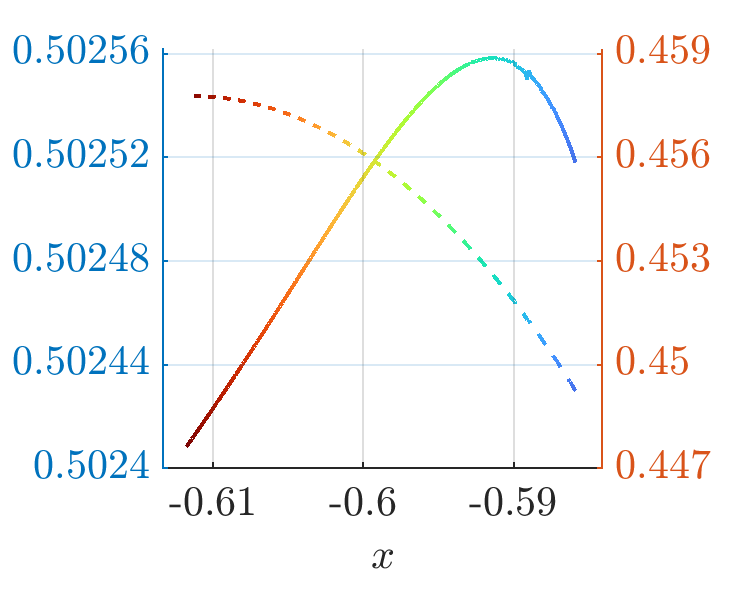}
    \includegraphics[width=0.48\linewidth]{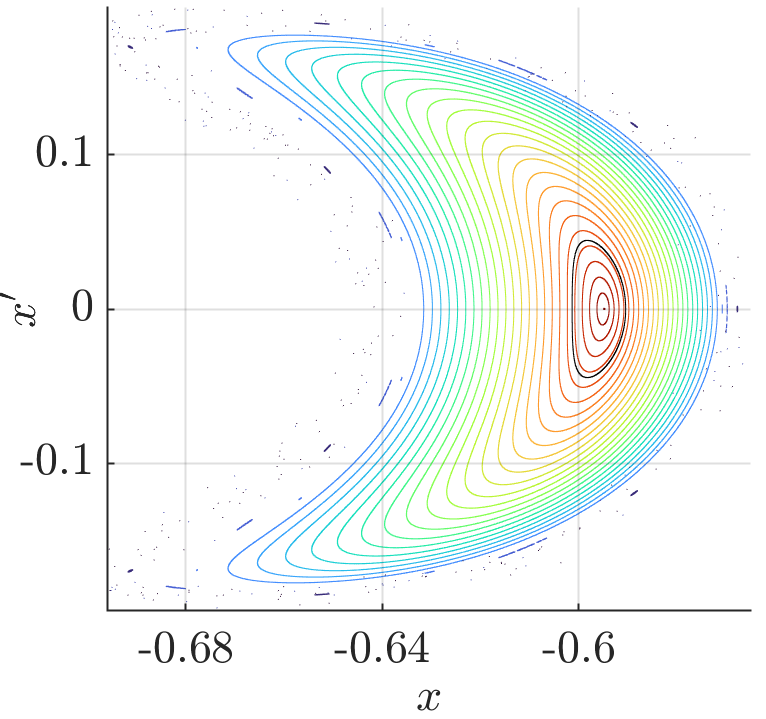}
    \includegraphics[width=0.48\linewidth]{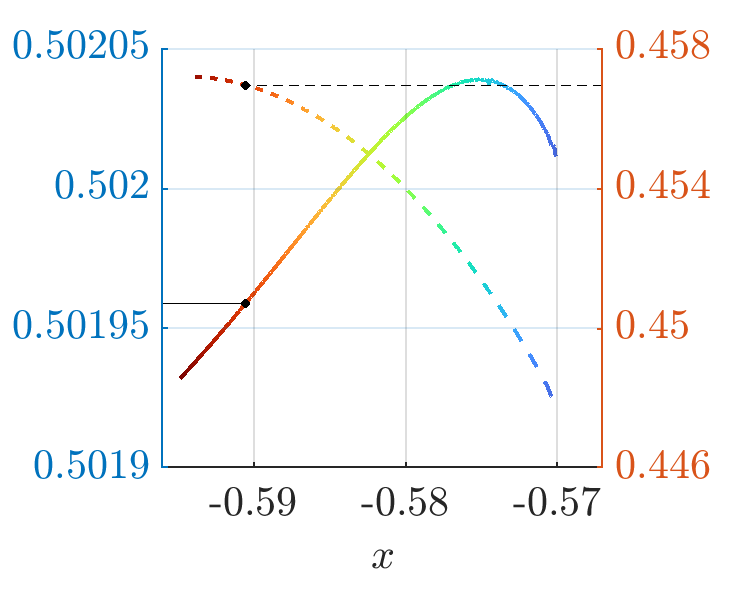}
    \includegraphics[width=0.48\linewidth]{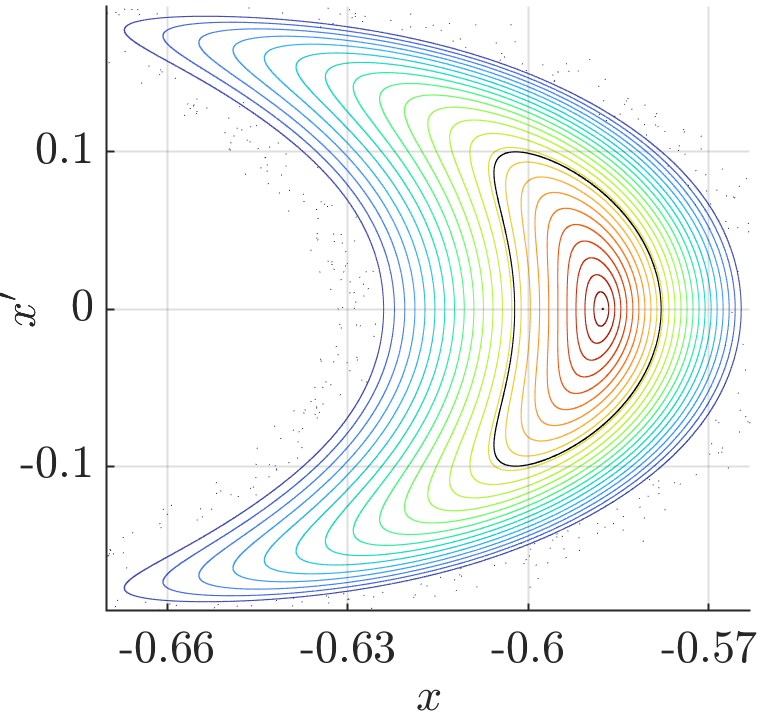}
    \includegraphics[width=0.48\linewidth]{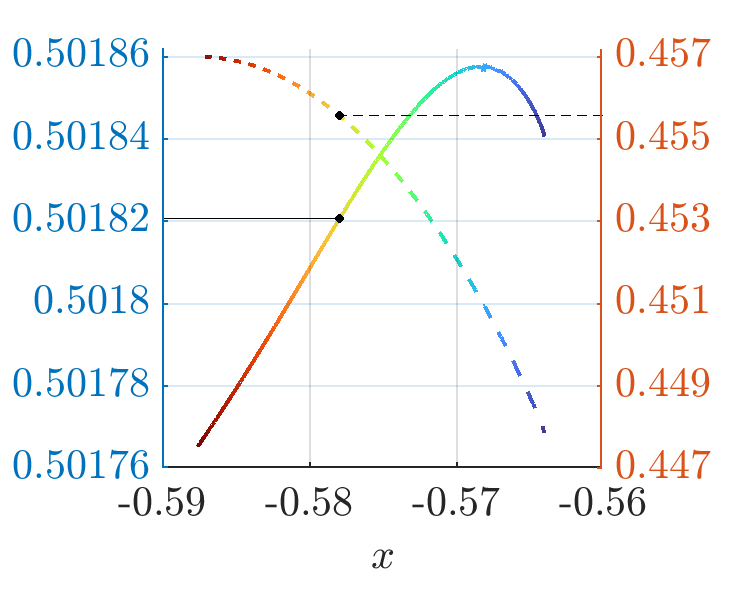}
    \caption{Left column, double sections plots generated using $\mathcal{C}$ value of asteroids (153) Hilda, (1911) Schubart and (499) Venusia. Black curves correspond to the intersections with $\Sigma_0$ of said asteroids. Right column, frequency analysis of the 3D quasi-periodic solutions at the left. See text for further details.}
    \label{fig:doublePSP_frecs1}
\end{figure}

\begin{figure}
    \centering
    \includegraphics[width=0.48\linewidth]{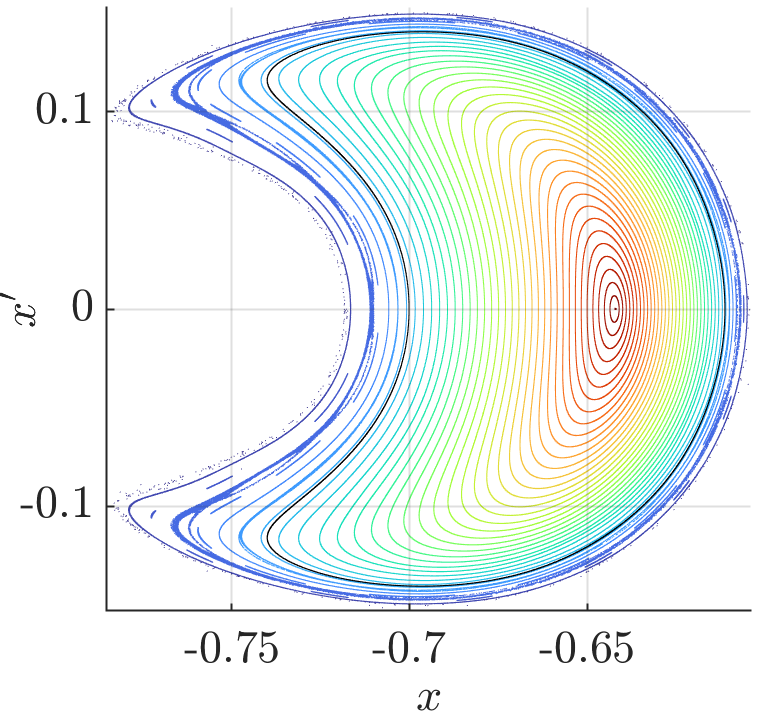}
    \includegraphics[width=0.48\linewidth]{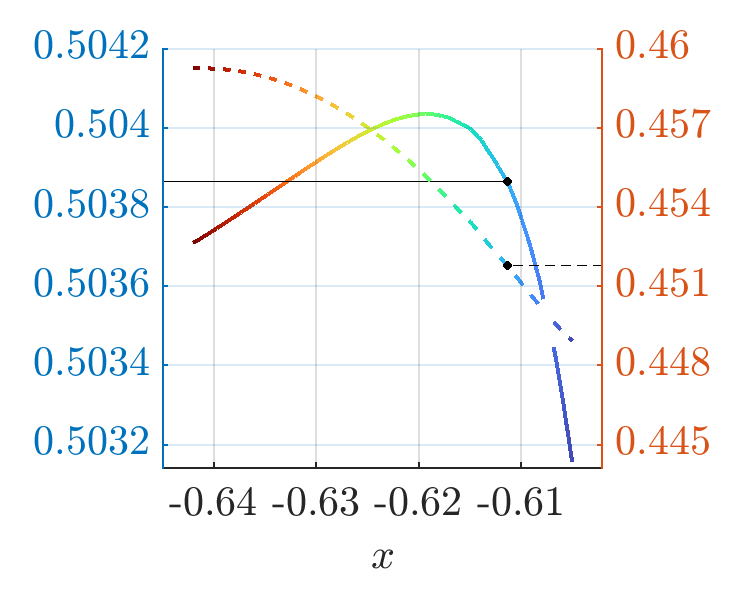}
    \includegraphics[width=0.48\linewidth]{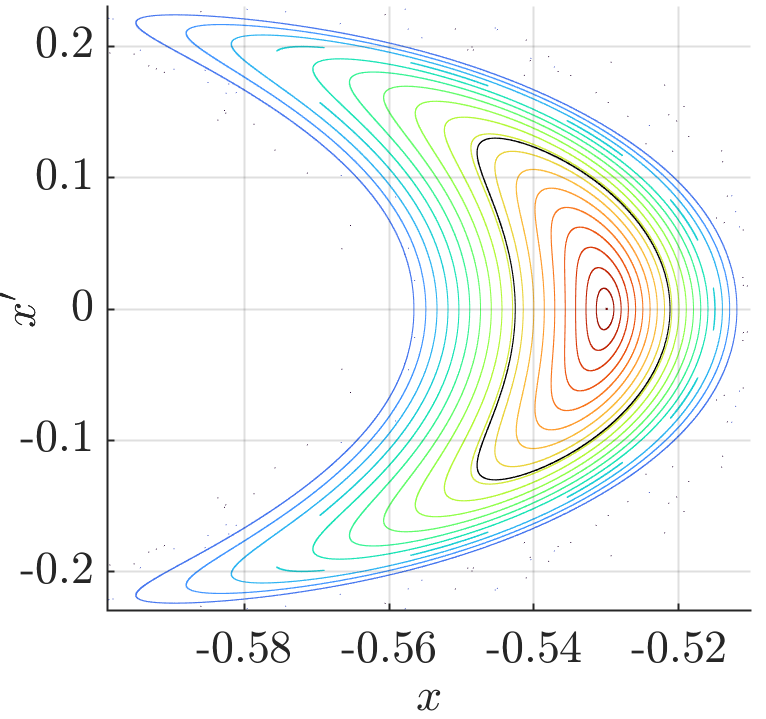}
    \includegraphics[width=0.48\linewidth]{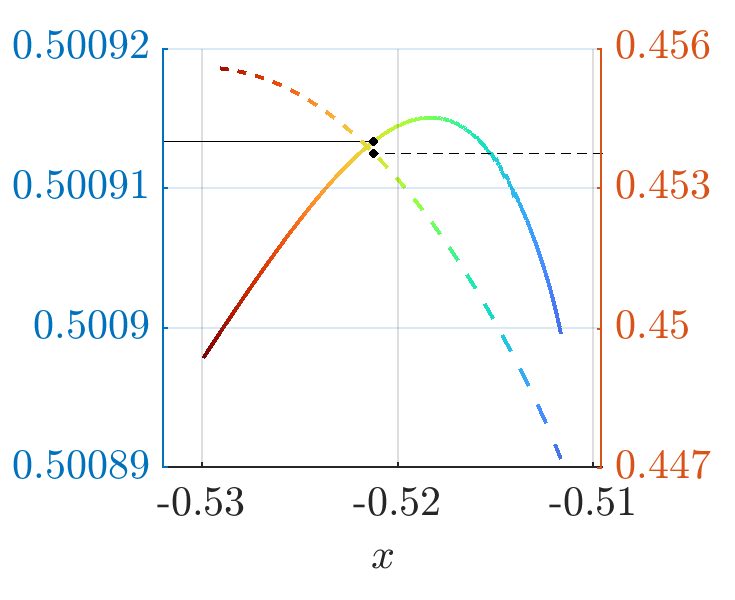}
    \includegraphics[width=0.48\linewidth]{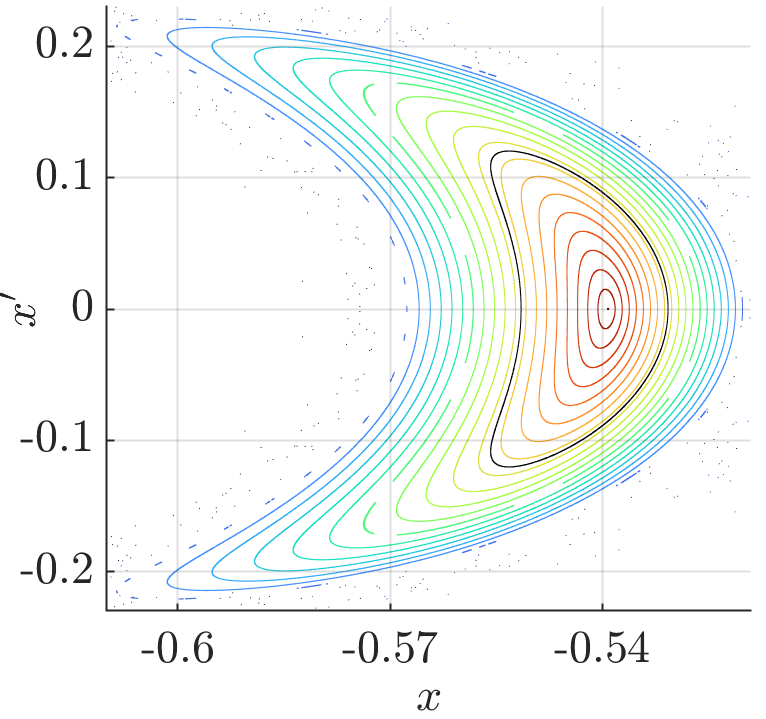}
    \includegraphics[width=0.48\linewidth]{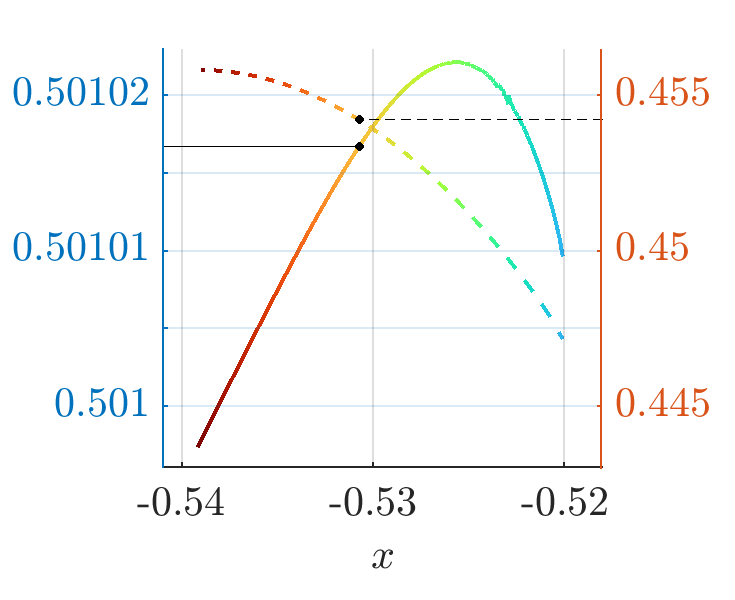}
    \caption{
    Left column, double sections plots generated using $\mathcal{C}$ value of asteroids (1038) Tuckia, (4446) Carolyn and (2483) Guinevere. Black curves correspond to the intersections with $\Sigma_0$ of said asteroids. Right column, frequency analysis of the 3D quasi-periodic solutions at the left. See text for further details.
    }
    \label{fig:doublePSP_frecs2}
\end{figure}

In Figures~\ref{fig:doublePSP_frecs1} and ~\ref{fig:doublePSP_frecs2}, the double section plots for the six families of three-dimensional quasi-periodic solutions are shown. In both figures, left collumn show in black intersections of the trajectories of the six asteroids with the double section $\Sigma_0$, and the other points, in colour, correspond to the families generated through the procedure explained above. 
Moreover, we have also computed the frequencies of these quasi-periodic solutions (see Appendix~\ref{sec:fran} for details).
At the right of each PDSP, the variation of the two main frequencies of the quasi-periodic solutions for each family is shown as they cut the $x$-axis, that is when the invariant curves intersect with $x' = 0$. Left vertical axis corresponds to the first main frequency, plotted in solid line, and right vertical axis corresponds to the second main frequency, plotted in dash line. Notice that one of the three frequencies defining these 3D quasi-periodic solutions, has to be equal to frequency of the perturbation, that is $1$ as shown in Table~\ref{tab:frecus}. In the frequency plots, the horizontal black lines are used to indicate the two main frequencies of the asteroids. Note that, as (153) Hilda asteroid does not cut the $x$-axis in $\Sigma_0$, its two main frequencies can not be pointed out in the plot.

Similarly as we observed in Section~\ref{subsec:Quasi_C}, the motion of the six representative asteroids of the Hilda group seems to be contained in families of quasi-periodic orbits, whose main frequencies are restricted to small ranges of values close to $0.5$, that is, close to the $4\pi$ resonance.

\section{Conclusions}\label{sec:Conclu}

Actual classification of asteroids is based on their orbital elements, analysed through the two-body problem. This is adequate when the orbit of the asteroid is close to Keplerian, however, as some asteroids are significantly perturbed by other massive bodies in the Solar system, their orbital elements suffer variations in time. This can lead to the misclassification of these asteroids.

A classification of asteroids according to their dynamics in the Circular Restricted Three-Body Problem seems to be a better approximation to reality than the two-body problem. Therefore this paper is a first step in this direction.
Through the definitions of different Poincaré maps, we have presented how dynamics of Hilda group of asteroids is contained in islands of two-dimensional quasi-periodic solutions surrounded by a chaotic sea. A frequency analysis has allowed us to determine the frequencies of these quasi-periodic motions.

In order to validate the results obtained in the frame of the planar Circular RTBP, we have analysed the dynamics of Hilda asteroids in the frame of a planar Elliptical RTBP, to account for the main perturbation on the CRTBP, that is the non circular motion of the primaries. ERTBP constitutes a 2$\pi$ periodic perturbation of the CRTBP, where the family of 2D invariant tori studied in the the circular case become (under non resonance conditions) 3D invariant tori. 
Our study of these solutions through double sections reveals that the dynamics of Hilda-type asteroids remain confined within the described islands, maintaining the same principal frequencies as those derived from the circular case. This similarity in frequencies between the tori of the ERTBP and the CRTBP suggests that, despite not accounting for Jupiter's eccentricity, we can still preserve essential information about the invariant families of objects.

Furthermore, this dynamical behavior and the underlying information cannot be adequately explained by a model that neglects the perturbative influence exerted by Jupiter, that is, with a two-body model. The preservation of these invariant structures, despite the added complexity of eccentric motion, reinforces our belief that the analysis performed for the CRTBP provides a robust framework for identifying the members of the asteroid family. Thus, the transition from the circular to elliptical model not only validates our earlier findings but also emphasizes the critical role of Jupiter's perturbation in understanding the dynamics of these group of asteroids.

We are currently working in an extension of this procedure to analyse other groups of asteroids. Besides, we are interested in developing a similar study for spatial models, in order to take inclinations of asteroids into account. 
Some asteroids exhibit significant inclinations, making it essential to incorporate a three-dimensional model into our analysis. We believe that this three-dimensional framework will allow us to explore the complexities introduced by substantial inclination values, which can significantly affect the motion of the asteroids. This approach will enable us to compare the results of our dynamical analysis, which accounts for Jupiter's perturbation, with the current orbital elements used nowadays.

\begin{acknowledgments}
The authors want to thank J.M. Mondelo for sharing with them
his frequency analysis code \cite{GomezMondeloFourier1,GomezMondeloFourier2}
that has proven to be both efficient and accurate.
The project has been supported with the Spanish grants
PID2021-123968NB-I00 and PID2021-125535NB-I00 (MICIU/AEI/10.13039/501100011033/FEDER/UE), and the Catalan grant 2021 SGR 01072.
This work has also been also funded through the Severo Ochoa and
Mar\'ia de Maeztu Program for Centers and Units of Excellence in R\&D (CEX2020-001084-M).
\end{acknowledgments}

\section*{Data Availability Statement}

The data that support the findings of this study are available from the corresponding author upon reasonable request.

\appendix

\section{Frequency analysis}\label{sec:fran}
Frequency analysis is a well-known numerical method to identify quasi-periodic
motions and to compute their frequencies. It was introduced by J. Laskar\cite{Laskar90}
and since then it has been used, among other things,
to distinguish regular and chaotic zones, and to estimate the amount of
diffusion through resonant channels \cite{Laskar93,Laskar99}. Here we have used
the refinement of the frequency analysis method that was introduced in
\cite{GomezMondeloFourier1,GomezMondeloFourier2,MondeloPhD}. This refined method is an
iterative procedure that, at each iteration, computes the frequencies with
amplitude greater than a given threshold. This threshold is decreased
at each iteration, up to a prescribed minimum value. Each iteration consists
of three steps: in the first step the frequencies are approximated by
looking at the peaks of the moduli of the discrete Fourier transform of the
signal (using a Hanning filter to reduce leakage), in the second step the
amplitudes corresponding to the previous frequencies are computed by a
collocation method (that  leads to solve a linear system) and, in the third
step, the same collocation method as before is used again but now to 
recompute both frequencies and amplitudes (now by solving a nonlinear
system by means of a Newton method). At each iteration,
the first step is applied to the initial signal minus its current
approximation by a trigonometric polynomial, while in the second and third
steps the collocation is done using the full signal and all the
frequencies and amplitudes found at the moment. See~\cite{GomezMondeloFourier1,GomezMondeloFourier2} for more details.
To detect which asteroids seem to have a quasi-periodic motion and to
obtain their frequencies, we have proceeded as follows. Given the initial
positions and velocities of the asteroid, we have integrated its motion
(in the CRTBP and ERTBP) by means of a Taylor method \cite{JorbaZ05}
with a local threshold of $10^{-16}$, producing a table of values
for the positions of the asteroid each unit of  (adimensional) time.
If the trajectory collides with Sun or Jupiter (here by collision we mean
that the position of the asteroid is ``below the surface'' of Sun or Jupiter,
that is, its distance is lower than their physical radius) we stop the
integration and we discard the asteroid as having quasi-periodic motion.
If the asteroid does not collide, we keep the integration up to a final
time of $2^{20}\approx 10^6$ units, which means that we
have produced a table of $2^{20}$ positions, that is
given to the frequency analysis algorithm. If the output of the algorithm
is a set of frequencies and amplitudes that match the input data, we accept
that the motion is quasi-periodic.

If the motion is quasi-periodic we focus on the
10 frequencies with larger amplitude. In the ERTBP, we expect
to find the frequency 1 (this is the perturbing frequency due
to the eccentricity of Jupiter) plus two other frequencies
that come from the CRTBP, being the remaining frequencies
integer combinations of these three. It is easy to check the
accuracy of these integer combinations as an additional test
of accuracy of the computed frequencies.
As example, Table~\ref{tab:frecus} contains the first 8 frequencies
for (153) Hilda in the CRTBP and the ERTBP. Let us focus first on
the CRTBP case. If we choose $\omega_{1,2}$ as basic frequencies
of the motion, the remaining frequencies have to be integer
combinations of these two. We look for these combinations
to have an estimate of the accuracy of the computations. So,
$|\omega_3-(2\omega_1-\omega_2)|\approx 10^{-15}$,
$|\omega_4-2\omega_1|\approx 10^{-15}$,
$|\omega_5-(3\omega_1-2\omega_2)|\approx 5\times 10^{-15}$,
$|\omega_6-(-\omega_1+2\omega_2)|\approx 10^{-15}$,
$|\omega_7-(\omega_1+\omega_2)|\approx 2\times 10^{-15}$ and
$|\omega_8-(5\omega_1-\omega_2)|\approx 10^{-15}$.
Similarly, in the Elliptic RTBP it is natural to choose, as basic
frequencies, $\omega_{1,2}$ plus $\omega_5$: $\omega_{1,2}$ are close to
the frequencies of (153) Hilda in the circular RTBP and $\omega_5$ is 1,
the frequency of the eccentricity of Jupiter in the Elliptic model.
Here, we have that
$|\omega_3-(2\omega_1-\omega_2)|\approx 3.8\times 10^{-14}$,
$|\omega_4-2\omega_1|\approx 2\times 10^{-15}$,
$|\omega_6-(-\omega_1-1)|\approx 10^{-15}$,
$|\omega_7-(-\omega_1+2\omega_2)|\approx 2.5\times 10^{-14}$ and
$|\omega_8-(3\omega_1-2\omega_2)|\approx 9.2\times 10^{-14}$,
where in the combination to approximate the frequency $\omega_6$ we have
written 1 instead of $\omega_5$. The total computation time
(numerical integration plus frequency analysis) to obtain
Table~\ref{tab:frecus} is of 10 seconds for the RTBP column
and of 15 seconds for the Elliptic RTBP column on an
Intel(R) Xeon(R) Silver 4214R CPU at 2.40GHz.

\begin{table}
   \centering
   {\tt
    \begin{tabular}{|c|c||c|c|}\hline
    \multicolumn{2}{|c||}{Circular RTBP} & \multicolumn{2}{c|}{ Elliptic RTBP}\\ \hline
    $\omega_1$ & 0.505380468452381 & $\omega_1$ & 0.502553031126153\\ \hline
    $\omega_2$ & 0.455882184798376 & $\omega_2$ & 0.451787121788286\\ \hline
    $\omega_3$ & 0.554878752106385 & $\omega_3$ & 0.553318940463982\\ \hline
    $\omega_4$ & 1.010760936904762 & $\omega_4$ & 1.005106062252308\\ \hline
    $\omega_5$ & 0.604377035760396 & $\omega_5$ & 1.000000000000000\\ \hline
    $\omega_6$ & 0.406383901144371 & $\omega_6$ & 0.497446968873846\\ \hline
    $\omega_7$ & 0.961262653250759 & $\omega_7$ & 0.401021212450444\\ \hline
    $\omega_8$ & 2.071020157463529 & $\omega_8$ & 0.604084849801795\\ \hline
    \end{tabular}
    }
    \caption{The first eight frequencies (sorted by amplitude)
    for (153) Hilda in the Circular RTBP (left) and the Elliptic RTBP (right).
    See the text for details.}
    \label{tab:frecus}
\end{table}

\nocite{*}
\bibliography{aipsamp}

\end{document}